\documentclass[12pt]{article}
\usepackage{graphicx,amsmath, amsthm, amssymb, natbib, color,graphicx, hyperref, amsfonts, sectsty, caption, subcaption,
mathrsfs,bbm,float}
\usepackage{psfrag,epsf}
 \usepackage[dvipsnames]{xcolor}
\usepackage{enumitem,tcolorbox,soul}

\def\spacingset#1{\renewcommand{\baselinestretch}%
{#1}\small\normalsize} \spacingset{1}

\usepackage{marquis_general}

\usepackage{tikz}
\usetikzlibrary{arrows}

\usetikzlibrary{fit}

\usepackage{sectsty}
\subsubsectionfont{\underline}

\usepackage{SSL-dense-glm-notations}

\begin{document}

\title{Surrogate Assisted Semi-supervised Inference for High Dimensional Risk Prediction}
\author{Jue Hou, Zijian Guo and Tianxi Cai}%

\maketitle

\begin{abstract}
Risk modeling with EHR data is challenging due to a lack of direct observations on the disease outcome,
and the high dimensionality of the candidate predictors
. In this paper, we develop a surrogate assisted semi-supervised-learning (SAS) approach to risk modeling with high dimensional predictors, leveraging a large
unlabeled  data on candidate predictors and surrogates of outcome, as well as a small labeled data
with annotated outcomes. The SAS procedure borrows information from surrogates along with candidate predictors to impute the unobserved outcomes via a sparse working imputation model with moment conditions to achieve robustness against mis-specification in the imputation model and a one-step bias correction to enable interval estimation for the predicted risk. We demonstrate that the SAS procedure provides valid inference for the predicted risk derived from a high dimensional working model,
even when the underlying risk prediction model is dense and the risk model is mis-specified. We present an extensive simulation study to demonstrate the superiority of our SSL approach compared to existing supervised methods. We apply the method to derive genetic risk prediction of type-2 diabetes mellitus using a EHR biobank cohort.
\end{abstract}

\clearpage
\newpage
\spacingset{1.45} 

\section{Introduction}

Precise risk prediction is vitally important for successful clinical care. High risk patients can be assigned to more intensive monitoring or intervention to improve outcome. Traditionally, risk prediction models are developed based on cohort studies or registry data. Population-based disease registries, while remain a critical source for epidemiological studies, collect information on a relatively small set of pre-specified variables and hence may limit researchers' ability to develop comprehensive risk prediction models \citep{warren2015challenges}. Most clinical care is delivered in healthcare systems \citep{thompson2015linking}, and electronic health records (EHR) embedded in healthcare systems accrue rich clinical data in broad patient populations. EHR systems centralize the data collected during routine patient care including structured elements such as codes for International Classification of Diseases  (ICD), medication prescriptions, and medical procedures, as well as free-text narrative documents such as physician notes and pathology reports that can be processed through natural language processing (NLP) for analysis. EHR data is also often linked with biobanks which provide additional rich molecular information to assist in developing comprehensive risk prediction models for a broad patient population.

Risk modeling with EHR data, however, is challenging due to several reasons. First, precise information on clinical outcome of interest, $Y$, is often embedded in free-text notes and requires manual efforts to extract accurately. Readily available surrogates of $Y$, $\bS$, such as the diagnostic codes or mentions of the outcome are at best good approximations to the true outcome $Y$. For example, using EHR data from Mass General Brigham (MGB), we found that the PPV was only 0.48 and 0.19 for having at least 1 ICD code of T2DM and for having at least 1 NLP mention of T2DM, respectively. Directly using these EHR proxies as true disease status to derive risk models may lead to substantial biases. On the other hand, extracting precise disease status requires manual chart review which is not feasible at a large scale.  It is thus of great interest to develop risk prediction models under a semi-supervised learning (SSL) framework using both a large unlabeled dataset of size $N$ containing information on predictors $\bX$ along with surrogates $\bS$ and a small labeled dataset of size $n$ with additional observations on $Y$ curated via chart review.

Additional challenges arise from the high dimensionality of the predictor vector $\bX$, and the potential model mis-specifications. Although much progress has been made in high dimensional regression in recent years, there is a paucity of literature on high dimensional inference under the SSL setting. Precise estimation of the high dimensional risk model is even more challenging if the risk model is not  sparse. Allowing the risk model to be dense is particularly important when $\bX$ includes genomic markers since a large number of genetic markers appear to contribute to the risk of complex traits \citep{frazer2009human}. For example, \cite{vujkovic2020discovery} recently identified 558 genetic variants as significantly associated with T2DM risk. An additional challenge arises when the fitted risk model is mis-specified, which occurs frequently in practice especially in the high dimensional setting. Model mis-specifications can also lead to the fitted model of $Y \mid \bX$ to be dense. There are limited methods currently available to make inference about high dimensional risk prediction models in the SSL setting especially under a possibly mis-specified dense model.  In this paper, we fill in the gap by proposing an efficient surrogate assisted SSL (SAS) prediction  procedure that  leverages the fully observed surrogates $\bS$ to make inference about a high dimensional risk model under such settings.

Under the supervised setting where both $Y$ and $\bX$ are fully observed, much progress has been made in recent years in the area of high dimensional inference. High dimensional regression methods have been developed for commonly used generalized linear models (GLM) under sparsity assumptions on the regression parameters \citep{vdGeerBuhlmann09,NRWY2010TR,HuangZhang12}. Recently, \cite{ZhuBradic18ejs} studied the inference of linear combination of coefficients under dense linear model and sparse precision matrix. Inference procedures have also been developed for both sparse \citep{ZhangZhang14,JavanmardMontanari14,vdGeerEtal14} and dense combinations of the regression parameters \citep{CaiCaiGuo19,ZhuBradic18}.
High-dimensional inference under the logistic regression model has also been studied recently \citep{vdGeerEtal14,MaCaiLi20,guo2020inference}.

Under the SSL setting with $n\ll N$, however, there is a paucity of literature on high dimensional inference. Although the SSL can be viewed as a missing data problem, it differs from the standard missing data setting in a critical way. Under the SSL setting, the missing probability tends to 1, which would violate a key assumption required in the missing data literature  \citep[e.g.]{BangRobins05,SmuclerEtal19,ChakraborttyEtal19}.
Existing work on SSL with high-dimensional covariates largely focuses on the post-estimation inference on the global parameters under sparse linear models with examples including SSL estimation of population mean \citep{ZhangEtal2019,ZhangBradic19}, the explained variance \citep{CaiGuo2018}, and the average treatment effect  \citep{ChengEtal2018,KallusMao2020arxiv}. To the best of our knowledge,
our SAS procedure is the first to conduct the semi-supervised inference of the high-dimensional coefficient and the individual prediction in the high-dimensional dense and possibly mis-specified risk prediction model.

Our proposed estimation and inference procedures are as follows.
For estimation, we first use the labelled data to fit a regularized imputation model with surrogates and high-dimensional covariates; then we impute the missing outcomes for the unlabeled
data and fit the risk model using the imputed outcome and high-dimensional predictors. 
For inference, we devise a novel bias correction method, which corrects the bias due to the regularization for both imputation and estimation.
For our proposed methods, we allow the fitted risk model for $Y \mid \bX$ to be both mis-specified and potentially dense but only require sparsity on the fitted imputation model of $Y \mid \bS, \bX$. The sparsity assumption on the imputation model is less stringent since we anticipate that most information on $Y$ can be well captured by the low dimensional $\bS$ while the fitted model of $Y \mid \bX$ might be dense especially under possible model mis-specifications.

The remainder of the paper is organized as follows. We introduce our population parameters and model assumptions in Section \ref{section:setup}. In Section \ref{section:method}, we propose the SAS estimation method along with its associated inference procedures. In Section \ref{section:theory}, we state the theoretical guarantees of the SAS procedures, whose proofs are provided in the Supplementary Materials. We also remark on the sparsity relaxation and the efficiency gain of the SSL. In Section \ref{section:simulation},
we present simulation results highlighting finite sample performance of the SAS estimators and comparisons to existing methods.
In Section \ref{section:data}, we apply the proposed method to derive individual risk prediction for T2DM using EHR data from MGB. 

\section{Settings and Notations}\label{section:setup}
For the $i$-th observation, $Y_i\in \R$ denotes the outcome variable, $S_i\in \R^q$ denotes the surrogates for $Y_i$ and $\bX_i\in \R^{p+1}$ denotes the high-dimensional covariates with the first element being the intercept.
Under the SSL setting, we observe $n$ independent and identically distributed (i.i.d.) labeled observations, $\Lscr = \{ (Y_i, \bX_i\trans,\bS_i\trans)\trans, i = 1, ..., n\}$ and $(N-n)$ i.i.d unlabeled observations, $\Uscr = \{\bW_i=(\bX_i\trans,\bS_i\trans)\trans, i = n+1, ..., N\}$.
We assume that the labeled subjects are randomly sampled by design and the proportion of labelled sample is $n/N = \rho \in (0,1)$ with $\rho \to 0$ as $n\to \infty$.
We focus on the high-dimensional setting where dimensions $p$ and $q$ grow with $n$ and allow $p+q$ to be larger than $n$.
Motivated by our application, our main focus is on the setting $N$ much larger than $p$, but our approach can be extended to $N\leq p$ under specific conditions.

To predict $Y_i$ with $\bX_i$, we consider a possibly mis-specified {\em working} conditional mean model with a known monotone and smooth link function $g$,
\begin{equation}\label{model:Y_X}
  E(Y_i \mid \bX_i) = g(\bbeta\trans \bX_i) .
\end{equation}
Our procedure generally allows for a wide range of link functions and detailed requirements on $g(\cdot)$ and its anti-derivative $G$ are given in Section \ref{section:theory}. In our motivating example, $Y$ is a binary indicator of T2DM status and $g(x) =1/(1+e^{-x})$ with  $G(x) = \log(1+e^x)$.
Our goal is to accurately estimate the high-dimensional parameter $\bbeta_0,$ defined as the solution of the estimation equation
\begin{equation}\label{model:Y_X-moment}
\E[ \bX_i \{ Y_i -  g(\bbeta\subo \trans \bX_i)\}] = 0.
\end{equation}
We shall further construct confidence intervals for
$g(\bbeta\subo\trans \bx\subnew)$ with any $\bx\subnew \in \R^{p+1}$.
The predicted outcome $g(\bbeta\subo\trans \bx\subnew)$ is the condition mean of $Y$ given $\bX_i=\bx\subnew$ when \eqref{model:Y_X} holds
and can be interpreted as the best generalized linear prediction with link $g$ even when \eqref{model:Y_X} fails to hold. To enable estimation of $\bbeta\subo$ under possible model mis-specification, we define a pseudo log-likelihood (PL) function
 \begin{equation}\label{def:PL}
 \ell(y, x) = y x - G(x)
 \end{equation}
 such that $\partial \ell(y,\bbeta\trans\bX_i) /\partial \bbeta = \bX_i\trans\{y - g(\bbeta\trans\bX_i)\}$ corresponds to the moment condition \eqref{model:Y_X-moment}.
We make no assumption on the sparsity of $\bbeta\subo$ and hence it is not feasible to perform valid supervised learning for $\bbeta\subo$ when $\sbeta = \|\bbeta\subo\|\subo > n$.

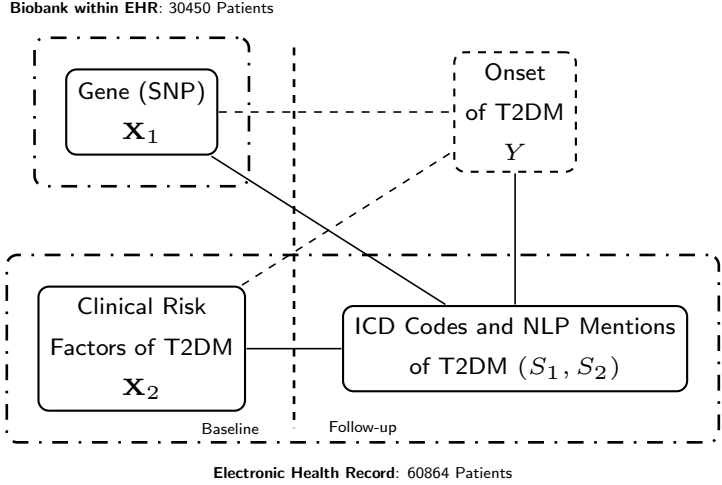
\begin{figure}
\centering
\begin{tikzpicture}[
  font=\sffamily\tiny,
  every matrix/.style={ampersand replacement=\&,column sep=0.5cm,row sep=1.5cm},
  observed/.style={draw,thick,rounded corners, fill=white,inner sep=.1cm,scale=1.5},
  unobserved/.style={draw,thick,dashed,rounded corners, fill=white,inner sep=.1cm,scale=1.5},
  dots/.style={gray,scale=2},
  to/.style={->,>=stealth',shorten >=1pt,semithick,font=\sffamily\footnotesize},
  link/.style={semithick,font=\sffamily\footnotesize},
  every node/.style={align=center}]
  \matrix{
     \node[observed] (X) {Gene (SNP) \\ $\bX_1$ }; \& \node[] (E1) {};
      \&   \node[unobserved] (Y) {Onset  \\of T2DM \\ $Y$}; \\
   \node[observed] (X2) {Clinical Risk\\ Factors of T2DM\\ $\bX_2$}; \& \node[] (E2) {};
   \&  \node[observed] (M) {ICD Codes and   NLP Mentions \\ of T2DM  $(S_1, S_2)$};\\
  };
  	\draw[link] (X) --  (M);
  	\draw[link] (X2) --  (M);
	\draw[link,dashed] (X) -- (Y);
	\draw[link,dashed] (X2) -- (Y);
  	\draw[link] (M) --  (Y);

  \node at (E1.north) [above, inner sep=3.5mm](E0){};
\node at (E2.south) [below, inner sep=4.5mm](E3){};
	\draw[dashed,line width=1pt] (E0.north) -- (E3.south);
\node at (E3.south west) [left, inner sep=0mm](B){Baseline};
\node at (E3.south east) [right, inner sep=0mm](F){Follow-up};

  \tikzset{grey dotted/.style={draw, line width=1pt,
                               dash pattern=on 1pt off 4pt on 6pt off 4pt,
                                inner sep=4mm, rectangle, rounded corners}};
  \tikzset{crimson dotted/.style={draw, line width=1pt,
                               dash pattern=on 1pt off 4pt on 6pt off 4pt,
                                inner sep=4mm, rectangle, rounded corners}};

  \node (first dotted box) [grey dotted,
                            fit = (X)] {};
  \node (second dotted box) [crimson dotted,
                            fit =(X2) (M)] {};

  \node at (first dotted box.north) [above, inner sep=3mm]
  {\textbf{Biobank within EHR}: 30450 Patients};
  \node at (second dotted box.south) [below, inner sep=3mm]
  {\textbf{Electronic Health Record}: 60864 Patients};
\end{tikzpicture}
\caption{A dense prediction model (graph with dashed lines) can be compress to a sparse imputation model
(through graph with solid lines)
when the effect of most baseline covariates are reflected in a few variables
in the EHR monitoring the development of the event of interest.
}
\label{fig:dag}
\end{figure}

We shall derive an efficient SSL estimate for $\bbeta\subo$ by leveraging $\Uscr$. To this end, we fit a {\em working} imputation model
\begin{equation}\label{model:Y_W}
  E(Y_i \mid \bW_i) = g(\bgamma\trans \bW_i)
\end{equation}
with the population parameter $\bgamma\subo$ defined as
\begin{equation}\label{eq:gamma-score-pop}
   \bgamma\subo: \; \E[\bW_i\{Y_i-g(\bgamma\subo\trans \bW_i)\}] = 0.
\end{equation}
The definition of $\bgamma$ guarantees
\begin{equation}\label{def:bar-gamma}
  \E[\bX_i\{Y_i-g(\bgamma\subo\trans \bW_i)\}] = 0.
\end{equation}
and hence if we impute $Y_i$ as $\Ybar_i = g(\bgamma\subo\trans\bW_i)$, we have $\E[\bX_i\{\Ybar_i - g(\bbeta\subo\trans\bX_i)\}] = 0$ regardless the adequacy of the imputation model \eqref{model:Y_W}. It is thus feasible to carry out an SSL procedure by first deriving an estimate for $\Ybar_i$ using the labelled data $\Lscr$ and then regressing the estimated $\Ybar_i$  against $\bX_i$ using the whole data $\Lscr\cup\Uscr$.  Although we do not require $\bbeta\subo$ to be sparse or any of the fitted models to hold, we do assume that $\bgamma\subo$ defined in \eqref{eq:gamma-score-pop} to be sparse. When the surrogates $\bS$ are strongly predictive for the outcome, the sparsity assumption on $\bgamma\subo$ is reasonable since the majority of the information in $Y$ can be captured in $\bS$.

\noindent {\bf Notations.}
We focus on the setting where $\min\{n,p+q,N\}\rightarrow \infty.$ For convenience, we shall use $n\rightarrow \infty$ in the asymptotic analysis.
For two sequences of random variables $A_n$ and $B_n$, we use $A_n = O_p(B_n)$ and $A_n=o_p(B_n)$
to denote $\lim_{c \to \infty}\lim_{n\to \infty} \P(|A| \ge c |B|) = 0$
and $\lim_{c \to 0}\lim_{n\to \infty} \P(|A| \ge c |B|) = 0$, respectively. For two positive sequences $a_n$ and $b_n$,  $a_n=O(b_n)$ or $b_n\gtrsim a_n$ means that $\exists C > 0$ such that $a_n \leq C b_n$ for all $n$;
$a_n \asymp b_n$ if $a_n=O(b_n)$ and $b_n=O(a_n)$, and $a_n \ll b_n$ or $a_n=o(b_n)$ if $\limsup_{n\rightarrow\infty} {a_n}/{b_n}=0$. We use $Z_n \stackrel{\mathscr{L}}{\rightarrow} N(0,1)$
to denote the sequence of random variables $Z_n$ converges in distribution
to a standard normal random variable.

\section{Methodology}\label{section:method}
\subsection{SAS Estimation of $\bbeta\subo$}
The SAS estimation procedure for $\bbeta\subo$ consists of two key steps: (i) fitting the imputation model to $\Lscr$ to obtain estimate $\bgammahat$ for $\bgamma\subo$ defined in \eqref{eq:gamma-score-pop}; and (ii) estimating $\bbeta\subo$ in \eqref{model:Y_X-moment} by fitting imputed outcome $\Yhat_i = g(\bgammahat\trans\bW_i)$ against $\bX_i$ to $\Uscr$. 

In Step (i), we estimate $\bgamma\subo$ by the $L_1$ regularized PL estimator $\bgammahat,$ defined as
\begin{equation}\label{def:gamma-lasso}
  \hat{\bgamma} = \argmin_{\bgamma \in \R^{p+q+1}}
 \implik (\bgamma)
  + \lamg \|\bgamma_{-1}\|_1 \quad \text{with}
\quad
  \lamg \asymp \sqrt{\log(p+q)/n},
\end{equation}
where $\ba_{-1}$ denotes the sub-vector of all the coefficients except for the intercept and \begin{equation}\label{def:gamma-dev}
  \implik (\bgamma) = \frac{1}{n}\sum_{i=1}^{n} \ell(Y_i, \bgamma\trans \bW_i) \quad \text{with}\quad \ell(y,x) \; \text{defined in }\; \eqref{def:PL}.
\end{equation}
The imputation loss \eqref{def:gamma-dev} corresponds to the negative log-likelihood when $Y$ is binary and the imputation model holds with $g$ being anti-logit.
With $\hat{\bgamma}$, we impute the unobserved outcomes for subjects in $\Uscr$ as
$\Yhat_i = g(\hat{\bgamma}\trans \bW_i)$, for $n+1\leq i \leq N$.

In Step (ii), we estimate $\bbeta\subo$ by $\bbetahat =\bbetahat(\bgammahat),$ defined as,
\begin{equation}\label{def:bhat}
\bbetahat(\bgammahat) =  \argmin_{\bbeta \in \R^{p+1}} \ell^\dag(\bbeta;\bgammahat) + \lamb\|\bbeta_{-1}\|_1 \quad \text{with}
\quad \lamb\asymp \sqrt{\log(p)/N},
\end{equation}
where $\ell^\dag(\bbeta;\hat{\bgamma})$ is the imputed PL:
\begin{align}
  \ell^\dag(\bbeta;\hat{\bgamma})
  =  \frac{1}{N}\sum_{i>n}\ell(\Yhat_i, \bbeta\trans \bX_i) + \frac{1}{N}\sum_{i=1}^{n}\ell(Y_i, \bbeta\trans \bX_i) \quad \text{with}\quad \ell(y,x) \; \text{defined in }\; \eqref{def:PL}.\label{def:beta-dev}
\end{align}

We denote the complete data PL of the full data as
\begin{equation}\label{def:beta-dev-comp}
  \ell\subPL(\bbeta) =  \frac{1}{N}\sum_{i=1}^{N}\ell(Y_i, \bbeta\trans \bX_i) .
\end{equation}
and define the gradients of the various losses \eqref{def:gamma-dev}-\eqref{def:beta-dev-comp} as
\begin{equation}\label{def:dot-grad}
  \impscore(\bgamma) = \nabla \implik(\bgamma), \;
   \dell\subPL(\bbeta) = \nabla \ell\subPL(\bbeta), \;
   \dell^\dag(\bbeta;\bgamma) = \frac{\partial }{\partial \bbeta} \ell^\dag(\bbeta;\bgamma).
\end{equation}

\subsection{SAS Inference for Individual Prediction}\label{section:inf}

Since $g(\cdot)$ is specified, the inference on $g(\bx\subnew\trans \bbeta)$ immediately follows from the inference on $\bx\subnew\trans \bbeta$.
We shall consider the inference on standardized linear prediction $\bx\substd\trans \bbeta$ with
the standardized covariates
\begin{equation*}
  \bx\substd = \bx\subnew/\|\bx\subnew\|_2
\end{equation*}
and then scale the confidence interval back.
This way, the scaling with $\|\bx\subnew\|_2$ is made explicit in the expression of the confidence interval.

The estimation error of $\hat{\bbeta}$ can be decomposed into two components corresponding to the respective errors associated with
\eqref{def:gamma-lasso} and \eqref{def:bhat}. Specifically, we  write
\begin{equation}\label{eq:bias+var}
  \hat{\bbeta} - \bbeta\subo = \{\bar{\bbeta}(\hat{\bgamma})-\bbeta\subo\} +  \{\hat{\bbeta}-\bar{\bbeta}(\hat{\bgamma})\},
\end{equation}
where $\bar{\bbeta}(\hat{\bgamma})$ is defined as the minimizer of the expected imputed loss conditionally on the labeled data, that is,
\begin{equation}\label{def:bar-beta}
  \bar{\bbeta}(\hat{\bgamma})
  = \argmin_{\bbeta \in \R^{p+1}}\E[\ell^\dag(\bbeta;\hat{\bgamma}) \mid \Lscr].
\end{equation}
The term $\bar{\bbeta}(\hat{\bgamma})-\bbeta\subo$ denotes
the error from the imputation model in \eqref{def:gamma-lasso} while
the term $\hat{\bbeta} - \bar{\bbeta}(\hat{\bgamma})$ denotes
the error from the prediction model in \eqref{def:bhat}
given the imputation model parameter $\hat{\bgamma}$.
As $\ell_1$ penalization is involved in both steps, we shall correct the regularization bias
from the two sources.
Following from
the typical one-step debiasing LASSO \citep{ZhangZhang14}, the bias $\hat{\bbeta} - \bar{\bbeta}(\hat{\bgamma})$ is estimated by
$\hPrec \dell^\dag(\hat{\bbeta};\hat{\bgamma})$
where $\hPrec$ is an estimator of $[\E\{g'(\bbeta\subo\trans \bX_i)\bX_i\bX_i\trans\}]^{-1}$,
the inverse Hessian of $\ell^\dag(\cdot;\hat{\bgamma})$ at $\bbeta=\bbeta\subo$.

The bias correction for $\bar{\bbeta}(\hat{\bgamma})-\bbeta\subo$
requires some innovation since we need to conduct the bias correction for a nonlinear functional $\bar{\bbeta}(\cdot)$ of LASSO estimator $\hat{\bgamma},$ which has not been studied in the literature.
  We identify $\bar{\bbeta}(\hat{\bgamma})$ and $\bbeta\subo$
by the first order moment conditions,
\begin{align}
\bar{\bbeta}(\hat{\bgamma})&: \E_{i > n}[\bX_i\{g(\bar{\bbeta}(\hat{\bgamma})\trans \bX_i) - g(\hat{\bgamma}\trans \bW_i)\}\mid \Lscr] \approx 0, \notag \\
\bbeta\subo &: \E[\bX_i\{g(\bbeta\subo\trans \bX_i) - Y_i\}] = \E[\bX_i\{g(\bbeta\subo\trans \bX_i) - g(\bgamma\subo\trans \bW_i)\}] = 0. \label{eq:1st-beta-bar}
\end{align}
Here $\E_{i > n}[\cdot \mid \Lscr]$ denotes the conditional expectation of a single copy of the unlabeled data
given the labelled data.
By equating the two estimating equations in \eqref{eq:1st-beta-bar}, we apply the first order approximation and approximate the difference $\bar{\bbeta}(\hat{\bgamma})-\bbeta\subo$
by
\begin{equation}
\begin{aligned}
   \bar{\bbeta}(\hat{\bgamma}) - \bbeta\subo
 \approx & - [\E\{g'(\bbeta\subo\trans \bX_i)\bX_i\bX_i\trans\}]^{-1} \E_{i > n}\left[\bX_i\{g(\bgamma\subo\trans \bW_i)-g(\hat{\bgamma}\trans \bW_i)\}\mid \Lscr\right] 
 \end{aligned}
\end{equation}
Together with the bias correction for $\bar{\bbeta}(\hat{\bgamma})-\bbeta\subo$, this motivates the debiasing procedure
$$
  \hat{\bbeta}- \frac{1-\rho}{n}\sum_{i=1}^{n}\hPrec \bX_i\left\{ g(\hat{\bgamma}\trans \bW_i)
-Y_i\right\}-\hPrec \dell^\dag(\hat{\bbeta};\hat{\bgamma}).
$$
The $1-\rho$ factor, which tends to one when $n$ much smaller than $N$, comes from the proportion of unlabeled data whose missing outcome are imputed.

For theoretical considerations,
we devise a cross-fitting scheme in our debiasing process.
We split the labelled and unlabeled data
into $\Kn$   folds of approximately equal size, respectively. The number of folds does not grow with dimension (e.g. $\Kn=10$).
We denote the indices sets for each fold of the labelled data $\Lscr$
as $\mathcal{I}_1,\dots, \mathcal{I}_{\Kn}$,
and those of the unlabeled data $\Uscr$ as $\mathcal{J}_1,\dots, \mathcal{J}_{\Kn}$.
We denote the respective sizes of each fold in the labelled data and full data as $n_k = |\Ical_k|$ and $N_k = n_k + |\Jcal_k |$, where $|\Ical|$ denotes the carnality of $\Ical$.
Define $\Ical^c_k = \{1,\dots,n\}\setminus \Ical_k$ and $\Jcal^c_k = \{n+1,\dots,N\}\setminus \Jcal_k.$
For each labelled fold $k$, we fit the imputation model with out-of-fold
labelled samples:
\begin{equation}\label{def:gammak}
  \hat{\bgamma}\supk = \argmin_{\bgamma \in \R^{p+q+1}}
 \frac{1}{n-n_k}\sum_{i\in \Ical^c_k}\ell(Y_i, \bgamma\trans\bW_i)
  + \lamg \|\bgamma_{-1}\|_1.
\end{equation}
Using $\hat{\bgamma}\supk$, we fit the prediction model with the out-of-fold data
$\Ical^c_k \cup \Jcal^c_k$:
\begin{equation}\label{def:betak}
  \hat{\bbeta}\supk = \argmin_{\bbeta \in \R^{p+1}} \frac{1}{N-N_k}\left[\sum_{i\in\Jcal_k^c}\ell(g(\hat{\bgamma}\supkt\bW_i), \bbeta\trans \bX_i) + \sum_{i\in\mathcal{I}_k^c}\ell(Y_i, \bbeta\trans \bX_i)\right]
  + \lamb\|\bbeta_{-1}\|_1.
\end{equation}
To estimate the
projection
\begin{equation}\label{def:u0}
  \bu\subo = \E\{g'(\bbeta\subo\trans \bX_i)\bX_i\bX_i\trans\}]^{-1} \bx\substd,
\end{equation}
we propose an $L_1$-penalized estimator
\begin{equation}\label{def:uk}
  \hat{\bu}\supk =  \argmin_{\bu \in \R^p} \frac{1}{N-N_k} \sum_{k'\neq k}\sum_{i \in \Ical_{k'}\cup\Jcal_{k'}}\left[ \frac{1}{2}g'\left(\hat\bbeta\supkkt\bX_i\right)(\bX_i\trans\bu)^2 -\bu\trans\bx\substd + \lamu\|\bu\|_1\right],
\end{equation}
where $\hat{\bbeta}\supkk$ is trained with samples out of folds $k$ and $k'$,
\begin{gather}
  \hat{\bbeta}\supkk = \argmin_{\bbeta \in \R^{p+1}} \frac{\sum_{i\in(\Jcal_k\cup\Jcal_{k'})^c}\ell\left(g\left(\hat{\bgamma}\supkkt\bW_i\right), \bbeta\trans \bX_i\right) + \sum_{i\in(\Ical_k\cup\Ical_{k'})^c}\ell(Y_i, \bbeta\trans \bX_i)}{N-N_k-N_{k'}}
  + \lamb\|\bbeta_{-1}\|_1,\label{def:est-kk} \\
\text{with}\quad
\hat{\bgamma}\supkk = \argmin_{\bgamma \in \R^{p+q+1}}
 \frac{\sum_{i\in \Ical^c_k \cap \Ical^c_{k'}}\ell(Y_i, \bgamma\trans\bW_i)}{n-n_k-n_{k'}}
  + \lamg \|\bgamma_{-1}\|_1. \notag
\end{gather}
The estimators in \eqref{def:est-kk}  take similar forms as those in \eqref{def:gammak} and \eqref{def:betak}
except that their training samples exclude two folds of data $\Ical_k \cup \Jcal_{k}$ and $\Ical_{k'} \cup \Jcal_{k'}$.
In the summand of \eqref{def:uk},
the data $(Y_i, \bX_i, \bS_i)$ in fold $k'$ $\Ical_{k'} \cup \Jcal_{k'}$ is independent of  $\hat{\bbeta}\supkk$ trained without folds $k$ and $k'$.
The estimation of $\bu$ requires an estimator of $\bbeta$ and both estimators are subsequently used for the debiasing step.
Using the same set of data multiple times for $\hat{\bbeta}$, $\hat{\bu}$,
 debiasing and variance estimation
may induce over-fitting bias,
so we implemented the cross-fitting scheme to reduce the over-fitting bias.
As a remark, cross-fitting might not be necessary for theory with additional assumptions
and/or empirical process techniques.

We obtain the cross-fitted debiased estimator for $\bx\substd\trans \bbeta$ as $\widehat{\bx\substd\trans \bbeta}$, defined as
\begin{equation}\label{xbet-crossfit}
\begin{split}
& \frac{1}{\Kn}\sum_{k=1}^{\Kn} \bx\substd\trans \hat{\bbeta}\supk
-\frac{1}{N} \sum_{k=1}^{\Kn}\sum_{i\in\Jcal_k}\hat{\bu}\supkt\bX_i\{g(\hat{\bbeta}\supkt\bX_i)-g(\hat{\bgamma}\supkt  \bW_i)\}\\
&-\frac{1}{n} \sum_{k=1}^{\Kn} \sum_{i \in \mathcal{I}_k} \hat{\bu}\supkt \bX_i\left\{(1-\rho)\cdot g(\hat{\bgamma}\supkt  \bW_i) + \rho\cdot g(\hat{\bbeta}\supkt\bX_i)
-Y_i\right\} .
\end{split}
\end{equation}
The second term is used to correct the bias $\bar{\bbeta}(\hat{\bgamma})-\bbeta_0$ and the third term is used to correct the bias $\hat{\bbeta} - \bar{\bbeta}(\hat{\bgamma})$.
The corresponding variance estimator is
\begin{align}
  \hVSAS
  =&  \frac{1}{n}\sum_{k=1}^{\Kn}\sum_{i \in \Ical_k}
  (\hat{\bu}\supkt \bX_i)^2\left\{(1-\rho) \cdot g(\hat{\bgamma}\supkt  \bW_i)
  +\rho \cdot g(\hat{\bbeta}\supkt \bX_i) -Y_i\right\}^2 \notag \\
  &+ \frac{\rho^2}{n} \sum_{k=1}^{\Kn}\sum_{i\in \Jcal_k} (\hat{\bu}\supkt\bX_i)^2 \left\{g(\hat{\bbeta}\supkt\bX_i) - g(\hat{\bgamma}\trans \bW_i)\right\}^2\label{def:var-cf}
\end{align}
Through the link $g$ and the scaling factor $\|\bx\subnew\|_2$,
we estimate $g(\bx\subnew\trans \bbeta_0)$ by $g\left(\|\bx\subnew\|_2\widehat{\bx\substd\trans \bbeta} \right)$ and construct the $(1-\alpha) \times 100\%$ confidence interval for
$g(\bx\subnew\trans \bbeta_0)$ as
\begin{equation}\label{def:CI}
  \left[g\left\{ \|\bx\subnew\|_2\left(\widehat{\bx\substd\trans \bbeta}  - \mathcal{Z}_{\alpha/2} \sqrt{\hVSAS /n}\right\}\right), \,
  g\left\{ \|\bx\subnew\|_2 \left(\widehat{\bx\substd\trans \bbeta}  + \mathcal{Z}_{\alpha/2} \sqrt{\hVSAS /n} \right)\right\}\right] ,
\end{equation}
where $Z_{\alpha/2}$ is the $1-\alpha/2$ quantile of the standard normal distribution.

\section{Theory}\label{section:theory}

We introduce assumptions required for both estimation and inference in Section \ref{section:assume}.
We state our
theories for estimation
and inference, respectively in Sections \ref{section:thm-est} and \ref{section:thm-inf}.

\subsection{Assumptions}\label{section:assume}

We assume the complete data consist of i.i.d. copies of $(Y_i, \bX_i, \bS_i)$,
for $i=1,\dots, N$.
For our focused SSL settings, only the first $n$ outcome labels $Y_i,\dots,Y_n$ are observed. Under the i.i.d assumption, our SSL setting
is equivalent to the missing completely at random (MCAR) assumption.
The sparsities of $\bgamma\subo$, $\bbeta\subo$ and $\bu\subo$ are denoted as
\begin{equation*}
    \sgamma = \|\bgamma\subo\|_0,\; \sbeta= \|\bbeta\subo\|_0,\; \su= \|\bu\subo\|_0.
\end{equation*}
We focus on the setting with $n,p+q,N \rightarrow \infty$ with $n$ being allowed to be smaller than $p+q$. We allow that $\sgamma,\sbeta$ and $\su$ grow with $n,p+q,N$ and
satisfy $\sgamma \ll n$ and $\sbeta+\su \ll N$.
To achieve the sharper dimension conditions,
we consider the sub-Gaussian design as in \cite{Portnoy1984,Portnoy1985,NRWY2010TR}.
We denote the sub-Gaussian norm for random variables and
random vectors both as $\|\cdot\|_{\psi_2}$.
The detailed definition is given in Appendix \ref{asection:detail}. 
\begin{assumption}
\label{assume:model}
For constants $\nu_1$, $\nu_2$ and $M$ independent of $n,p$ and $N$,
\begin{enumerate}[label = \alph*), ref = \ref{assume:model}\alph*]
    \item \label{assume:Y-tail} the residuals $Y_i - g(\bgamma\subo\trans \bW_i)$ and $Y_i - g(\bbeta\subo\trans \bX_i)$ are sub-Gaussian
    random variables with sub-Gaussian norm bounded by $\left\|Y_i - g(\bgamma\subo\trans \bW_i)\right\|_{\psi_2} \le \nu_1$ and $\left\|Y_i - g(\bbeta\subo\trans \bX_i)\right\|_{\psi_2} \le \nu_2$;  
     \item \label{assume:link} The link function $g$ satisfies the monotonicity and
         smoothness conditions: $\inf_{x\in\R} g'(x) \ge 0$,  $\sup_{x\in\R} g'(x) < M$
         and $\sup_{x\in\R} g^{\prime\prime}(x) < M$.
\end{enumerate}
\end{assumption}

Under our motivating example with a binary $Y_i$ and
$
g(x)=e^{x}/(1+e^x)
$, \ref{assume:Y-tail} and \ref{assume:link} are satisfied. The condition is also satisfied for
the probit link function and the identity link function.
Condition \ref{assume:Y-tail} is universal for high-dimensional regression.
Admittedly, Lipschitz requirement in \ref{assume:link} rules out some GLM links
with unbounded derivatives
like the exponential link, but we may substitute the condition
by assuming a bounded $\|\bX_i\|_\infty$. 

\begin{assumption}
\label{assume:X}
For constants $\sigma\submax^2$ and $\sigma\submin^2$ independent of $n,p,N$,
\begin{enumerate}[label = \alph*), ref = \ref{assume:X}\alph*]
  \item \label{assume:X-subG} $\bW_i$ is a sub-Gaussian vector
  with sub-Gaussian norm
  $\|\bW_i\|_{\psi_2}\leq \sigma\submax/\sqrt{2}$;
   \item \label{assume:sigmin}
    The weak overlapping condition at the population parameter $\bbeta\subo$ and $\bgamma\subo$,
    \begin{enumerate}[label = (\roman*), ref =\ref{assume:sigmin}-\roman* ]
        \item \label{assume:sigmin-X} $\inf_{\|\bv\|_2 = 1} \bv\trans \E([g'(\bbeta\subo\trans \bX_i) \wedge 1]\bX_i\bX_i\trans) \bv \ge \sigma^2\submin$,
        \item \label{assume:sigmin-W} $\inf_{\|\bv\|_2 = 1} \bv\trans \E[\{g'(\bgamma\subo\trans \bW_i)\wedge 1\}\bW_i\bW_i\trans] \bv \ge \sigma^2\submin$;
    \end{enumerate}
    \item \label{assume:nondeg_var}
      The non-degeneracy of average residual variance:
      $$\inf_{\|\bv\|_2 = 1} \E [ \{Y_i - (1-\rho)\cdot g(\bgamma\subo\trans \bW_i)-\rho\cdot g(\bbeta\subo\trans\bX_i)\}^2 (\bX_i\trans \bv)^2]  \ge \sigma^2\submin.$$
\end{enumerate}
\end{assumption}

Assumption \ref{assume:X-subG} is typical for high-dimensional
regression \citep{NRWY2010TR}, which also implies the
bounded maximal eigenvalue of the second moment
$$\sup_{\|\bv\|_2 = 1} \bv\trans \E[\bW_i\bW_i\trans] \bv \le \sigma^2\submax.$$
Notably, we do not require two common conditions under high-dimensional generalized linear models \citep{HuangZhang12,vdGeerEtal14}:
1) the upper bound on $\sup_{i=1,\dots,N}\|\bX_i\|_\infty$;
2) the lower bound on  $\inf_{i=1,\dots,N} g'(\bbeta\subo\trans\bX_i)$, often
known as the overlapping condition for logistic regression model.
Compared to the overlapping condition under logistic regression
that $g(\bbeta\subo\trans\bX_i)$ and $g(\bgamma\subo\trans \bW_i)$ are bounded away from zero,
our Assumptions \ref{assume:sigmin} and \ref{assume:nondeg_var} are weaker
because they are implied by
the typical minimal eigenvalue condition
$$\inf_{\|\bv\|_2 = 1} \bv\trans \E(\bW_i\bW_i\trans) \bv \ge \sigma^2\submin$$ plus the overlapping condition.

\subsection{Consistency of the SAS Estimation}\label{section:thm-est}

We now state the $L_2$ and $L_1$ convergence rates of our proposed SAS estimator. 
\begin{theorem}[Consistency of SAS estimation]\label{thm:consistency}
Under Assumptions \ref{assume:model}, \ref{assume:X} and with
\begin{equation}
\sgamma=o(n/\log(p+q)), \, \sbeta=o(N/\log(p)), \, \lamb \gtrsim \sqrt{\log(p)/N},
\end{equation}
we have
\begin{gather*}
  \|\hat{\bbeta}-\bbeta\subo\|_2 = O_p\left(\sqrt{\sbeta}\lamb+ (1-\rho)\sqrt{\sgamma \log(p+q)/n}\right), \\
\|\hat{\bbeta}-\bbeta\subo\|_1 = O_p\left(\sbeta\lamb + (1-\rho)^2\sgamma \log(p+q)/(n\lamb)\right).
\end{gather*}
\end{theorem}

A few remarks are in order for Theorem \ref{thm:consistency}. First, the dimension requirement for our SAS estimator achieving $L_2$ consistency significantly weakens the existing dimension requirement in the supervised setting \citep{NRWY2010TR,HuangZhang12,BuhlmannVDGeer2011,BickelRT09}  With $\lamb \asymp \sqrt{\log(p)/N},$
Theorem \ref{thm:consistency} implies the $L_2$ consistency of $\hat{\bbeta}$ under the dimension condition,
\begin{equation}
 (1-\rho)^2 \sgamma\log(p+q)/n +  \sbeta\log(p)/N = o(1).  \label{assume:dim-est}
  \end{equation}
  For the setting $N \gg n,$
  our requirement on the sparsity of $\beta$, $\sbeta = o(N/\log(p))$
  is significantly weaker than $\sbeta = o(n/\log(p)),$ which is known as the fundamental sparsity limit to identify the high-dimensional regression vector in the supervised setting.
  Theorem \ref{thm:consistency} indicates that with assistance from observed $\bS \in \Uscr$, the SAS procedure allows $\sbeta > n$ provided that $N$ is sufficiently large and the imputation model is sparse. This distinguishes our result
from most estimation results in high-dimensional supervised settings.

Second, we briefly discuss the $L_1$ consistency. If the $L_1$ consistency is of interest, the penalty levels are chosen as
\begin{equation}\label{def:lamb-opt}
\lamb \asymp \max\left\{\sqrt{\log(p)/N}, \sqrt{\sgamma/\sbeta} \lamg\right\},
\end{equation}
which produces the $L_1$ estimation rate from Theorem \ref{thm:consistency}
$$
\|\hat{\bbeta}-\bbeta\subo\|_1 = O_p\left(\sbeta\sqrt{\log(p)/N}+ \sqrt{\sgamma\sbeta\log(p)/n}\right).
$$
Compared to the condition for $L_1$ consistency under supervised learning,
$\sbeta= o\left(\sqrt{n/\log(p)}\right)$,
the condition from SAS estimation $\sbeta = o\left((n/\sgamma+N)/\log(p)\right)$
allows a denser $\bbeta\subo$ in the setting with a very sparse $\bgamma\subo$ and a large unlabeled data.
On the other hand, the $L_2$ estimation rate in Theorem \ref{thm:consistency} remains the same if $$\sqrt{\log(p)/N}\lesssim \lamb\lesssim \max\left\{\sqrt{\log(p)/N}, \sqrt{\sgamma/\sbeta} \lamg\right\}.$$ We shall point out that our subsequent theory on the SAS inference procedure is based the $L_2$ consistency, instead of $L_1$ consistency.

Theorem \ref{thm:consistency} implies the following prediction consistency result.
\begin{corollary}[Consistency of individual prediction]\label{cor:individual}
Suppose $\bx\subnew$ is sub-Gaussian random vector satisfying
$\sup_{\|\bv\|_2 = 1} \bv\trans \E[\bx\subnew \bx\subnew\trans] \bv \le \sigma^2\submax$.
Under the conditions of Theorem \ref{thm:consistency},
we have
\begin{equation*}
  g\left(\hat{\bbeta}\trans \bx\subnew\right) - g\left(\bbeta\subo\trans \bx\subnew\right)= O_p\left(\|\hat{\bbeta}-\bbeta\subo\|_2\right)=o_p(1).
\end{equation*}
\end{corollary}
The concentration result of Corollary \ref{cor:individual} is established with respect to the joint distribution of the data and the new observation $\bx\subnew$. This is in a sharp contrast to the individual prediction conditioning on any new observation $\bx\subnew.$ If the goal is to conduct inference for any given $\bx\subnew$, the theoretical justification is provided in the following Theorem \ref{thm:inference}
and Corollary \ref{cor:CI}.

\subsection{$\sqrt{n}$-inference with Debiased SAS Estimator}\label{section:thm-inf}

We state the validity of our SSL inference in Theorem \ref{thm:inference}.
We use to $A \stackrel{\mathscr{L}}{\rightarrow} B$
to denote that random variable $A$ converges in distribution
to a distribution $B$.
\begin{theorem}[SAS Inference]\label{thm:inference}
Let $\bx\subnew$ be the random vector representing the covariate of a new individual.
Under Assumptions \ref{assume:model}, \ref{assume:X} and the dimension condition
\begin{equation}
  (1-\rho)^4 \frac{\sgamma^2\log(p+q)^2}{n} + \frac{\rho(\sbeta^2+\sbeta\su)\log(p)^2}{N}
  +
  (1-\rho)^2 \frac{\sgamma\su\log(p+q)\log(p)}{N}
   = o(1),
   \label{assume:dim-inf}
\end{equation}
we draw inference on $\bx\subnew\trans \bbeta\subo$ conditionally on $\bx\subnew$
according to
\begin{equation*}
  \sqrt{n}\hVSAS ^{-1/2}\left(\widehat{\bx\substd\trans \bbeta}  - \frac{\bx\subnew\trans \bbeta\subo}{\|\bx\subnew\|_2}\right) \mid \bx\subnew
\stackrel{\mathscr{L}}{\rightarrow} N(0,1),
\end{equation*}
where $\hVSAS ^2$
defined in \eqref{def:var-cf} is the estimator of the asymptotic variance
\begin{gather}
\begin{aligned}
\VSAS =  & \E[(\bu\subo\trans \bX_i)^2\{Y - (1-\rho)\cdot g(\bgamma\subo\trans \bW_i)-\rho\cdot g(\bbeta\subo\trans \bX_i)\}^2] \\
& + \rho(1-\rho)\E[(\bu\subo\trans \bX_i)^2\{g(\bgamma\subo\trans \bW_i) - g(\bbeta\subo\trans \bX_i)\}^2],\end{aligned} \notag \\
\text{with}\quad   \bu\subo =  \Prec \frac{\bx\subnew}{\|\bx\subnew\|_2} =  [\E\{g'(\bbeta\subo\trans \bX_i)\bX_i\bX_i\trans\}]^{-1}\frac{\bx\subnew}{\|\bx\subnew\|_2}.\label{def:sig0}
\end{gather}
\end{theorem}

By the Young's inequality, the condition \eqref{assume:dim-inf} is implied by
\begin{equation}\label{eq:dim-inf-Young}
     (1-\rho)^4 \frac{\sgamma^2\log(p+q)^2}{n} + \frac{\sqrt{\rho}(\sbeta+\su)\log(p)}{\sqrt{N}}
   = o(1),
\end{equation}
When $p$ is much smaller than the full sample size $N$, our condition \eqref{eq:dim-inf-Young} allows the sparsity levels of
 $\bbeta\subo$ and $\bu\subo$ to be as large as $p$. Even if $p$ is larger than $N$, our SAS inference procedure is valid if
$\sbeta + \su \lesssim \sqrt{N}/\log(p).$ In the literature on confidence interval construction in  high-dimensional supervised setting, the valid inference procedure for a single regression coefficient in the linear regression requires $\sbeta \lesssim \sqrt{n}/\log(p)$ \citep{ZhangZhang14,JavanmardMontanari14,vdGeerEtal14}. Such a sparsity condition has been shown to be necessary to construct a confidence interval of a parametric rate \citep{cai2017confidence}. We have leveraged the unlabeled data to significantly relax the fundamental limit of statistical inference from $\sbeta \lesssim \sqrt{n}/\log(p)$ to $\sbeta \lesssim \sqrt{N}/\log(p).$ The amount of labelled data validates the statistical inference for a dense model in high dimensions.

The sparsity of $\bu\subo$ is determined by $\bx\subnew$ and the precision matrix $\Prec.$
In the supervised learning setting, for confidence interval construction for a single regression coefficient, \cite{vdGeerEtal14} requires $\su \lesssim {n}/\log(p)$ is required. According to \eqref{eq:dim-inf-Young}, our SAS inference requires
$\su \lesssim \sqrt{N}/\log(p),$ which can be weaker than $\su \lesssim {n}/\log(p)$ if the amount of unlabeled data is larger than $n^2.$
Theorem \ref{thm:inference} implies that our proposed CI in \eqref{def:CI} is valid in terms of coverage, which is summarized in the following corollary. 

\begin{corollary}\label{cor:CI}
Under Assumptions \ref{assume:model} and \ref{assume:X}, as well as \eqref{assume:dim-inf},
the CI defined in \eqref{def:CI} satisfies, 
\begin{align*}
&\P \left\{g\left( \|\bx\subnew\|_2\left(\widehat{\bx\substd\trans \bbeta}  -\mathcal{Z}_{\alpha/2} \sqrt{\hVSAS /n}\right)\right)
\le g\left(\bx\subnew\trans \bbeta\subo\right)\right.\\
& \qquad \le \left.g\left(\left(\|\bx\subnew\|_2 \widehat{\bx\substd\trans \bbeta}  +\mathcal{Z}_{\alpha/2} \sqrt{\hVSAS /n}\right)\right)\right\}
= 1 - \alpha + o(1).
\end{align*}
$$
2g'(\bbeta\subo\trans\bx\subnew)\|\bx\subnew\|_2\mathcal{Z}_{\alpha/2}\sqrt{\VSAS/n} \lesssim \|\bx\subnew\|_2/\sqrt{n},
$$
where $\VSAS$ is the the asymptotic variance defined in \eqref{def:sig0}.
\end{corollary}

Confidence interval construction for $g\left(\bx\subnew\trans \bbeta\subo\right)$ in high-dimensional supervised setting has been recently studied in \cite{guo2020inference}.  \cite{guo2020inference} assumes the prediction model to be correctly specified as a high-dimensional sparse logistic regression and the inference procedure is valid if $\sbeta\lesssim \sqrt{n}/\log p.$ In contrast, we leverage the unlabeled data to allow for mis-specified prediction model and a dense regression vector, as long as the dimension requirement in \eqref{assume:dim-inf} is satisfied.

\subsection{Efficiency comparison of SAS Inference}\label{section:eff}

Efficiency in high-dimensional setting or SSL setting in which the proportion of labelled data decays to zero
is yet to be formalized.
Here we use the efficiency bound in the classical low-dimensional with a fixed $\rho$
as the benchmark.
Apart from the relaxation of various sparsity conditions,
we illustrate next that our SAS inference achieves a decent efficiency with properly specified imputation model compared to the supervised learning
and the benchmark.

Similar to the phenomenon discovered by \cite{ChakraborttyCai2018}, if the imputation model is correct,
we can guarantee the efficiency gain by SAS inference in comparison to the asymptotic variance of the supervised learning,
\begin{equation}\label{def:var-super}
    \VSL = \E[(\bu\subo\trans\bX_i)^2\{Y_i - g(\bbeta\subo\trans\bX_i)\}^2].
\end{equation}
\begin{proposition}\label{prop:eff-vs-super}
If $\E(Y_i \mid \bS_i, \bX_i) = g(\bgamma\subo\trans\bW_i)$,
we have $\VSL \ge \VSAS$.
\end{proposition}
Moreover, we can show that our SAS inference attains the benchmark efficiency
derived from classical fixed $\rho$ setting \citep{Tsiatis2007book}.
To simplify the derivation, we describe the missing-completely-at-random mechanism
through the binary observation indicator $R_i$, $i=1,\dots,N$,
independent of $Y_i$, $\bX_i$ and $\bS_i$.
We still denote the proportion of labelled data as $\rho = \E(R_i)$.
The unsorted data take the form
$$
\Dscr = \left\{\bD_i = (\bX_i\trans, \bS_i\trans, R_i, R_iY_i)\trans, i=1,\dots, N\right\}.
$$
We consider the following class of complete data semi-parametric models
\begin{align}
    \Mcomp = \bigg\{ & f_{\bX,Y,\bS, R} (\bx,y,\bs, r) = f_{\bX}(\bx) \rho^r(1-\rho)^{1-r} f_{Y|\bS,\bX}(y|\bs,\bx) f_{\bS|\bX}(\bs|\bx): \notag
    \\ & \qquad  f_{Y|\bS,\bX}, f_{\bX}, f_{\bS|\bX} \text{ are arbitrary density}\bigg\}, \label{def:semiM}
\end{align}
and establish the efficiency bounds for RAL estimators under $\Mcomp$
by deriving the associated efficient influence function in the following proposition.
We denote the nuisance parameters for $f_{Y|\bS,\bX}$, $f_{\bX}$ and $f_{\bS|\bX}$
as $\bge$. We use $\bge\subo$ to denote the true underlying nuisance parameter that generates the data.
The parameter of interest $\bbeta_0$ is not part of the model $\Mcomp$ but defined by the implicit function
through the moment condition \eqref{model:Y_X-moment}.
\begin{proposition}\label{prop:eff-infl}
The efficient influence function for $\theta = \bx\substd\trans\bbeta$ under $\Mcomp$ is
$$
\infleff(\bD_i; \theta\subo, \bge\subo) =  \frac{R_i}{\rho} \bu\subo\trans\bX_i\{Y_i-
\E(Y_i \mid \bS_i, \bX_i)\} - \bu\subo\trans\bX_i\{
\E(Y_i \mid \bS_i, \bX_i)-g(\bbeta\subo\trans\bX_i)\}.
$$
\end{proposition}
Under the Assumptions of Theorem \ref{thm:inference}
and additionally $\E(Y_i \mid \bS_i, \bX_i) = g(\bgamma\subo\trans\bW_i)$,
our SAS debiased estimator admits the same influence function
$$
\widehat{\bx\substd\trans \bbeta}  - \frac{\bx\subnew\trans \bbeta\subo}{\|\bx\subnew\|_2}
= \frac{1}{N}\sum_{i=1}^N\infleff(\bD_i; \theta\subo, \bge\subo) + o_p\left((\rho N)^{-1/2}\right)
$$
according to Appendix \ref{asection:proof-inf} Step 2
\eqref{eq:eff-approx}.

\section{Simulation}\label{section:simulation}
We have conducted extensive simulation studies to evaluate the finite sample performance of the SAS estimation and inference procedures
under various scenarios. Throughout, we let $p=500$, $q=100$,  $N=20000$ and consider $n = 500$. The signals in $\bbeta$ are varied to be approximately sparse or fully dense with a mixture of strong
 and weak signals. The surrogates $\bS$ are either moderately and strongly predictive of $Y$ as specified below. For each configuration, we summarize the results based on 500 simulated datasets.

 To mimic the zero-inflated discrete distribution of EHR features, we first generate $Z^{\scriptscriptstyle x}_{i,1},$ $\ldots,Z^{\scriptscriptstyle x}_{i,p}, Z^u_i, Z^s_{i,1}, \dots, Z^s_{i,q}$ independently from $N(0,25)$. Then we construct $\bX_i$ from $\{Z_{i}^u, \bZ^{\scriptscriptstyle x}_{i}=(Z^{\scriptscriptstyle x}_{i,1}, ..., Z^{\scriptscriptstyle x}_{i,1})\trans\}$
via the transformation $\varsigma(z)=\lfloor \log\{1+\exp(z)\}\rfloor$:
\begin{gather*}\textstyle
X_{i,1} = \left\{ \varsigma\left(\sum_{j=2}^p2 X_{i,j}/\sqrt{p-1}+ Z^{\scriptscriptstyle x}_{i,1}/\sqrt{2}\right)- \mu_{\scriptscriptstyle \bX}\right\}/\sigma_{\scriptscriptstyle \bX},\\
X_{i,j} = [\varsigma(Z^{\scriptscriptstyle x}_{i,j}\sqrt{1-p^{-1}}+ Z^u_i/\sqrt{p}) - \mu_{\scriptscriptstyle \bX}]/\sigma_{\scriptscriptstyle \bX}, \; j = 2, \dots, p.
\end{gather*}
We standardize $X_{i,j}$ to roughly mean zero and unit variance with $\mu_{\scriptscriptstyle \bX} = 1.80$ and $\sigma_{\scriptscriptstyle \bX} = 2.74$. The shared term $Z^u_i$ induces correlation among the covariates.

For $\bS$ and $Y$, we consider two scenarios under which the imputation model is either correctly or incorrectly specified.
We present the ``Scenario I: neither the risk prediction model nor
the imputation model is correctly specified'' in the main text
and the ``Scenario II: The imputation model is correctly specified and
exactly sparse'' in Section \ref{asection:sim} of the Supplementary materials.

\paragraph{Scenario I: neither the risk prediction model nor
the imputation model is correctly specified.} In this scenario, we first generate $Y_i$ from  the probit model
$$
\mathbb{P}(Y_i=1| \bZ_i^x)=
\Phi(\bga\trans \bZ_i^x) \quad \text{with}\quad
\Phi(x) = \int_{-\infty}^x (2\pi)^{-1/2} e^{-x^2/2} dx,
$$
and then generate $\bS$ from
$$
S_{i,1} = \left\{\varsigma (Z^s_{i,1}/2+ \theta Y_i) - \mu_{\scriptscriptstyle \bS}\right\}\sigma_{\scriptscriptstyle \bS}^{-1} + \bgx\trans \bX_i , \quad \mbox{and} \quad S_{i,j} = \{\varsigma(Z^s_{i,j}) - \mu_{\scriptscriptstyle \bX}\}\sigma_{\scriptscriptstyle \bX}^{-1}, \ j=2, \ldots, p.
$$
We chose $\mu_{\scriptscriptstyle \bS} $ and $\sigma_{\scriptscriptstyle \bS}$ depending on $\bga$ such that $S_{i,1}$ is roughly mean 0 and variance 1. Under this setting, a logistic imputation model would be misspecified but nevertheless approximately sparse with appropriately chosen  $\bgx$.
The coefficients $\bga$ control the optimal prediction accuracy of $\bX$ for $Y$ while $\theta$  controls the optimal prediction accuracy of $\bS$ for $Y$.  We consider two $\bga$ of different sparsity patterns,
which also determine the rest of parameters
\begin{alignat*}{2}
& \text{Sparse }(\salpha = 3): &\quad& \bga = (0.45,0.318,0.318,\mathbf{0}_{497\times 1}\trans)\trans, \,
\mu_{\scriptscriptstyle \bS} = 1.82, \,
\sigma_{\scriptscriptstyle \bS} = 2.01, \\
&\text{Dense }(\salpha = 500): &\quad& \bga = (0.316,\mathbf{0.059}_{29 \times 1}\trans,\mathbf{0.007}_{470 \times 1}\trans)\trans, \,
\mu_{\scriptscriptstyle \bS} = 2.71, \,
\sigma_{\scriptscriptstyle \bS} = 2.68,
\end{alignat*}
where $\ba_{k\times 1}=(a, ..., a)_{k\times 1}\trans$ for any $a$.
The sparsity of $\bga$ affects the approximate sparsity of
$\bgb$ subsequently (Table \ref{tab:sim-auc}),
which we measured by the squared ratio between $\ell_1$ norm and
$\ell_2$ norm
\begin{equation}\label{def:approx_S}
    \Scal(\bgb) = \|\bgb\|_1^2 / \|\bgb\|_2^2, \;
  \min_{j:\beta_j\neq 0} |\beta_j|   \le  \Scal(\bgb)/ \|\bgb\|_0 \le 1.
\end{equation}
We consider two $\theta$: (a) $\theta = 0.6$  for $\bS$ to be moderately predictive of $Y$; and (b) $\theta = 1$ for strong surrogates. The parameter $\bgx$ depends on both the choices of $\bga$ and $\theta$:
\begin{alignat*}{4}
& \salpha = 3, \, \theta = 0.6 : &\quad& \bgx = (0.407,0.330,0.330,\mathbf{0.005}_{497\times 1}\trans)\trans, \\
& \salpha = 3, \, \theta = 1 : &\quad& \bgx = (0.199,0.163,0.163,\mathbf{0.002}_{497\times 1}\trans)\trans, \\
& \salpha = 500, \, \theta = 0.6 : &\quad& \bgx = (0.350,\mathbf{0.064}_{29 \times 1}\trans,\mathbf{0.011}_{470 \times 1}\trans)\trans, \\
& \salpha = 500, \, \theta = 1 : &\quad& \bgx = (0.169,\mathbf{0.032}_{29 \times 1}\trans,\mathbf{0.005}_{470 \times 1}\trans)\trans.  \\
\end{alignat*}

Due to the complexity of the data generating process
and the noncollapsibility of the logistic regression
models, we cannot analytically express the true $\bbeta\subo$
in both scenarios.
Instead, we numerically evaluate $\bbeta\subo$
with a large simulated data using the oracle knowledge of
the ex-changeability among covariates
according to the model
$$
\mathrm{logit}\{\mathbb{P}(Y_i=1| S_{i,1})\} \sim \eta\subo + \eta_1 X_{i,1} + \eta_2\sum_{j=2}^{\ssbeta } X_{i,j} +
\eta_3 \sum_{j=\ssbeta+1}^p X_{i,j}.
$$
We derive the true $\bgb\subo$ as
$$
\bgb\subo = (\eta_0, \eta_1, (\boldsymbol{\eta_2})\trans_{\ssbeta \times 1}, (\boldsymbol{\eta_3})\trans_{(p-\ssbeta) \times 1} )\trans.
$$

We report the simulation settings under Scenario I in Table \ref{tab:sim-auc},
where we present the predictive power of the oracle estimation and
the lasso estimation.
We also report the average area-under-curve (AUC) of the receiver operating characteristic (ROC)
curve
for oracle $\bbeta\subo$, supervised LASSO (SLASSO)
and the proposed SAS estimation.
Our SAS estimation achieves a better AUC compared to supervised LASSO across all
scenarios, and is comparable to the AUC with the true coefficient $\bbeta\subo$.
Besides, we observe that the AUC of supervised LASSO is sensitive to the
approximate sparsity $\Scal(\bgb\subo)$, while the AUC of
SAS estimation does not seem to be affected by $\Scal(\bgb\subo)$.

\begin{table}
\begin{center}
\caption{AUC Table for simulations with 500 labels under Scenario I. The AUCs are evaluated on an independent testing set of size $100$. We approximately measure the sparsity
by $\Scal(\bv) = \|\bv\|^2_1/\|\bv\|_2^2$.
}\label{tab:sim-auc}
\scriptsize
\begin{tabular}{lll|lll}
\hline
\hline
 \multicolumn{3}{c|}{Scenario}&  \multicolumn{3}{c}{Prediction Accuracy (AUC)}\\
Surrogate & $\Scal(\bbeta\subo)$
& $\Scal(\bgamma\subo)$
& Oracle & SLASSO
& SAS \\
  \hline
Strong & 174 & 1.32 & 0.724 & 0.660  & 0.711  \\
  Moderate & 174 & 1.26 & 0.724 & 0.660 &  0.713  \\
  Strong & 28.3 & 1.33 & 0.719 & 0.694  & 0.713  \\
  Moderate & 28.3 & 1.24 & 0.719 & 0.694  & 0.711  \\
   \hline
\hline
\end{tabular}
\end{center}
\end{table}

To evaluate the SAS inference for the individualized prediction, we consider six different choices of $\bx\subnew$. We first select $\{\bx\subnew\supL,\bx\subnew\supM,\bx\subnew\supH\}$ from a random sample of $\bx\subnew$ generated from the distribution of $\bX_i$ such that their predicted risks are around $0.2$, $0.5$, and $0.7$,  corresponding to low, moderate and high risk. We additionally consider three sets of $\bx\subnew$ with different levels of sparsity:
\begin{alignat*}{2}
&\text{Sparse: } &\quad&  \bx\subnew\supS = (1,1,\mathbf{0}_{499\times 1}\trans)\trans; \\
& \text{Intermediate: }  &\quad&  \bx\subnew\supI =  (1,\mathbf{0.183}_{30\times 1}\trans,\mathbf{0}_{470\times 1}\trans)\trans; \\
&\text{Dense: } &\quad&  \bx\subnew\supD = (1,\mathbf{0.045}_{500\times 1}\trans)\trans.
\end{alignat*}
In Table \ref{table:sim-est-1},
we compare our SAS estimator of $\bx\subnew\trans\bgb_0$ with the corresponding SLASSO across all settings under Scenario I.
The root mean-squared-error (rMSE)
of the SAS estimation decays proportionally with the sample size,
while the rMSE of the supervised LASSO provides evidence of inconsistency
for moderate and dense  deterministic $\bx\subnew$.
The bias of the supervised LASSO is also significantly larger than that of
the SAS estimation.
The performance
of the SAS estimation is insensitive to sparsity of $\bbeta\subo$,
while that of supervised LASSO severely deteriorate with dense $\bbeta\subo$.
The improvement from the supervised LASSO
to the SAS estimation is regulated by the surrogate strength.

\begin{table}
\begin{center}
\caption{Comparison of SAS Estimation to the  supervised  LASSO (SLASSO) with Bias, Empirical standard error (ESE) and root mean-squared error (rMSE) of the linear predictions $\bx\subnew\trans \bbeta\subo$ under Scenario I 500 labels, moderate or large $\Scal(\bbeta\subo)$  and strong or moderate surrogates.
}\label{table:sim-est-1}
\scriptsize
\begin{tabular}{l|lll|lll|lll}
\hline
\hline
& \multicolumn{3}{c|}{SLASSO}&\multicolumn{3}{c|}{SAS: Moderate } &\multicolumn{3}{c}{SAS: Strong } \\
\hline
Type & Bias & ESE & rMSE & Bias & ESE & rMSE& Bias & ESE & rMSE \\
\hline
\multicolumn{10}{c}{Moderate $\Scal(\bbeta\subo)$ } \\
\hline
 $\bx\subnew\supL$ & 0.605 & 0.387 & 0.719 & 0.165 & 0.249 & 0.298 & 0.118 & 0.196 & 0.229 \\
  $\bx\subnew\supM$ & -0.083 & 0.337 & 0.347 & -0.008 & 0.246 & 0.246 & -0.016 & 0.195 & 0.196 \\
  $\bx\subnew\supH$ & -0.718 & 0.521 & 0.887 & -0.234 & 0.294 & 0.376 & -0.176 & 0.225 & 0.286 \\
   \hline
$\bx\subnew\supS$ & -0.072 & 0.144 & 0.161 & -0.080 & 0.094 & 0.123 & -0.018 & 0.078 & 0.080 \\
  $\bx\subnew\supI$ & -0.460 & 0.096 & 0.470 & -0.110 & 0.093 & 0.143 & -0.055 & 0.071 & 0.090 \\
  $\bx\subnew\supD$ & -0.413 & 0.091 & 0.423 & -0.110 & 0.089 & 0.141 & -0.114 & 0.069 & 0.133 \\

\hline
\multicolumn{10}{c}{Large $\Scal(\bbeta\subo)$} \\
\hline
 $\bx\subnew\supL$ & 0.389 & 0.275 & 0.477 & 0.161 & 0.215 & 0.269 & 0.133 & 0.264 & 0.296 \\
  $\bx\subnew\supM$ & -0.017 & 0.280 & 0.280 & -0.014 & 0.213 & 0.213 & -0.017 & 0.268 & 0.268 \\
  $\bx\subnew\supH$ & -0.600 & 0.481 & 0.769 & -0.251 & 0.271 & 0.370 & -0.164 & 0.296 & 0.339 \\
   \hline
$\bx\subnew\supS$ & -0.202 & 0.140 & 0.246 & -0.074 & 0.097 & 0.122 & -0.009 & 0.078 & 0.079 \\
  $\bx\subnew\supI$ & -0.178 & 0.098 & 0.203 & -0.075 & 0.086 & 0.115 & -0.071 & 0.075 & 0.103 \\
  $\bx\subnew\supD$ & -0.185 & 0.090 & 0.206 & -0.109 & 0.084 & 0.138 & -0.113 & 0.073 & 0.135 \\

\hline
\hline
\end{tabular}
\end{center}
\end{table}

In Table \ref{table:sim-inf-1},
we compare our SAS inference  with supervised  debiased LASSO
across the settings under Scenario I.
Our SAS inference procedure attains approximately honest coverage of 95 \% confidence intervals
for all types of $\bx\subnew$ under all scenarios.
Unsurprisingly, the debiased SLASSO has under coverage for the deterministic $\bx\subnew$ as the consequence of violation to the sparsity assumption for $\bbeta\subo$
and precision matrix.
Under our design, the first covariate $X_1$ has the strongest dependence upon the other covariates, whose associated row in the precision matrix is thus densest.
Consequently, the inference for $\bbeta\trans\bx\subnew\supS = \beta_0 + \beta_1$ 
The debiased SLASSO also has an acceptable coverage for random $\bx\subnew\supL$, $\bx\subnew\supM$, $\bx\subnew\supH$
sampled from the covariate distribution despite the presence of  substantial bias,
which we attribute to the even larger variance that dominates the bias.
In contrast, our  SAS inference has small bias across all scenarios
and improved variance from the strong surrogate.

\begin{table}
\begin{center}
\caption{Bias, Empirical standard error (ESE), average of the estimated standard error (ASE) along with empirical coverage of the 95\% confidence intervals (CP) for the debiased supervised  LASSO (SLASSO) and debiased SAS estimator of linear predictions $\bx\subnew\trans \bbeta\subo$ under Scenario I with 500 labels,  moderate or large $\Scal(\bbeta\subo)$   and strong or moderate surrogates.
}\label{table:sim-inf-1}
\scriptsize
\begin{tabular}{l|llll|llll|llll}
\hline
\hline
& &&&& \multicolumn{8}{c}{Debiased SAS} \\
& \multicolumn{4}{c|}{Debiased SLASSO}&\multicolumn{4}{c}{Moderate Surrogates}
&\multicolumn{4}{c}{Strong Surrogates} \\
\hline
Type & Bias & ESE & ASE & CP & Bias & ESE & ASE & CP & Bias & ESE & ASE & CP\\
\hline
\multicolumn{13}{c}{Risk prediction model approximatedly sparse} \\
\hline
 $\bx\subnew\supL$ & -0.290 & 1.901 & 1.896 & 0.948 & 0.021 & 1.873 & 1.864 & 0.949 & 0.018 & 1.531 & 1.531 & 0.950 \\
  $\bx\subnew\supM$ & -0.091 & 1.994 & 1.981 & 0.947 & -0.007 & 1.961 & 1.954 & 0.950 & -0.015 & 1.560 & 1.570 & 0.953 \\
  $\bx\subnew\supH$ & 0.348 & 2.106 & 2.074 & 0.942 & -0.050 & 2.036 & 2.039 & 0.950 & -0.011 & 1.632 & 1.623 & 0.950 \\
   \hline
$\bx\subnew\supS$ & 0.171 & 0.157 & 0.128 & 0.694 & -0.019 & 0.149 & 0.150 & 0.950 & -0.001 & 0.132 & 0.125 & 0.924 \\
  $\bx\subnew\supI$ & -0.001 & 0.129 & 0.125 & 0.938 & -0.013 & 0.123 & 0.116 & 0.932 & 0.010 & 0.101 & 0.094 & 0.920 \\
  $\bx\subnew\supD$ & 0.141 & 0.137 & 0.138 & 0.812 & -0.011 & 0.123 & 0.118 & 0.944 & -0.001 & 0.096 & 0.095 & 0.940 \\

\hline
\multicolumn{13}{c}{Large $\Scal(\bbeta\subo)$} \\
\hline
 $\bx\subnew\supL$ & -0.134 & 1.918 & 1.914 & 0.951 & 0.018 & 1.875 & 1.878 & 0.951 & 0.018 & 1.529 & 1.524 & 0.948 \\
  $\bx\subnew\supM$ & -0.056 & 1.970 & 1.962 & 0.948 & -0.020 & 1.911 & 1.927 & 0.952 & 0.005 & 1.603 & 1.597 & 0.950 \\
  $\bx\subnew\supH$ & 0.109 & 2.051 & 2.029 & 0.945 & -0.022 & 1.997 & 1.991 & 0.950 & -0.040 & 1.671 & 1.668 & 0.951 \\
   \hline
$\bx\subnew\supS$ & 0.029 & 0.155 & 0.127 & 0.892 & -0.008 & 0.153 & 0.147 & 0.946 & -0.013 & 0.133 & 0.131 & 0.938 \\
  $\bx\subnew\supI$ & 0.002 & 0.131 & 0.125 & 0.930 & 0.001 & 0.122 & 0.114 & 0.936 & 0.002 & 0.101 & 0.098 & 0.936 \\
  $\bx\subnew\supD$ & 0.113 & 0.135 & 0.139 & 0.874 & -0.007 & 0.119 & 0.116 & 0.938 & -0.003 & 0.099 & 0.097 & 0.960 \\

\hline
\hline
\end{tabular}
\end{center}
\end{table}

According to Tables \ref{tab:sim-auc-2},
\ref{table:sim-est-2} and \ref{table:sim-inf-2}
in the Appendix \ref{asection:sim},
the results under Scenario II are consistent with our findings
under Scenario I.

\section{Application of SAS to EHR Study}\label{section:data}

We applied the proposed SAS method to the risk prediction of Type II Diabetes Mellitus (T2DM)
using EHR and genomic data of participants of the Mass General Brigham (MGB) Biobank study. To define the study cohort, we  extracted from the EHR of each patient their date of first EHR encounter ($t_{ini}$), follow up period ($C$), the counts and dates for the ICD codes and NLP mentions of clinical concepts related to T2DM as well as its risk factors. We only included patients who do not have any ICD code or NLP mention of T2DM up to baseline, where the baseline time is defined as 1990 if $t_{ini}$ is prior to 1990 and as their first year if $t_{ini}\ge 1990$. Although neither the ICD code nor NLP mention of T2DM is sufficiently specific, they are highly sensitive and can be used to accurately remove patients who have already developed T2DM at baseline. This exclusion criterion resulted in $N=20216$  patients who are free of T2DM at baseline and have both EHR and genomics features for risk modeling. Among those, we have a total of $n=271$ patients whose T2DM status during follow up, $Y$, has been obtained via manual chart review. The prevalence of T2DM was about 14\% based on labeled data.

We aim to develop a risk prediction model for $Y$ by fitting a working model $P(Y = 1 \mid \bX) = g(\bbeta_0\trans\bX)$, where the baseline covariate vector $\bX$ includes age, gender,
 indicator for occurrence of ICD code and NLP counts for obesity, hypertension, coronary artery disease (CAD), hyperlipidemia during the first year window,  as well as a total of 49 single nucleotide polymorphism (SNP) previously reported as associated with T2DM in \cite{Mahajan2018} with odds ratio greater than 1.1. We additionally adjust for follow up by including $\log(C)$ and allow for non-linear effects by including two-way interactions between the SNPs and other baseline covariates. All variables with less than 10 nonzero values within  the labelled set are removed, resulting the final covariates to be of dimension $p = 260$. We standardize the covariates to have mean 0 and variance 1.  To impute the outcome, we used the predicted probability of T2DM derived from the unsupervised phenotyping method MAP \citep{LiaoEtal19MAP}, which achieves an AUC of $0.98$, indicating a strong surrogate.
In addition to the proposed SAS procedure, we derive risk prediction models based on the supervised LASSO with both the same set of covariates.  We let $\Kn=5$ in cross-fitting and use 5-fold cross-validation for tuning parameter selection. To compare the performance of different risk prediction models, we use 10-fold cross-validation to estimate the out-of-sample AUC.
We repeated the process 10 times and took average of predicted probabilities
across the repeats for each labelled sample and method in comparison.

\begin{figure}
\centering
\includegraphics[width= 0.6\textwidth]{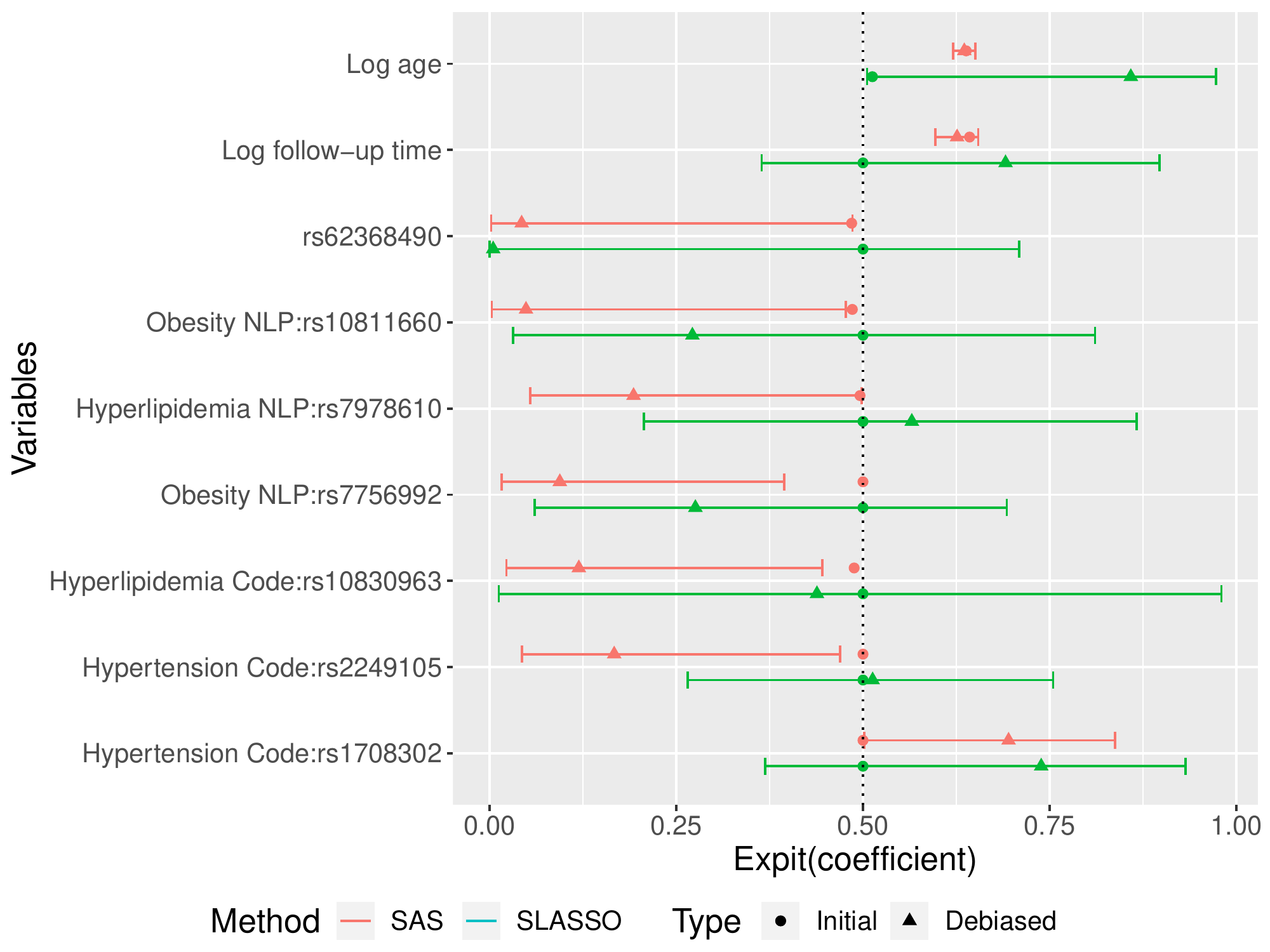}
\caption{Point and 95\% confidence interval estimates for the coefficients with nominal p-value $<0.05$ from SAS inference. The horizontal bars indicate the estimated 95\% confidence intervals. The solid points indicate the (initial) estimates, and the
triangles indicate debiased estimates. Colors red and green
indicate different methods, SAS and SLASSO, respectively.
}\label{fig:data-coef}
\end{figure}

In Figure \ref{fig:data-coef}, we present the estimated $\bbeta$ coefficients for the covariates that received p-value less than $0.05$ from the SAS inference.
The confidence intervals are generally narrower from the SAS inference.
For the coefficients of baseline age and follow-up time, the SAS inference produced much narrower confidence interval than debiased SLASSO, which are expected to have a positive effect on the T2DM onset status during the observation.
In addition, the SAS inference identified one global genetic risk factor
and 6 other subgroup genetic risk factors while SLASSO identified none of these.

\begin{table}
\begin{center}
\caption{The cross-validated (CV) AUC
the estimated risk prediction models with high dimensional EHR and genetic features  based on SAS and supervised LASSO. Shown also are the AUC of the imputation model derived for the SAS procedure.
}\label{table:data-auc}
\footnotesize
\begin{tabular}{l|ccc}
\hline
\hline
Method & Imputation &  SAS & SLASSO\\
\hline
CV AUC & 0.928  &  0.763& 0.488  \\
\hline
\hline
\end{tabular}
\end{center}
\end{table}
In Table \ref{table:data-auc}, we present the AUCs of the estimated risk prediction models using the high dimensional $\bX$ . It is important to note that AUC is a measurement of prediction accuracy, so debiasing
might lead to worse AUC by accepting larger variability for
reduced bias.
The AUC from SLASSO is very poor, probably due to the over-fitting bias with the small sample sizes
of the labeled set.
With the information from a large unlabeled data, SAS produced the significantly higher AUC than the SLASSO.

\def\bgahat{\widehat{\bga}}

\begin{figure}
\centering
\includegraphics[width= \textwidth]{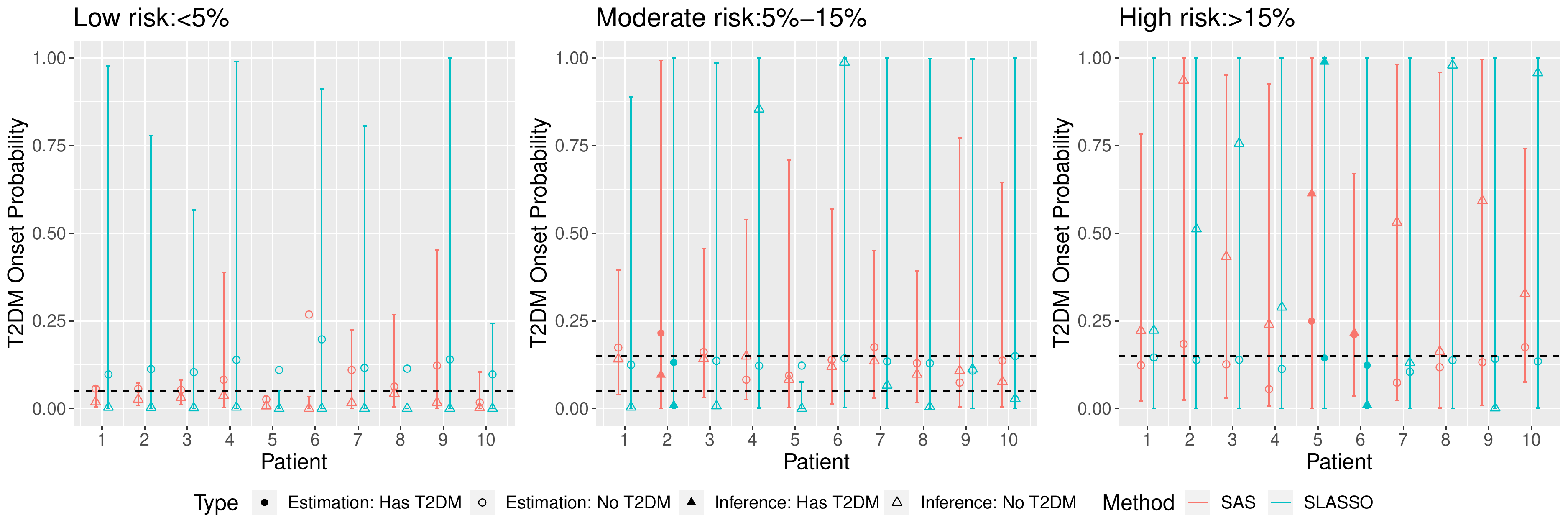}
\caption{Point and 95\% confidence interval estimates for the predicted risks of 30 randomly selected patients. The vertical bars indicate the estimated 95\% confidence intervals. The circle and the triangle shapes correspond to (initial) estimation and debiased estimation, correspondingly.
Solid points indicate the observed T2DM cases.
Colors red and green
indicate different methods, SAS and SLASSO.
}\label{fig:data-individual}
\end{figure}

For illustration, we present in Figure \ref{fig:data-individual} the individual risk predictions with 95\% confidence intervals for three sets of 10 patients with each set randomly selected from low $(<5\%)$, medium $(5\% \sim 15\%)$ or high risk $(> 15\%)$  subgroups. These risk groups are constructed for illustration purposes and a patient with $\bx\subnew$ classified to low, medium and high risk if $\expit(\hat{\bgb}\trans\bx\subnew)$ belongs to the low, medium and high tertiles of $\{\expit(\hat{\bgb}\trans\bX_i),i=1,...,N\}$.
We observe that the confidence intervals for patients with predicted
The debiased SLASSO inference is not very informative with most error bars
stretching from zero to one.
The contrast between SAS CIs and SLASSO CIs
demonstrates the improved efficiency as the result of leveraging information from
the unlabeled data through predictive surrogates.

\bibliographystyle{chicago}
\bibliography{SSL_GLM_imputation}

\begin{appendix}
\renewcommand\thefigure{A\arabic{figure}}
\renewcommand\theequation{A.\arabic{equation}}
\renewcommand\thetable{A\arabic{table}}
\renewcommand\thelemma{A\arabic{lemma}}
\renewcommand\thesubsection{\Alph{section}\arabic{subsection}}
\renewcommand\thedefinition{A\arabic{definition}}
\renewcommand\theremark{A\arabic{remark}}

\setcounter{equation}{0}
\setcounter{lemma}{0}
\setcounter{table}{0}

\clearpage
\newpage

\section*{Supplementary Material}

We present the simulation Scenario II in which the imputation model
is correctly specified and exactly sparse in Appendix \ref{asection:sim}.
The proofs of Theorems \ref{thm:consistency}, \ref{thm:inference}, Corollary \ref{cor:individual}
and Propositions \ref{prop:eff-vs-super} and \ref{prop:eff-infl}
are given in Appendix \ref{asection:main}.
The technical details are put in Appendix \ref{asection:aux}.
Definitions and existing results are stated in Appendix \ref{asection:detail}.

\section{Additional Simulation}\label{asection:sim}

\begin{table}[H]
\begin{center}
\caption{AUC Table for simulations with 500 labels under Scenario II. The AUCs are evaluated on an independent testing set of size $100$. We approximately measure the sparsity
by $\Scal(\bv) = \|\bv\|^2_1/\|\bv\|_2^2$.
}\label{tab:sim-auc-2}
\scriptsize
\begin{tabular}{lll|lll}
\hline
\hline
 \multicolumn{3}{c|}{Scenario}&  \multicolumn{3}{c}{Prediction Accuracy (AUC)}\\
Surrogate & $\Scal(\bbeta\subo)$
& $\Scal(\bgamma\subo)$
& Oracle & SLASSO
& SAS \\
  \hline
Strong & 159 & 1.10 & 0.715 & 0.660 & 0.702  \\
  Moderate & 128 & 1.06 & 0.715 & 0.665  & 0.704  \\
  Strong & 26.4 & 1.09 & 0.710 & 0.691  & 0.708 \\
  Moderate & 18.4 & 1.03 & 0.709 & 0.693  & 0.707  \\
   \hline
\hline
\end{tabular}
\end{center}
\end{table}

\paragraph{Scenario II: The imputation model is correctly specified and
exactly sparse.} In the second scenario, we first generate $\bS_i$ from
$$
S_{i,1} = \left[\varsigma \{\nu Z^s_{i,1}+ \bga\trans (\bZ^{\scriptscriptstyle x}_{i}\sqrt{1-p^{-1}}+ Z^u_i/\sqrt{p})\} - \mu_{\scriptscriptstyle \bS}\right]/\sigma_{\scriptscriptstyle \bS} .
$$
and $S_{i,j} = \{\varsigma(Z^s_{i,j}) - \mu_{\scriptscriptstyle \bX}\}/\sigma_{\scriptscriptstyle \bX}$ for $j=2, \ldots, p$, and then generate $Y_i$ from a sparse model
$$
\mathbb{P}(Y_i=1| \bX_i)=
\mathrm{expit}(\theta S_{i,1}).
$$
We chose $\mu_{\scriptscriptstyle \bS} \approx 0.66$ and $\sigma_{\scriptscriptstyle \bS} \approx 1$ such that $S_{i,1}$ is roughly mean 0 and variance 1.
Under this setting, the imputation model holds with $\sgamma=1$.
The factor $\nu$ and the coefficients $\bga$ control the predictiveness of $\bX$ for $S_1$ and $Y$ while $\theta$ controls the predictiveness of $\bW$ for $Y$. We consider two $\bga$ of different sparsity patterns,
\begin{alignat*}{2}
& \text{Sparse }(\salpha = 3): &\quad& \bga = (0.3,0.212,0.212,\mathbf{0}_{497\times 1}\trans)\trans \\
&\text{Dense }(\salpha = 500): &\quad& \bga = (0.211,\mathbf{0.039}_{29 \times 1}\trans,\mathbf{0.004}_{470 \times 1}\trans)\trans,
\end{alignat*}
where $\ba_{k\times 1}=(a, ..., a)_{k\times 1}\trans$ for any $a$.
Similar to Scenario I, the sparsity of $\bga$ regulates the
approximate sparsity of $\bgb$ measured by \eqref{def:approx_S}
(See Table \ref{tab:sim-auc}).
We consider two sets of $(\nu,\theta)$ to allow $\bW$ to be either moderately or strongly predictive of $Y$:
$$
\text{Moderate: } \nu = 0.4, \, \theta = 2; \; \qquad \mbox{and}\qquad
\text{Strong: } \nu = 0.6, \, \theta = 3.7.
$$

The layouts of Tables \ref{table:sim-est-2} and \ref{table:sim-inf-2}
are different from those of \ref{table:sim-est-1} and Table \ref{table:sim-inf-1}
because of the different data generating mechanism.
The distribution of $Y_i \mid  \bX_i$ is not affected by the distribution of $\bS_i$
in Scenario I, while the property does not hold in Scenario II.

\begin{table}
\begin{center}
\caption{Comparison of SAS Estimation to the  supervised  LASSO (SLASSO) with Bias, Empirical standard error (ESE) and root mean-squared error (rMSE) of the linear predictions $\bx\subnew\trans \bbeta\subo$ under Scenario II with 500 labels, moderate or large $\Scal(\bbeta\subo)$   and strong or moderate surrogates.
}\label{table:sim-est-2}
\scriptsize
\begin{tabular}{l|lll|lll|lll|lll}
\hline
\hline
& \multicolumn{6}{c|}{Moderate Surrogates}&\multicolumn{6}{c}{Strong Surrogates}\\
\hline
& \multicolumn{3}{c}{SLASSO}&\multicolumn{3}{c|}{SAS}&
\multicolumn{3}{c}{SLASSO} &\multicolumn{3}{c}{SAS} \\
\hline
Type & Bias & ESE & rMSE & Bias & ESE & rMSE& Bias & ESE & rMSE& Bias & ESE & rMSE \\
\hline
\multicolumn{13}{c}{Risk prediction model approximatedly sparse} \\
\hline
 $\bx\subnew\supL$ & 0.505 & 0.378 & 0.631 & 0.163 & 0.283 & 0.327 & 0.349 & 0.278 & 0.446 & 0.085 & 0.222 & 0.238 \\
  $\bx\subnew\supM$ & -0.140 & 0.331 & 0.359 & -0.047 & 0.272 & 0.276 & -0.113 & 0.282 & 0.304 & -0.058 & 0.217 & 0.225 \\
  $\bx\subnew\supH$ & -0.713 & 0.512 & 0.878 & -0.262 & 0.313 & 0.408 & -0.678 & 0.469 & 0.825 & -0.210 & 0.241 & 0.320 \\
   \hline
$\bx\subnew\supS$ & -0.111 & 0.143 & 0.181 & -0.058 & 0.081 & 0.100 & -0.190 & 0.142 & 0.237 & -0.037 & 0.063 & 0.072 \\
  $\bx\subnew\supI$ & -0.437 & 0.098 & 0.448 & -0.119 & 0.076 & 0.141 & -0.155 & 0.098 & 0.183 & -0.065 & 0.061 & 0.089 \\
  $\bx\subnew\supD$ & -0.349 & 0.093 & 0.361 & -0.138 & 0.078 & 0.158 & -0.150 & 0.093 & 0.176 & -0.112 & 0.063 & 0.129 \\

\hline
\multicolumn{13}{c}{Large $\Scal(\bbeta\subo)$} \\
\hline
 $\bx\subnew\supL$ & 0.366 & 0.266 & 0.453 & 0.142 & 0.224 & 0.265 & 0.482 & 0.398 & 0.625 & 0.117 & 0.300 & 0.322 \\
  $\bx\subnew\supM$ & -0.060 & 0.275 & 0.282 & -0.035 & 0.213 & 0.216 & -0.199 & 0.337 & 0.391 & -0.082 & 0.299 & 0.310 \\
  $\bx\subnew\supH$ & -0.656 & 0.475 & 0.810 & -0.272 & 0.257 & 0.374 & -0.749 & 0.503 & 0.903 & -0.214 & 0.325 & 0.389 \\
   \hline
$\bx\subnew\supS$ & -0.236 & 0.139 & 0.274 & -0.087 & 0.079 & 0.117 & -0.054 & 0.138 & 0.148 & 0.003 & 0.063 & 0.063 \\
  $\bx\subnew\supI$ & -0.173 & 0.096 & 0.197 & -0.078 & 0.077 & 0.109 & -0.409 & 0.097 & 0.420 & -0.094 & 0.057 & 0.110 \\
  $\bx\subnew\supD$ & -0.144 & 0.092 & 0.171 & -0.101 & 0.080 & 0.129 & -0.359 & 0.093 & 0.371 & -0.154 & 0.060 & 0.166 \\

\hline
\hline
\end{tabular}
\end{center}
\end{table}

\begin{table}
\begin{center}
\caption{Bias, Empirical standard error (ESE) along with empirical coverage of the 95\% confidence intervals (CP) for the debiased supervised  LASSO (SLASSO) and debiased SAS estimator of linear predictions $\bx\subnew\trans \bbeta\subo$ under Scenario II with 500 labels,  moderate or large $\Scal(\bbeta\subo)$   and strong or moderate surrogates.
}\label{table:sim-inf-2}
\scriptsize
\begin{tabular}{l|llll|llll}
\hline
\hline
& \multicolumn{4}{c|}{Debiased SLASSO}&\multicolumn{4}{c}{Debiased SAS}\\
\hline
Type & Bias & ESE & ASE & CP & Bias & ESE  & ASE & CP \\
\hline
\multicolumn{9}{c}{Risk prediction model approximatedly sparse, moderate surrogates} \\
\hline
 $\bx\subnew\supL$ & -0.236 & 1.936 & 1.915 & 0.947 & 0.014 & 1.786 & 1.771 & 0.950 \\
  $\bx\subnew\supM$ & -0.044 & 2.031 & 1.997 & 0.944 & -0.028 & 1.873 & 1.853 & 0.947 \\
  $\bx\subnew\supH$ & 0.364 & 2.110 & 2.084 & 0.944 & -0.045 & 1.943 & 1.924 & 0.947 \\
   \hline
$\bx\subnew\supS$ & 0.133 & 0.156 & 0.127 & 0.784 & -0.028 & 0.133 & 0.130 & 0.944 \\
  $\bx\subnew\supI$ & 0.004 & 0.124 & 0.126 & 0.942 & -0.014 & 0.102 & 0.100 & 0.936 \\
  $\bx\subnew\supD$ & 0.149 & 0.121 & 0.139 & 0.848 & -0.014 & 0.104 & 0.105 & 0.948 \\

\hline
\multicolumn{9}{c}{Risk prediction model approximatedly sparse, strong surrogates} \\
\hline
 $\bx\subnew\supL$ & -0.070 & 1.953 & 1.935 & 0.947 & 0.021 & 1.371 & 1.366 & 0.949 \\
  $\bx\subnew\supM$ & -0.031 & 2.019 & 1.986 & 0.946 & -0.026 & 1.408 & 1.401 & 0.949 \\
  $\bx\subnew\supH$ & 0.148 & 2.073 & 2.055 & 0.948 & -0.010 & 1.458 & 1.444 & 0.949 \\
   \hline
$\bx\subnew\supS$ & 0.029 & 0.153 & 0.127 & 0.894 & -0.016 & 0.103 & 0.096 & 0.928 \\
  $\bx\subnew\supI$ & 0.018 & 0.134 & 0.126 & 0.938 & -0.004 & 0.081 & 0.079 & 0.944 \\
  $\bx\subnew\supD$ & 0.134 & 0.128 & 0.141 & 0.842 & -0.007 & 0.081 & 0.083 & 0.956 \\

\hline
\multicolumn{9}{c}{Large $\Scal(\bbeta\subo)$, moderate surrogates} \\
\hline
 $\bx\subnew\supL$ & -0.092 & 1.942 & 1.925 & 0.950 & 0.004 & 1.796 & 1.792 & 0.951 \\
  $\bx\subnew\supM$ & -0.034 & 1.995 & 1.969 & 0.947 & -0.018 & 1.852 & 1.835 & 0.951 \\
  $\bx\subnew\supH$ & 0.082 & 2.061 & 2.036 & 0.946 & -0.027 & 1.912 & 1.890 & 0.948 \\
   \hline
$\bx\subnew\supS$ & -0.009 & 0.155 & 0.125 & 0.876 & -0.027 & 0.131 & 0.125 & 0.922 \\
  $\bx\subnew\supI$ & 0.000 & 0.126 & 0.125 & 0.952 & -0.009 & 0.104 & 0.103 & 0.950 \\
  $\bx\subnew\supD$ & 0.119 & 0.126 & 0.139 & 0.894 & -0.012 & 0.108 & 0.108 & 0.940 \\

\hline
\multicolumn{9}{c}{Large $\Scal(\bbeta\subo)$, strong surrogates} \\
\hline
 $\bx\subnew\supL$ & -0.221 & 1.929 & 1.926 & 0.949 & 0.022 & 1.353 & 1.349 & 0.951 \\
  $\bx\subnew\supM$ & 0.032 & 2.047 & 2.017 & 0.947 & -0.003 & 1.427 & 1.414 & 0.950 \\
  $\bx\subnew\supH$ & 0.442 & 2.137 & 2.104 & 0.940 & -0.039 & 1.479 & 1.469 & 0.951 \\
   \hline
$\bx\subnew\supS$ & 0.176 & 0.150 & 0.128 & 0.698 & -0.018 & 0.094 & 0.099 & 0.946 \\
  $\bx\subnew\supI$ & 0.030 & 0.128 & 0.129 & 0.936 & -0.002 & 0.079 & 0.077 & 0.952 \\
  $\bx\subnew\supD$ & 0.167 & 0.125 & 0.142 & 0.804 & -0.008 & 0.082 & 0.080 & 0.954 \\

\hline
\hline
\end{tabular}
\end{center}
\end{table}

\newpage

\section{Proofs of Main Results}\label{asection:main}
\allowdisplaybreaks

We first summarize below notations used Section \ref{section:method} for the conditional expectations given different part of the data.
\begin{definition}\label{def:cond-exp}
The conditional expectation for samples with index in set $\Scal$
conditionally on subset of the data $\Dscr$ is denoted as
\begin{equation*}
  \E_{i\in \Scal}\{ f(Y_i,\bX_i,\bS_i) \mid \Dscr\}, \;
  \Scal \subseteq \{1,\dots, n+N\},  \Dscr \subset \Lscr \cup \Uscr.
\end{equation*}
We denote the conditional expectation of unlabeled data given labelled data by $\E_{i>n}\{f(\bW_i) \mid \Lscr \}$ and the conditional probability  of new copy of data given current data by $\P\subnew\{f(\bW_i) \mid \Dscr \}$. With $\Lscr$ and $\Uscr$ partitioned into $\Kn$ folds indexed respectively by $\{\Ical_k, k = 1, ..., \Kn\}$ and $\{\Jcal_k, k = 1, ..., \Kn\}$, we denote the conditional expectation of fold-k labelled data and unlabeled data given the out-of-fold data respectively by
  \begin{gather*}
  \E_{i \in \Ical_k} \{ f(Y_i, \bX_i, \bS_i) \mid \Dscr^c_k \}\quad \mbox{and}\quad \E_{i \in \Jcal_k} \{ f(\bW_i) \mid \Dscr^c_k \}, \\ \mbox{where}\quad
  \Dscr^c_k = \{\bS_i, \bX_i, i \in \Jcal_k^c\} \cup \{Y_i, \bS_i, \bX_i , i \in \Ical_k^c\}.
  \end{gather*}
\end{definition}
\subsection{Proof of Theorem \ref{thm:consistency}}

Our proof shares the general steps with the the restricted strong convexity framework
laid down in \cite{NRWY2010TR}
while we have a delicate analysis of the symmetrized Bregman  divergence
to establish the improved rate of estimation under semi-supervised learning setting.
To bound $\bbetahat$ through the symmetrized Bregman  divergence $(\hat{\bbeta}-\bbeta_0)^\top \dell^\dag(\bbeta\subo;\hat{\bgamma})$, instead of directly applying the H\"{o}lder's bound, we first split it into two parts,
\begin{align}
  (\hat{\bbeta}-\bbeta_0)^\top \dell^\dag(\bbeta\subo;\hat{\bgamma})
  = &\underbrace{(\hat{\bbeta}-\bbeta_0)^\top \left[\dell^\dag(\bbeta\subo;\hat{\bgamma}) - \E\{\dell^\dag(\bbeta\subo;\hat{\bgamma})\mid \Lscr\}+\E\{\dell^\dag(\bbeta\subo;\bgamma\subo)\mid \Lscr\}
 \right]}_{\text{variance from unlabeled data}} \notag \\
 & +
 \underbrace{(\hat{\bbeta}-\bbeta_0)^\top\E\left\{\dell^\dag(\bbeta\subo;\hat{\bgamma}) -\dell^\dag(\bbeta\subo;\bgamma\subo)\mid \Lscr\right\}}_{\text{bias from $\hat{\bgamma}$}}\label{eq:thm1-decomp}
\end{align}
and discuss which part dominates the estimation error.  When the first variance term in \eqref{eq:thm1-decomp} is dominant, the bias from $\hat{\bgamma}$ becomes eligible. Then, we should recover the usual error bound for LASSO as if $\bgamma\subo$ is used. When the second bias term in \eqref{eq:thm1-decomp} is dominant, the  error bound of $\hat{\bbeta}$ can be controlled by the error bound of $\hat{\bgamma}$.
Combining the error bounds in the two cases, we obtain the oracle inequalities.
\begin{lemma}\label{lemma:beta-oracle}
On event
\begin{align*}
\Omega =
\Big\{&\ell\subPL(\bbeta\subo+\Delta ) - \ell\subPL(\bbeta\subo)
  - \Delta^\top\dell\subPL(\bbeta\subo) \\
   &\ge \crsc{1} \|\Delta\|_2\{\|\Delta\|_2-\crsc{2}\sqrt{\log(p)/N}\|\Delta\|_1\},  \forall \|\Delta\|_2 \le 1 \Big\},
\end{align*}
setting $\lamb \gtrsim \sqrt{\log(p)/N}$ such that
$$
\lamb\ge  3 \left\|\dell^\dag(\bbeta\subo;\hat{\bgamma}) - \E\{\dell^\dag(\bbeta\subo;\hat{\bgamma})\mid \Lscr\}+\E\{\dell^\dag(\bbeta\subo;\bgamma\subo)\mid \Lscr\}
 \right\|_\infty  + \crsc{1}\crsc{2}\sqrt{\frac{\log(p)}{N}},
$$
we have the oracle inequalities for estimation error of $\hat{\bbeta}$,
\begin{gather*}
  \|\hat{\bbeta}-\bbeta\subo\|_2 \le \max\left\{14\sqrt{\sbeta}\lamb/\crsc{1},
   (1-\rho) 7 M \sigma^2\submax \|\bgamma\subo-\hat{\bgamma}\|_2/\crsc{1}\right\}, \\
     \|\hat{\bbeta}-\bbeta\subo\|_1 \le \max\left\{84\sbeta\lamb/\crsc{1},
   (1-\rho)^2 21 M^2 \sigma^4\submax \|\bgamma\subo-\hat{\bgamma}\|_2^2/(\crsc{1}\lamb)\right\}.
\end{gather*}
The constants $\crsc{1},\crsc{2}$ are the restrictive strong convexity parameters
specified in Lemma \ref{lemma:rsc}.
\end{lemma}

We next prove the oracle inequalities. First, we note that by the definition of $\hat{\bbeta}$,
\begin{equation}\label{eq:loss-min}
   \ell^\dag(\hat{\bbeta};\hat{\bgamma}) + \lamb\|\hat{\bbeta}\|_1 \le \ell^\dag(\bbeta\subo;\hat{\bgamma}) + \lamb\|\bbeta\subo\|_1.
\end{equation}
Denote the standardized estimation error as $\bgd = (\hat{\bbeta}-\bbeta\subo)/\|\hat{\bbeta}-\bbeta\subo\|_2$.
Due to convexity of the loss function, we have for $t = \|\hat{\bbeta}-\bbeta\subo\|_2 \wedge 1$
\begin{equation}\label{eq:loss-min-dir}
   \ell^\dag(\bbeta\subo + t\bgd;\hat{\bgamma}) + \lamb\|\bbeta\subo + t\bgd\|_1 \le \ell^\dag(\bbeta\subo;\hat{\bgamma}) + \lamb\|\bbeta\subo\|_1.
\end{equation}
By the triangle inequality $\|\bbeta\subo\|_1 -  \|\bbeta\subo + t\bgd\|_1 \le t\|\bgd\|_1$, we have from \eqref{eq:loss-min-dir}
\begin{equation}\label{eq:loss-dir-triangle}
  \ell^\dag(\bbeta\subo + t\bgd;\hat{\bgamma})-\ell^\dag(\bbeta\subo;\hat{\bgamma}) \le t \lamb \|\bgd\|_1
\end{equation}
To apply the restricted strong convexity of the complete data loss \eqref{def:beta-dev-comp} established in Lemma \ref{lemma:rsc}, we show that the second order approximation error of the imputed loss is equivalent to that of the complete data loss,
\begin{align*}
  & \ell^\dag(\bbeta\subo + t\bgd;\hat{\bgamma})-\ell^\dag(\bbeta\subo;\hat{\bgamma}) - t\bgd^\top \dell^\dag(\bbeta\subo;\hat{\bgamma})
=  \ell\subPL(\bbeta\subo+t\bgd) - \ell\subPL(\bbeta\subo) - t\bgd^\top\dell\subPL(\bbeta\subo).
\end{align*}
Then by applying the restricted strong convexity event $\Omega$, we obtain
\begin{align}
  \ell^\dag(\bbeta\subo + t\bgd;\hat{\bgamma})-\ell^\dag(\bbeta\subo;\hat{\bgamma}) - t\bgd^\top \dell^\dag(\bbeta\subo;\hat{\bgamma}) 
  \ge t^2\crsc{1} -t\crsc{1}\crsc{2}\sqrt{\log(p)/N}\|\bgd\|_1. \label{eq:rsc-imputed}
\end{align}
Applying \eqref{eq:rsc-imputed} to \eqref{eq:loss-dir-triangle}, we have with large probability
\begin{equation*}
t \bgd^\top \dell^\dag(\bbeta\subo;\hat{\bgamma})
+ t^2\crsc{1} -t\crsc{1}\crsc{2}\sqrt{\log(p)/N}\|\bgd\|_1
   \le t \lamb\|\bgd\|_1
\end{equation*}
where $\|\bgd\|_2 = 1$ from definition.
Thus, we have reach
\begin{equation}\label{eq:thm1-bound}
 t\crsc{1}
   \le  \lamb\|\bgd\|_1 - \bgd^\top \dell^\dag(\bbeta\subo;\hat{\bgamma})+\crsc{1}\crsc{2}\sqrt{\log(p)/N}\|\bgd\|_1.
\end{equation}

Next, we analyze $\bgd^\top \dell^\dag(\bbeta\subo;\hat{\bgamma})$
by the decomposition
\begin{align}
  \left|\bgd^\top \dell^\dag(\bbeta\subo;\hat{\bgamma})\right|
= &\bgd^\top\left[\dell^\dag(\bbeta\subo;\hat{\bgamma}) - \E\{\dell^\dag(\bbeta\subo;\hat{\bgamma})\mid \Lscr\}+\E\{\dell^\dag(\bbeta\subo;\bgamma\subo)\mid \Lscr\}
 \right] \notag \\
 &+ \bgd^\top\E\left\{\dell^\dag(\bbeta\subo;\hat{\bgamma})-\dell^\dag(\bbeta\subo;\bgamma\subo)\mid \Lscr\right\} \notag \\
  \le & \|\bgd\|_1 \left\|\dell^\dag(\bbeta\subo;\hat{\bgamma}) - \E\{\dell^\dag(\bbeta\subo;\hat{\bgamma})\mid \Lscr\}+\E\{\dell^\dag(\bbeta\subo;\bgamma\subo)\mid \Lscr\}
 \right\|_\infty \notag \\
&  + (1-\rho)\left\|\E_{i > n}[\bX_i\{g(\bbeta\subo^\top \bX_i) - g(\hat{\bgamma}^\top \bW_i)\}\mid \Lscr]\right\|_2.
  \label{eq:thm1-decomp-score}
\end{align}
To establish the rate for $L_2$-norm of $\E\{\dell^\dag(\bbeta\subo;\hat{\bgamma})\mid \Lscr\}$, we note that
\begin{equation}
 \E_{i > n}[\bX_i\{g(\bbeta\subo^\top \bX_i) - g(\hat{\bgamma}^\top \bW_i)\}\mid \Lscr]
  =
 \E_{i > n}\left[\bX_i\left\{Y_i - g\left(\hat{\bgamma}^\top \bW_i\right)\right\}\mid \Lscr\right].\label{eq:bias-1}
\end{equation}
 By the characterization of $\bgamma\subo$ as in \eqref{def:bar-gamma}, we may rewrite \eqref{eq:bias-1}
as
\begin{align}
 \E_{i > n}\left[\bX_i\left\{Y_i - g\left(\hat{\bgamma}^\top \bW_i\right)\right\}\mid \Lscr\right]
  = & \E_{i > n}\left[g'(\bgamma_u^\top \bW_i)\bX_i\bW_i^\top\mid \Lscr\right] (\bgamma\subo-\hat{\bgamma}), \label{eq:bias-2}
\end{align}
where $\bgamma_u = u\hat{\bgamma}+(1-u)\bgamma\subo$ for some $u\in [0,1]$.
Under Assumptions \ref{assume:link} and \ref{assume:X-subG},
as well as the fact that $\bX_i$ is a sub-vector of $\bW_i$,
we have
\begin{align}
\left\|  \E\{\dell^\dag\left(\bbeta\subo;\hat{\bgamma}\right)\mid \Lscr \}\right\|_2 \le & \left\|\E_{i\in \Ucal}\left[g'(\bgamma_u^\top \bW_i)\bW_i^{\otimes 2}\mid \Lscr\right]\right\|_2  \|\bgamma\subo-\hat{\bgamma}\|_2 \notag \\
 \le &   M \left\| \E\left( \bW_i^{\otimes 2} \right)\right\|_2  \|\bgamma\subo-\hat{\bgamma}\|_2
  \le M \sigma^2\submax \|\bgamma\subo-\hat{\bgamma}\|_2,
  \label{eq:thm1-decomp-score}
\end{align}
where for any vector $\bx$, $\bx^{\otimes 2}=\bx\bx\trans$.
By the bound for (\ref{eq:thm1-decomp-score}) and the definition of $\lamb$, we have the bound from \eqref{eq:thm1-bound}
\begin{equation}
  t\crsc{1}
   \le  2\lamb\|\bgd\|_1 +   (1-\rho) M \sigma^2\submax \|\bgamma\subo-\hat{\bgamma}\|_2.  \label{eq:thm1-case}
\end{equation}
Hence, we can reach
an immediate bound for estimation error
from \eqref{eq:thm1-case} without considering the sparsity of $\bbeta\subo$.
We shall proceed to derive a sharper bound
that involves the sparsity of $\bbeta\subo$.
We separately analyze two cases.
\begin{equation*}
\mbox{\underline{\bf Case 1:} }\quad
 (1-\rho) M \sigma^2\submax \|\bgamma\subo-\hat{\bgamma}\|_2
  \ge  \|\bgd\|_1\lamb/3 .
\end{equation*}
In this case, the estimation error is dominated by $\hat{\bgamma}-\bgamma\subo$.
We simply have from \eqref{eq:thm1-case}
\begin{gather*}
  t\crsc{1}
   \le  7 (1-\rho) M \sigma^2\submax \|\bgamma\subo-\hat{\bgamma}\|_2, \\
   t\crsc{1} \|\bgd\|_1\lamb/3 \le 7 (1-\rho)^2 M^2 \sigma^4\submax \|\bgamma\subo-\hat{\bgamma}\|_2^2.
   \end{gather*}
Thus, we have
\begin{gather}\label{eq:thm1-rate1}
  \|\hat{\bbeta}-\bbeta\subo\|_2 \le (1-\rho) 7 M \sigma^2\submax \|\bgamma\subo-\hat{\bgamma}\|_2/\crsc{1}, \notag \\
   \|\hat{\bbeta}-\bbeta\subo\|_1 \le (1-\rho)^2 21 M^2 \sigma^4\submax \|\bgamma\subo-\hat{\bgamma}\|_2^2/(\crsc{1}\lamb).
\end{gather}
If case 1 does not hold, then instead
\begin{equation}\label{eq:thm1-case2}
\mbox{\underline{\bf Case 2:} }\quad
  (1-\rho)M \sigma^2\submax \|\bgamma\subo-\hat{\bgamma}\|_2
  \le  \|\bgd\|_1\lamb/3 .
\end{equation}
In this case, the estimation error is comparable to that when we have the true $\bgamma\subo$
for the imputation.
Thus, the sparsity of $\bbeta\subo$ may affect the estimation error.

Following the typical approach to establish the cone condition for $\bgd$,
we analyze the symmetrized Bregman's divergence,
\begin{equation}\label{eq:thm1-breg}
  (\hat{\bbeta} - \bbeta\subo)^\top  \left\{ \dell^\dag(\hat{\bbeta};\hat{\bgamma})
  - \dell^\dag(\bbeta\subo;\hat{\bgamma}) \right\}
  = \|\hat{\bbeta} - \bbeta\subo\|_2 \bgd^\top \left\{ \dell^\dag(\hat{\bbeta};\hat{\bgamma})
  - \dell^\dag(\bbeta\subo;\hat{\bgamma}) \right\}.
\end{equation}
Due to the convexity of the loss $\dell^\dag(\cdot;\hat{\bgamma})$
under Assumption \ref{assume:link},
the symmetrized Bregman's divergence \eqref{eq:thm1-breg} is nonnegative
through a mean-value theorem,
$$
   (\hat{\bbeta} - \bbeta\subo)^\top  \left\{ \dell^\dag(\hat{\bbeta};\hat{\bgamma})
  - \dell^\dag(\bbeta\subo;\hat{\bgamma}) \right\}
  \ge   \inf_{\bbeta \in \R^{p+1}} \frac{1}{N}\sum_{i > n} g'(\bbeta^\top \bX_i) \{ (\hat{\bbeta} - \bbeta\subo)^\top \bX_i\}^2
  \ge  0.
$$
Denote the indices set of nonzero coefficient in $\bbeta\subo$ as
$\nzb = \{j: \bbeta_{0,j} \neq 0\}$.
We denote the $\bgd_{\nzb}$ and $\bgd_{\nzbc}$
as the sub-vectors for $\bgd$ at positions in $\nzb$ and at positions not in $\nzb$,
respectively.
The solution $\hat{\bbeta}$ satisfies the KKT condition
\begin{equation*}
 \| \dell^\dag(\hat{\bbeta};\hat{\bgamma})\|_\infty \le \lamb, \;
  \dot{\ell}^\dag(\hat{\bbeta};\hat{\bgamma})_j = -\lamb\sgn(\hat{\bbeta}_j), \, j:\hat{\beta}_j \neq 0.
\end{equation*}
From the KKT condition and the definitions of $\bgd$ and $\nzb$, we have
\begin{equation}\label{eq:KKT-breg}
  \delta_j \dot{\ell}^\dag(\hat{\bbeta};\hat{\bgamma})_j \le |\delta_j| \lamb, \, j \in \nzb; \;
  \delta_j \dot{\ell}^\dag(\hat{\bbeta};\hat{\bgamma})_j = \frac{-\hat{\beta}_j \lamb\sgn(\hat{\beta}_j) }{\|\hat{\bbeta}-\bbeta\subo\|_2} = -\lamb|\delta_j|, \, j \in \nzbc.
\end{equation}
Applying the \eqref{eq:KKT-breg} to \eqref{eq:thm1-breg}, we have the upper bound,
\begin{align*}
  \bgd^\top   \left\{ \dell^\dag(\bbeta\subo;\hat{\bgamma})
  - \dell^\dag(\bbeta\subo;\hat{\bgamma}) \right\}
  = & \sum_{j \in \nzb} \delta_j \dot{\ell}^\dag(\hat{\bbeta};\hat{\bgamma})_j
  + \sum_{j \in \nzbc} \delta_j \dot{\ell}^\dag(\hat{\bbeta};\hat{\bgamma})_j
  + \bgd^\top \dell^\dag(\bbeta\subo;\hat{\bgamma}) \notag \\
  \le &  \lamb\sum_{j \in \nzb} |\delta_j|
  -  \lamb\sum_{j \in \nzbc} |\delta_j|
  + \bgd^\top \dell^\dag(\bbeta\subo;\hat{\bgamma}) \notag \\
  \le & \lamb\|\bgd_{\nzb}\|_1 - \lamb\|\bgd_{\nzbc}\|_1
  + \left|\bgd^\top \dell^\dag(\bbeta\subo;\hat{\bgamma})\right|.
\end{align*}
Then, we apply \eqref{eq:thm1-decomp-score}, the definition of $\lamb$ and \eqref{eq:thm1-case2},
\begin{align*}
  0 \le & \lamb\|\bgd_{\nzb}\|_1 - \lamb\|\bgd_{\nzbc}\|_1
  + \frac{2}{3} \lamb\|\bgd\|_1 \quad \mbox{and}\quad
 \lamb\|\bgd_{\nzbc}\|_1  \le  5 \lamb\|\bgd_{\nzb}\|_1.
\end{align*}
Therefore, we can bound the $L_1$ norm of $\bgd$ by the cone property,
\begin{equation}\label{eq:thm1-cone}
\|\bgd\|_1 \le 6 \lamb\|\bgd_{\nzb}\|_1
\le 6\sqrt{\sbeta}\|\bgd\|_2 = 6\sqrt{\sbeta}.
\end{equation}
We then apply the cone condition \eqref{eq:thm1-cone} and the case condition \eqref{eq:thm1-case2}
to the bound \eqref{eq:thm1-case},
\begin{equation*}
   t\crsc{1}
   \le  14\sqrt{\sbeta}\lamb, \; \quad \mbox{and}\quad
   t\crsc{1} \|\bgd\|_1 \le  84\sbeta\lamb
\end{equation*}
Thus, we obtain the rate for estimation error
\begin{equation}\label{eq:thm1-rate2}
  \|\hat{\bbeta}-\bbeta\subo\|_2 \le 14\sqrt{\sbeta}\lamb/\crsc{1}, \;
  \quad \mbox{and}\quad
  \|\hat{\bbeta}-\bbeta\subo\|_1 \le 84\sbeta\lamb/\crsc{1}.
\end{equation}

Since Case 1 and Case 2 are the complement of each other,
one of them must occur.
Thus, the bound of estimation error is controlled by the larger bound in
the two cases,
\begin{align*}
  \|\hat{\bbeta}-\bbeta\subo\|_2 & \le \max\left\{14\sqrt{\sbeta}\lamb/\crsc{1},
   (1-\rho) 7 M \sigma^2\submax \|\bgamma\subo-\hat{\bgamma}\|_2/\crsc{1}\right\}, \\
     \|\hat{\bbeta}-\bbeta\subo\|_1 & \le \max\left\{84\sbeta\lamb/\crsc{1},
   (1-\rho)^2 21 M^2 \sigma^4\submax \|\bgamma\subo-\hat{\bgamma}\|_2^2/(\crsc{1}\lamb)\right\},
\end{align*}
which is our oracle inequality in Lemma \ref{lemma:beta-oracle}.

\paragraph{Consistency}
We next show that the oracle inequality leads to the consistency
under dimension condition \eqref{assume:dim-est}. To show
\begin{align*}
   &  \left\|\dell^\dag(\bbeta\subo;\hat{\bgamma}) - \E\{\dell^\dag(\bbeta\subo;\hat{\bgamma})\mid \Lscr\}+\E\{\dell^\dag(\bbeta\subo;\bgamma\subo)\mid \Lscr\}
 \right\|_\infty \\
 \le & \left\|\dell^\dag(\bbeta\subo;\hat{\bgamma}) - \E\{\dell^\dag(\bbeta\subo;\hat{\bgamma})\mid \Lscr\}\right\|_\infty +\left\|\E\{\dell^\dag(\bbeta\subo;\bgamma\subo)\mid \Lscr\}
 \right\|_\infty  = O_p\left(\sqrt{\log(p)/N}\right),
\end{align*}
we express the term of interest as the sum of the following empirical processes
\begin{gather*}
\begin{aligned}
\dell^\dag(\bbeta\subo;\hat{\bgamma}) - \E\{\dell^\dag(\bbeta\subo;\hat{\bgamma})\mid \Lscr\}
= &\frac{1}{N} \sum_{i>n} \Big(\bX_i \{g(\bbeta\subo^\top \bX_i) - Y_i + Y_i - g(\hat{\bgamma}^\top \bW_i)\} \\ & - \E_{i>n}\left[\bX_i \{g(\bbeta\subo^\top \bX_i) - Y_i + Y_i - g(\hat{\bgamma}^\top \bW_i)\}\mid \Lscr\right]\Big),
\end{aligned} \\
\E\{\dell^\dag(\bbeta\subo;\bgamma\subo)\mid \Lscr\} =  \frac{1}{N} \sum_{i=1}^{n} \bX_i \{g(\bbeta\subo^\top \bX_i) - Y_i\}.
\end{gather*}
Under Assumption \ref{assume:Y-tail} and \ref{assume:Y-tail},
$\bX_i$ and $g(\bbeta\subo^\top \bX_i) - Y_i$ are sub-Gaussian.
According to Lemma \ref{lemma:est-gamma}, the event $\|\hat{\bgamma}-\bgamma\subo\|_2 \le 1$ occurs with large probability, on which
we have a bound for the sub-Gaussian norm of $Y_i - g(\hat{\bgamma}^\top \bW_i)$
by Lemma \ref{lemma:subg-g}.
\begin{equation}\label{eq:Yg-subg}
  \|Y_i - g(\hat{\bgamma}^\top \bW_i)\|_{\psi_2} \le \max\{2\nu_1,
  M \sqrt{2}\sigma\submax\}, \, i > n
\end{equation}
Thus, we obtain from \eqref{eq:Yg-subg} that $Y_i - g(\hat{\bgamma}^\top \bW_i)$ is sub-Gaussian
with large probability.
Thus by the properties of sub-Gaussian random variables in Lemma \ref{lemma:sub-add} and \ref{lemma:sub-mg},
we have established that the elements in the summands of $\dell^\dag(\bbeta\subo;\hat{\bgamma})$
are all sub-exponential random variables conditionally on the labelled data.
We apply the Bernstein's inequality (Lemma \ref{lemma:sub-bern})
conditionally on the labelled data to obtain
\begin{align*}
\left\|\dell^\dag(\bbeta\subo;\hat{\bgamma}) - \E\{\dell^\dag(\bbeta\subo;\hat{\bgamma})\mid \Lscr\}
 \right\|_\infty & = O_p\left(\sqrt{(1-\rho)\log(p)/N}\right), \\
  \left\|\E\{\dell^\dag(\bbeta\subo;\bgamma\subo)\mid \Lscr\}
 \right\|_\infty & = O_p\left(\sqrt{\rho\log(p)/N}\right).
\end{align*}
This establishes the order for $\lambda_{\bbeta}$,
\begin{equation}\label{eq:rate-lamb}
\lambda_{\bbeta}  \gtrsim \sqrt{(1-\rho)\log(p)/N}+\sqrt{\rho\log(p)/N}
\asymp \sqrt{\log(p)/N}.
\end{equation}

By Lemma \ref{lemma:rsc} from \cite{NRWY2010TR}, we have that
the probability of restricted strong convexity event converges to one,
$$
\P(\Omega) \ge 1- \crsc{3}e^{\crsc{4}N} \to 1.
$$
Setting $\lamb \asymp \sqrt{\log(p)/N}$ for optimal $L_2$ estimation, we achieve the stated conclusion
\begin{align*}
  \|\hat{\bbeta}-\bbeta\subo\|_2 & = O_p\left(\sqrt{\sbeta\log(p)/N}+ (1-\rho)\sqrt{\sgamma \log(p+q)/n}\right), \\
  \sqrt{\log(p)/N} \|\hat{\bbeta}-\bbeta\subo\|_1 & = O_p\left(\sbeta\log(p)/N+ (1-\rho)^2\frac{\sgamma \log(p+q)}{n}\right),
\end{align*}
by applying the rates from Lemma \ref{lemma:est-gamma} and \eqref{eq:rate-lamb}.
For optimal $L_1$ estimation, we set a larger penalty
 $
 \lamb' \asymp \sqrt{\log(p)/N} \vee \sqrt{\sgamma \log(p+q)/(\sbeta n)} \gtrsim \lamb
 $
to achieve
$$
\|\hat{\bbeta}-\bbeta\subo\|_1 = O_p\left(\sbeta\sqrt{\log(p)/N}+ (1-\rho)^2\sqrt{\frac{\sgamma \sbeta\log(p+q)}{n}}\right).
$$

\subsection{Proof of Corollary \ref{cor:individual}}

Under Assumption \ref{assume:link}, we have
\begin{equation}\label{eq:cor2-g}
|g(\hat{\bbeta}\trans\bx\subnew) - g(\bbeta\subo\trans\bx\subnew)|
\le M |(\hat{\bbeta} - \bbeta\subo)\trans \bx\subnew |.
\end{equation}
Since $\bx\subnew$ satisfies Assumption \ref{assume:X-subG},
we have
$$
\|(\hat{\bbeta} - \bbeta\subo)\trans \bx\subnew \|_{\psi_2}
\le \|\hat{\bbeta} - \bbeta\subo\|_2 \sigma\submax/\sqrt{2}.
$$
The tail distribution is regulated by the sub-Gaussian norm by Lemma \ref{lemma:sub-tail},
\begin{equation}\label{eq:cor2-tail}
  \P\subnew\left(|(\hat{\bbeta} - \bbeta\subo)\trans \bx\subnew | \ge  t \mid \Dscr \right)
  \le 2 \exp \left(-t\sqrt{2}/\left\{\|\hat{\bbeta} - \bbeta\subo\|_2 \sigma\submax\right\}\right).
\end{equation}
Combining \eqref{eq:cor2-g} and \eqref{eq:cor2-tail}, we obtain
$$
 \P\subnew\left(|g(\hat{\bbeta}\trans\bx\subnew) - g(\bbeta\subo\trans\bx\subnew)|
 \ge t M \|\hat{\bbeta} - \bbeta\subo\|_2 \sigma\submax/\sqrt{2} \mid \Dscr\right)
 \le 2e^{-t}.
$$
Thus,
$$
\left|g(\hat{\bbeta}\trans\bx\subnew) - g(\bbeta\subo\trans\bx\subnew)\right|  = O_p\left(\|\hat{\bbeta} - \bbeta\subo\|_2\right).
$$

\subsection{Proof of Theorem \ref{thm:inference}}\label{asection:proof-inf}

\def\half{\frac{1}{2}}
\def\nnhalf{n^{-\half}}

Our proof is organized in five parts.
In Part 1, we establish the consistency of the cross-fitting estimator for precision matrix,
namely $\|\hat{\bu}\supk-\bu\subo\|_2 = o_p(1)$ with $\hat{\bu}\supk$ and $\bu\subo$ defined in \eqref{def:uk} and \eqref{def:u0}, respectively.
In Part 2, we show that the debiased estimator can be approximated by
the empirical process
\begin{align*}
\sqrt{n}\left(\widehat{\bx\substd\trans \bbeta} - \bx\substd\trans \bbeta\subo\right)
= & -(1-\rho)\sqrt{n}\bu\subo\trans\impscore(\bgamma_0)- \sqrt{n}\bu\subo\trans\dell^\dag(\bbeta\subo; \bgamma\subo)
+ o_p(1) \\
= & -\nnhalf \left[\sum_{i=1}^{n}\bu\subo\trans \bX_i \{(1-\rho)\cdot g(\bgamma\subo\trans \bW_i)+ \rho \cdot g(\bbeta\subo\trans \bX_i)-Y_i\}\right.\\
& \qquad \qquad \left.  +
\rho\sum_{i>n} \bX_i \{g(\bbeta\subo\trans \bX_i)-g(\bgamma\subo\trans \bW_i)\}
\right] + o_p(1)
\end{align*}
As long as the asymptotic variance $\VSAS$ defined in \eqref{def:sig0} is bounded and bounded away from zero,
we have the asymptotic normality of the leading term from the Central Limit Theorem
\begin{align*}
-\nnhalf\VSAS^{-1/2}\left[\sum_{i=1}^{n}\right.&\bu\subo\trans \bX_i \{(1-\rho)\cdot g(\bgamma\subo\trans \bW_i)+ \rho \cdot g(\bbeta\subo\trans \bX_i)-Y_i\} \\
& \left. \hspace{1in} +
\rho\sum_{i>n} \bX_i \{g(\bbeta\subo\trans \bX_i)-g(\bgamma\subo\trans \bW_i)\}
\right] \leadsto N(0,1).
\end{align*}
In Part 3, we deal with the asymptotic variance $\VSAS$ and the consistency of the variance estimator $\hVSAS $ defined in \eqref{def:var-cf}.
In Part 4, we reach the conclusion of the theorem based on,
\begin{align}
 & (1-\rho)\|\hat{\bgamma}\supk-\bgamma\subo\|_2 +\|\hat{\bu}\supk-\bu\subo\|_2\notag \\
 &+  \sqrt{n} \|\hat{\bbeta}\supk-\bbeta\subo\|_2\left(\|\hat{\bbeta}\supk-\bbeta\subo\|_2 +\|\hat{\bu}\supk-\bu\subo\|_2\right)= o_p(1),  \label{assume:est-rate}
\end{align}
for all $1\leq k\leq \Kn$.
Following Part 4, we show in Part 5 that \eqref{assume:dim-inf} implies \eqref{assume:est-rate}.

\subsubsection*{Part 1: Consistency of estimated precision matrix}

The definitions of $\bu\subo$ and $\hat{\bu}\supk$ are given in \eqref{def:u0} and
\eqref{def:uk}.
In this part, we show
\begin{align*}
\|\hat{\bu}\supk-\bu\subo\|_2 =  & O_p\left(\sqrt{(\su+\sbeta)\log(p)/(N-N_k)} + (1-\rho)\sqrt{\sgamma\log(p+q)/(n-n_k)}\right) \\
= & O_p\left(\sqrt{(\su+\sbeta)\log(p)/N} + (1-\rho)\sqrt{\sgamma\log(p+q)/n}\right).
\end{align*}
Since we set the number of folds $\Kn \le 10$ to be finite,
the estimation rate applies for $\hat{\bu}\supk$ for all $k= 1, \dots, \Kn$.

We denote the components in the quadratic loss function of \eqref{def:uk} and their derivatives as
\begin{gather}
    m\supkk(\bu;\bbeta) = \frac{1}{N_{k'}}\sum_{i \in \Ical_{k'}\cup\Jcal_{k'}} \frac{1}{2}g'(\bbeta\trans\bX_i)(\bX_i\trans\bu)^2- \bu\trans\bx\substd, \notag \\
    \dm\supkk(\bu;\bbeta) = \frac{\partial}{\partial \bu}  m\supkk(\bu;\bbeta), \;
     \ddm\supkk(\bbeta) =\frac{\partial}{\partial \bu}\dm\supkk(\bu;\bbeta) \label{def:m}
\end{gather}
for $k' \in \{1,\dots, \Kn\} \setminus \{k\}$.
We may express \eqref{def:uk} as
$$
\hat{\bu}\supk =   \argmin_{\bu \in \R^p} \sum_{k'\neq k}\frac{N_{k'}}{N-N_k} m\supkk\left(\bu;\hat{\bbeta}\supkk\right) + \lamu\|\bu\|_1,
$$
Similar to the proof of Theorem \ref{thm:consistency}, we establish the estimation rate for $\hat{\bu}$
through an oracle inequality,
\begin{lemma}\label{lemma:u-oracle}
Under Assumptions \ref{assume:model}, \ref{assume:X}, we establish
On event
\begin{align*}
\Omega\supk =
\bigcap_{k'\neq k}\Big\{\Delta^\top \ddm\supkk\left(\hat{\bbeta}\supkk\right) \Delta \ge \crsc{1}^* \|\Delta\|_2\{\|\Delta\|_2-\crsc{2}^*\sqrt{\log(p)/N_{k'}}\|\Delta\|_1\},  \forall \|\Delta\|_2 \le 1 \Big\},
\end{align*}
setting $\lamu \asymp \sqrt{\log(p)/N}$ such that
\begin{align*}
    \lamu\ge  3 \sum_{k'\neq k} \frac{N_{k'}}{N-N_k}\Bigg\{&\left\|\dm\supkk\left(\bu\subo;\hat{\bbeta}\supkk\right) - \E\left\{\dm\supkk\left(\bu\subo;\hat{\bbeta}\supkk\right)\mid \Dscr_{k'}^c\right\}
 \right\|_\infty  \\
& + \crsc{1}^*\crsc{2}^*\sqrt{\log(p)/N_{k'}}\Bigg\},
\end{align*}
we have the oracle inequality for estimation error of $\hat{\bbeta}$,
\begin{gather*}
  \|\hat{\bu}\supk-\bu\subo\|_2 \le \max\left\{14\sqrt{\su}\lamu/\crsc{1}^*,
   7 M \sigma\submax^3  \|\bu\subo\|_2 \sup_{k'\neq k}\left\|\bbeta\subo-\hat{\bbeta}\supkk\right\|_2/\crsc{1}^*\right\}, \\
     \|\hat{\bu}\supk-\bu\subo\|_1 \le \max\left\{84\su\lamu/\crsc{1}^*,
   21 M \sigma\submax^6  \|\bu\subo\|_2^2 \sup_{k'\neq k}\left\|\bbeta\subo-\hat{\bbeta}\supkk\right\|_2^2/(\crsc{1}^*\lamu)\right\}.
\end{gather*}
The constants $\crsc{1}^*,\crsc{2}^*$ are the restrictive strong convexity parameters
specified in Lemma \ref{lemma:u-rsc}.
\end{lemma}
The proof of Lemma \ref{lemma:u-oracle} repeats the proof of the oracle inequality for Theorem \ref{thm:consistency},
so we put the detail to Section \ref{asection:aux}.

To use Lemma \ref{lemma:u-oracle} for the estimation rate of $\hat{\bu}$,
we only need to verify two conditions.
First, the event $\Omega\supk$ occurs with probability tending to one.
Second, the oracle choice of $\lamu$ is of order $\sqrt{\log(p)/N}$.

Repeating Theorem \ref{thm:consistency} for each $\hat{\bbeta}\supkk$, we have under \eqref{assume:dim-inf}
$$
\left\|\hat{\bbeta}\supkk - \bbeta\subo\right\|_2 = o_p(1).
$$
Then by Lemma \ref{lemma:u-rsc}, the sets whose intersection forms $\Omega\supk$ each occurs with probability tending to one.
Since we set the number of fold finite $\Kn \le 10$,
we can take union bound to obtain that $\Omega\supk$ occurs with probability tending to one.

We may write
\begin{align}
& \dm\supkk\left(\bu\subo;\hat{\bbeta}\supkk\right) - \E\left\{\dm\supk\left(\bu\subo;\hat{\bbeta}\supkk\right)\mid \Dscr_{k'}^c\right\} \notag \\
= & \frac{1}{N_{k'}} \sum_{i \in \Ical_{k'}\cup\Jcal_{k'}} g'(\hat{\bbeta}\supkkt\bX_i) \bX_i\bX_i\trans\bu\subo
- \E_{i \in \Ical_{k'}\cup\Jcal_{k'}}\left\{g'(\hat{\bbeta}\supkkt\bX_i) \bX_i\bX_i\trans\bu\subo \mid \Dscr_{k'}^c\right\}. \label{eq:u-score}
\end{align}
Each element in \eqref{eq:u-score} is an empirical process.
Under Assumptions \ref{assume:link} and \ref{assume:X-subG},
we can show that
each summand is a sub-exponential random variable by Lemma \ref{lemma:sub-mb}, \ref{lemma:sub-mg},
$$
    \left\|g'(\hat{\bbeta}\supkkt\bX_i) X_{i,j}\bX_i\trans\bu\subo\right\|_{\psi_1}
    \le M \|X_{i,j}\|_{\psi_2}\|\bX_i\trans\bu\subo\|_{\psi_2}
    \le M \sigma\submax^2 \|\bu\subo\|_2/2.
$$
Hence, we can apply the Bernstein's inequality to show that
$$
\left\|\dm\supkk\left(\bu\subo;\hat{\bbeta}\supkk\right) - \E\left\{\dm\supkk\left(\bu\subo;\hat{\bbeta}\supkk\right)\mid \Dscr_{k'}^c\right\}
 \right\|_\infty  = O_p\left(\sqrt{\log(p)/N_{k'}}\right).
$$
Using the fact that $N_{k'} \asymp N$,
we obtain that the oracle $\lamu$ is of order $O_p\left(\sqrt{\log(p)/N}\right)$.

Therefore, we can apply Lemma \ref{lemma:u-oracle} to obtain
\begin{align*}
\|\hat{\bu}\supk - \bu\subo\|_2 = & O_p\left(\sqrt{\su\log(p)/N}+\sup_{k'\neq k}\left\|\hat{\bbeta}\supkk-\bbeta\subo\right\|_2 \right) \\
= & O_p\left(\sqrt{(\sbeta+\su)\log(p)/N} + (1-\rho)\sqrt{\sgamma \log(p+q)/n}\right).
\end{align*}

\subsubsection*{Part 2: Asymptotic approximation}

Under Assumption \ref{assume:sigmin-X}, we also have the tightness of $\|\hat{\bu}\supk\|_2$
from the bound of $\|\bu\subo\|_2$
\begin{equation}\label{eq:u0b-ukt}
\|\bu\subo\|_2 \le \|\Sigb^{-1}\|_2 \|\bx\substd\|_2 \le \sigma\submin^{-2}, \;
\|\hat{\bu}\supk\|_2 \le \|\bu\subo\|_2+  \|\hat{\bu}\supk - \bu\subo\|_2 = O_p(1).
\end{equation}

Define the scores of in-fold data as
\begin{gather}
  \dlkdag(\bbeta; \bgamma) = \frac{1}{N_k}\left[\sum_{i\in\Jcal_k}\bX_i\{g(\bbeta\trans \bX_i)-g(\bgamma\trans\bW_i)\} + \sum_{i\in\mathcal{I}_k}\bX_i \{g(\bbeta\trans \bX_i)-Y_i\}\right], \notag \\
  \dlkimp(\bgamma) = \frac{1}{n_k}\sum_{i\in\Ical_k}\bX_i \{g(\bgamma\trans\bW_i)-Y_i\}.
  \label{def:lk}
\end{gather}
Since $\widehat{\bx\substd\trans \bbeta}$ is the average over $\Kn$ (at most 10) cross-fitted estimators,
it suffices to study one of the cross-fitted estimators,
\begin{equation}
\widehat{\bx\substd\trans \bbeta}\supk  =  \bx\substd\trans \hat{\bbeta}\supk
-\dlkdag(\hat{\bbeta}\supk; \hat{\bgamma}\supk) - (1-\rho)\dlkimp(\hat{\bgamma}\supk), \;
\widehat{\bx\substd\trans \bbeta} = \frac{1}{\Kn}\sum_{k=1}^{\Kn} \widehat{\bx\substd\trans \bbeta}\supk.
\end{equation}

We denote the expected Hessian matrices of losses in \eqref{def:lk} as
\begin{gather}
  \Hess (\bbeta)
  = \E\left\{ g'(\bbeta\trans\bX_i)\bX_i\bX_i\trans\right\}, \; \Sigb = \Hess (\bbeta\subo),  \notag \\
  \Hessg (\bgamma)
  = \E\left\{g'(\bgamma\trans\bW_i)\bX_i\bW_i\trans\right\}, \; \Sigg = \Hessg (\bgamma\subo). \label{def:Hessg}
\end{gather}
Our analysis of the approximation error is based on the first order Mean Value Theorem
identity,
\begin{align}
 & \E\{\dlkdag(\hat{\bbeta}\supk; \hat{\bgamma}\supk)\mid \Dscr_k^c\}
  + (1-\rho) \E\{\dlkimp(\hat{\bgamma}\supk)\mid \Dscr_k^c\} \notag \\
   = & \underbrace{\E\{\dlkdag(\bbeta\subo; \bgamma\subo)\}}_{=0}
  + \Hess (\tilde{\bbeta}) \{\hat{\bbeta}\supk-\bbeta\subo\}  - (1-\rho)\Hessg(\tilde{\bgamma}) \{\hat{\bgamma}\supk- \bgamma\subo\} \notag \\
  & + (1-\rho)\underbrace{\E\{\dlkimp(\bgamma\subo)\}}_{=0}
  + (1-\rho) \Hessg(\tilde{\bgamma}) \{\hat{\bgamma}\supk- \bgamma\subo\} \notag \\
  = & \Hess (\tilde{\bbeta}) \{\hat{\bbeta}\supk-\bbeta\subo\} \label{eq:MVT}
\end{align}
for some $\tilde{\bbeta}$ on the path from $\bar{\bbeta}(\hat{\bgamma})$ to $\bbeta\subo$
  and some $\tilde{\bgamma}$ on the path from $\hat{\bgamma}$ to $\bgamma\subo$.
The conditional expectation notation is declared at Definition \ref{def:cond-exp}.
Based on \eqref{eq:MVT}, we analyze the approximation error for $\sqrt{n}\left(\widehat{\bx\substd\trans \bbeta}\supk - \bx\substd\trans \bbeta\subo\right)$ through the following decomposition,
\begin{align}
   & \sqrt{n}\left(\widehat{\bx\substd\trans \bbeta}\supk - \sqrt{n}\bx\substd\trans \bbeta\subo\right)
   + \sqrt{n}\bu\subo\trans\dlkdag(\bbeta\subo; \bgamma\subo) + \sqrt{n}(1-\rho)\bu\subo\trans\dlkimp(\bgamma\subo) \notag \\
   = & \sqrt{n}\bx\substd\trans(\hat{\bbeta}\supk - \bbeta_0)
   - \sqrt{n}\hat{\bu}\supkt\dlkdag(\hat{\bbeta}\supk; \hat{\bgamma}\supk) - \sqrt{n}(1-\rho)\hat{\bu}\supkt\dlkimp(\hat{\bgamma}\supk) \notag \\
   & + \sqrt{n}\bu\subo\trans\dlkdag(\bbeta\subo; \bgamma\subo) + (1-\rho)\bu\subo\trans\dlkimp(\bgamma\subo)
   + \sqrt{n}\hat{\bu}\supkt\E\{\dlkdag(\hat{\bbeta}\supk; \hat{\bgamma}\supk)\mid \Dscr_k^c\}\notag\\
 & + \sqrt{n}(1-\rho) \hat{\bu}\supkt\E\{\dlkimp(\hat{\bgamma}\supk)\mid \Dscr_k^c\}
   - \sqrt{n}\hat{\bu}\supkt\Hess (\tilde{\bbeta}) \{\hat{\bbeta}\supk-\bbeta\subo\} \notag \\
   = & \underbrace{\sqrt{n}\left\{\bx\substd - \Hess(\tilde{\bbeta})\bu\subo\right\}\trans(\hat{\bbeta}\supk - \bbeta\subo)}_{T_1}
   + \underbrace{\sqrt{n}\left(\bu\subo-\hat{\bu}\supk\right)\trans\Hess(\tilde{\bbeta})(\hat{\bbeta}\supk - \bbeta\subo)}_{T_2} \notag \\
   & + \underbrace{\sqrt{n} \hat{\bu}\supkt\left[ \E\{\dlkdag(\hat{\bbeta}\supk; \hat{\bgamma}\supk)\mid \Dscr_k^c\}
     - \left\{\dlkdag(\hat{\bbeta}\supk; \hat{\bgamma}\supk) - \dlkdag(\bbeta\subo; \bgamma\subo)\right\}\right]}_{T_3} \notag \\
   & + \underbrace{\sqrt{n} (1-\rho) \hat{\bu}\supkt\left[ \E\{\dlkimp(\hat{\bgamma}\supk)\mid \Dscr_k^c\}
     - \left\{\dlkimp(\hat{\bgamma}\supk) - \dlkimp(\bgamma\subo)\right\}\right]}_{T_4} \notag \\
   & + \underbrace{\sqrt{n}\left(\bu\subo-\hat{\bu}\supk\right)\trans\left\{\dlkdag(\bbeta\subo; \bgamma\subo)+(1-\rho)\dlkimp(\bgamma\subo)\right\}}_{T_5}
\end{align}

   Here we state the rates for $T_1$-$T_5$,
\begin{gather*}
  T_1 =  O_p\left(\sqrt{n}\|\hat{\bbeta}\supk-\bbeta\subo\|_2^2 \right), \;
  T_2 =  O_p\left(\sqrt{n}\|\hat{\bu}\supk-\bu\subo\|_2\|\hat{\bbeta}\supk-\bbeta\subo\|_2 \right), \\
  T_3 =  O_p\left(\rho\|\hat{\bbeta}\supk-\bbeta\subo\|_2+\sqrt{\rho(1-\rho)}\|\hat{\bgamma}\supk-\bgamma\subo\|_2\right), \\
  T_4 =  O_p\left((1-\rho)\|\hat{\bgamma}\supk-\bgamma\subo\|_2\right), \;
  T_5 =  O_p \left(\|\hat{\bu}\supk-\bu\subo\|_2\right).
\end{gather*}
With the assumed estimation rate in \eqref{assume:est-rate},
we have
$$
 T_1 + T_2 + T_3 + T_4 + T_5 = o_p(1).
$$
Thus, we have shown
\begin{align*}
\sqrt{n}\left(\widehat{\bx\substd\trans \bbeta}-\bx\substd\trans \bbeta\subo\right)
= &\frac{1}{\Kn}\sum_{k=1}^{\Kn}\sqrt{n}\left(\widehat{\bx\substd\trans \bbeta}-\bx\substd\trans \bbeta\subo\right) \\
= &\frac{1}{\Kn}\sum_{k=1}^{\Kn} -\sqrt{n}\bu\subo\trans\dlkdag(\bbeta\subo; \bgamma\subo) - \sqrt{n}(1-\rho)\bu\subo\trans\dlkimp(\bgamma\subo) + o_p(1)\\
= & -\sqrt{n}\bu\subo\trans\dell^\dag(\bbeta\subo; \bgamma\subo) - \sqrt{n}(1-\rho)\bu\subo\trans\impscore(\bgamma\subo) + o_p(1).
\end{align*}
Using the indicator $R_i = I(i \le n)$,
we can alternatively write
\begin{equation}\label{eq:eff-approx}
\widehat{\bx\substd\trans \bbeta}-\bx\substd\trans \bbeta\subo
= \frac{1}{N}\sum_{i=1}^N  \frac{R_i}{\rho} \bu\subo\trans\bX_i\{Y_i-
g(\bgamma\subo\trans\bW_i)\} - \bu\subo\trans\bX_i\{
g(\bgamma\subo\trans\bW_i)-g(\bbeta\subo\trans\bX_i)\}
+ o_p\left((\rho n)^{-1/2}\right).
\end{equation}
We provide the details of $T_1$-$T_5$ in Section \ref{asection:inf-remainders}.

\subsubsection*{Part 3: Variance estimation}
Finally, we show that asymptotic variance $\VSAS$ defined in \eqref{def:sig0}
is bounded from infinity and zero with the consistent estimator $\hVSAS$ defined in \eqref{def:var-cf}.

By the Cauchy-Schwartz inequality,
we have a bound for the variance
\begin{align*}
  \VSAS = & \E\left[(\bu\subo\trans \bX_i)^2 \{(1-\rho)\cdot g(\bgamma\subo\trans \bW_i)+\rho \cdot g(\bbeta\subo\trans\bX_i)-Y_i\}^2\right] \\
  &+ \rho(1-\rho) \E[(\bu\subo\trans \bX_i)^2 \{g(\bgamma\subo\trans \bW_i)- g(\bbeta\subo\trans\bX_i)\}^2] \notag \\
  \le & \sqrt{\E\left[(\bu\subo\trans \bX_i)^4\right] \E\left[\{(1-\rho)\cdot g(\bgamma\subo\trans \bW_i)+\rho \cdot g(\bbeta\subo\trans\bX_i)-Y_i\}^4\right]} \\
&  + \rho(1-\rho) \sqrt{\E\left[(\bu\subo\trans \bX_i)^4\right] \E\left[\{g(\bgamma\subo\trans \bW_i)- g(\bbeta\subo\trans\bX_i)\}^4\right]}
  \end{align*}
Under Assumptions \ref{assume:Y-tail}, \ref{assume:X-subG},
we have the sub-Gaussian and sub-exponential variables
\begin{gather*}
\|\bu\subo\trans \bX_i\|_{\psi_2} \le \|\bu\subo\|_2\sigma\submax/\sqrt{2} \le \sigma\submin^{-2}\sigma\submax/\sqrt{2}, \\
\|(1-\rho)\cdot g(\bgamma\subo\trans \bW_i)+\rho \cdot g(\bbeta\subo\trans\bX_i)-Y_i\|_{\psi_2}
\le  2 (\nu_1 \vee \nu_2), \\
\|g(\bgamma\subo\trans \bW_i)- g(\bbeta\subo\trans\bX_i)\|_{\psi_2}
\le 2 \{\| g(\bgamma\subo\trans \bW_i)-Y_i\|_{\psi_2}\vee \|g(\bbeta\subo\trans\bX_i)-Y_i\|_{\psi_2} \}
\le  2 (\nu_1 \vee \nu_2).
\end{gather*}
By the bound for the moments of sub-Gaussian and sub-exponential random variables stated in Lemma \ref{lemma:sub-moment}, we have
$$
\VSAS \le  8\sqrt{2}\sigma\submin^{-4}\sigma\submax^2(\nu_1\vee \nu_2)^2.
$$
Under Assumptions \ref{assume:link}, \ref{assume:X-subG}, \ref{assume:sigmin-X}
and \ref{assume:nondeg_var},
we have a lower bound for $\VSAS$,
\begin{align*}
\VSAS \ge & \bu\subo\trans \E[\bX_i\bX_i\trans \{(1-\rho)\cdot g(\bgamma\subo\trans \bW_i)+\rho \cdot g(\bbeta\subo\trans\bX_i)-Y_i\}^2] \bu\subo \\
 \ge & \|\bu\subo\|_2^2 \sigma\submin^4 \\
\ge & M^{-2}\sigma\submax^{-4}\sigma\submin^4\nu_3 \|\bx\substd\|_2^2,
\end{align*}
which is bounded away from zero.

We analyze the estimation error of variance $\hVSAS  - \VSAS$ through
the decomposition,
\begin{align*}
   & \quad \hVSAS  - \VSAS \\
& \left.\begin{aligned}
= & \sum_{k=1}^{\Kn}\frac{n_k}{n}\Bigg( \frac{1}{n_k}\sum_{i\in\Ical_k} (\hat{\bu}\supkt  \bX_i)^2 \{(1-\rho)\cdot g(\hat{\bgamma}\supkt  \bW_i)+\rho \cdot g(\hat{\bbeta}\supkt \bX_i)-Y_i\}^2 \\
& \qquad -\E_{i\in\Ical_k}\left[(\hat{\bu}\supkt  \bX_i)^2 \{(1-\rho)\cdot g(\hat{\bgamma}\supkt  \bW_i)+\rho \cdot g(\hat{\bbeta}\supkt \bX_i)-Y_i\}^2\mid \Dscr_k^c\right]\Bigg)
\end{aligned}\right\}{\scriptstyle T'_1}\\
& \left.\begin{aligned}
&+\sum_{k=1}^{\Kn}\frac{n_k}{n}\bigg( \E_{i\in\Ical_k}\left[(\hat{\bu}\supkt  \bX_i)^2 \{(1-\rho)\cdot g(\hat{\bgamma}\supkt  \bW_i)+\rho \cdot g(\hat{\bbeta}\supkt \bX_i)-Y_i\}^2\mid \Dscr_k^c\right] \\
& \qquad -\E\left[(\bu\subo\trans \bX_i)^2 \{(1-\rho)\cdot g(\bgamma\subo\trans  \bW_i)+\rho \cdot g(\bbeta\subo\trans \bX_i)-Y_i\}^2\right]\bigg)
\end{aligned}\right\} {\scriptstyle T'_2} \\
& \left.\begin{aligned}
&+ \rho(1-\rho)\sum_{k=1}^{\Kn}\frac{N_k-n_k}{N-n}\Bigg( \frac{\Kn}{N_k-n_k}\sum_{i\in\Jcal_k} (\hat{\bu}\supkt  \bX_i)^2 \{g(\hat{\bbeta}\supkt \bX_i)-g(\hat{\bgamma}\supkt  \bW_i)\}^2 \\
& \qquad -\E_{i\in\Ical_k}\left[(\hat{\bu}\supkt  \bX_i)^2 \{g(\hat{\bbeta}\supkt \bX_i)-g(\hat{\bgamma}\supkt  \bW_i)\}^2\mid \Dscr_k^c\right]\Bigg)
\end{aligned}\right\}{\scriptstyle T'_3}\\
& \left.\begin{aligned}
&+ \rho(1-\rho)\sum_{k=1}^{\Kn}\frac{N_k-n_k}{N-n}\bigg( \E_{i\in\Ical_k}\left[(\hat{\bu}\supkt  \bX_i)^2 \{g(\hat{\bbeta}\supkt \bX_i)-g(\hat{\bgamma}\supkt  \bW_i)\}^2\mid \Dscr_k^c\right] \\
& \qquad -\E\left[(\bu\subo\trans \bX_i)^2 \{ g(\bbeta\subo\trans \bX_i)-g(\bgamma\subo\trans  \bW_i)\}^2\right]\bigg)
\end{aligned}\right\} {\scriptstyle T'_4} \\
\end{align*}

   Here we state the rates for $T'_1$-$T'_4$,
\begin{gather*}
 T'_1 =  O_p\left( n^{-1/2} \right), \;
  T'_2 =  O_p\left(\|\hat{\bu}-\bu\subo\|_2+(1-\rho)\|\hat{\bgamma}-\bgamma\subo\|_2+\rho\|\hat{\bbeta} - \bbeta\subo\|_2 \right), \\
  T'_3 =   O_p\left(\rho(1-\rho) N^{-1/2} \right), \;
  T'_4 =  O_p \left(\rho(1-\rho)\left\{\|\hat{\bu}-\bu\subo\|_2+\|\hat{\bgamma}-\bgamma\subo\|_2+\|\hat{\bbeta} - \bbeta\subo\|_2\right\} \right).
\end{gather*}
With the assumed estimation rate in \eqref{assume:est-rate},
we have
$$
 T'_1 + T'_2 + T'_3 + T'_4  = o_p(1).
$$
We provide the details of $T'_1$-$T'_4$ in Section \ref{asection:inf-remainders}.

\subsubsection*{Part 4: Conclusion with estimation rates}

From the approximation in Part 2 and the boundedness and non-degeneracy of $\VSAS$ in Part 3, we have shown the asymptotic normality of the cross-fitted debiased estimator
$$
\sqrt{n}\VSAS^{-1/2}\left(\widehat{\bx\substd\trans \bbeta}  - \bx\substd\trans \bbeta\subo\right)
\leadsto N(0,1).
$$
Together with the consistency of $\hVSAS$ in Part 3,
we have
$$
\sqrt{n}\hVSAS^{-1/2}\left(\widehat{\bx\substd\trans \bbeta}  - \bx\substd\trans \bbeta\subo\right)
\leadsto N(0,1).
$$

\subsubsection*{Part 5: Sufficient dimension condition}

We have established the rate of estimation for $\hat{\bgamma}$, $\hat{\bbeta}$ and $\hat{\bu}$
from Lemma \ref{lemma:est-gamma}, Theorem \ref{thm:consistency} and Part 4 of this proof above.
Since we only keep one fold of the data away for the cross-fitted estimators,
they follow the same rates of estimation,
\begin{gather*}
\|\hat{\bgamma}\supk-\bgamma\subo\|_2 = O_p\left(\sqrt{\sgamma \log(p+q)/n}\right), \\
\|\hat{\bbeta}\supk-\bbeta\subo\|_2 = O_p\left(\sqrt{\sbeta\log(p)/N} + (1-\rho)\sqrt{\sgamma \log(p+q)/n}\right), \\
\|\hat{\bu}\supk-\bu\subo\|_2 = O_p\left(\sqrt{(\sbeta+\su)\log(p)/N} + (1-\rho)\sqrt{\sgamma \log(p+q)/n}\right).
\end{gather*}
Applying the rates of estimation,
we show dimension assumption \eqref{assume:dim-inf}
is sufficient for \eqref{assume:est-rate}.

\subsection{Efficiency of SAS Inference}

\subsubsection*{Relative Efficiency to Supervised Learning}

\begin{proof}[Proof of Proposition \ref{prop:eff-vs-super}]
We prove the Proposition by direct calculation
\begin{align*}
 & \VSL - \VSAS \\
 = & \E[(\bu\subo^\top \bX_i)^2 \{Y-g(\bbeta\subo\trans \bX_i)\}^2]
- \E[(\bu\subo^\top \bX_i)^2 \{Y- (1-\rho) \cdot \E(Y|\bS_i,\bX_i) \rho \cdot g(\bbeta\subo\trans \bX_i) \}^2] \\
= &  \E[(\bu\subo^\top \bX_i)^2 \{(1-\rho^2)g(\bbeta\subo\trans \bX_i)^2 - 2(1-\rho)g(\bbeta\subo\trans \bX_i) \E(Y|\bS_i,\bX_i)
+(1-\rho^2) \E(Y|\bS_i,\bX_i)^2\}] \\
= & (1-\rho)^2 \E[(\bu\subo^\top \bX_i)^2 \{ \E(Y|\bS_i,\bX_i)-g(\bbeta\subo\trans \bX_i)\}^2] \\
& + 2 \rho (1-\rho) \E[(\bu\subo^\top \bX_i)^2 \{ \E(Y|\bS_i,\bX_i)^2 + g(\bbeta\subo\trans \bX_i)^2\}].
\end{align*}
The last expression is the sum of expectations of complete squares,
so it must be non-negative.
Thus, we have shown that the SAS asymptotic variance is no greater than the supervised learning variance.
The equality holds only if
1) $\rho=1$ all samples are labelled;
2) or $\rho=0$  and $\bu\subo^\top \bX_i\{\E(Y|\bS_i,\bX_i)-g(\bbeta\subo\trans \bX_i)\} = 0$ almost surely.
\end{proof}

\subsubsection*{Efficiency Bound among Semi-parametric RAL Estimators}

\begin{proof}[Proof of Proposition \ref{prop:eff-infl}]
The proof follows the flow of Section D.2 in \cite{KallusMao2020arxiv}.
The semi-parametric model for the observed data is
\begin{align}
    \Mobs = \bigg\{ & f_{\bX,Y,\bS, R} (\bx,y,\bs, r) = f_{\bX}(\bx) f_{\bS|\bX} (\bs|\bx) \left\{\rho f_{Y|\bS,\bX}(y|\bs,\bx)\right\}^r   (1-\rho)^{1-r}: \notag
    \\ & \qquad  f_{\bX}, f_{\bS|\bX}, f_{Y\mid \bS,\bX} \text{ are arbitrary pdf/pmf}, \, \bigg\}. \label{def:obsM}
\end{align}
We consider the parametric sub-model
\begin{align}
    \Mpar = \bigg\{ & f_{\bX,Y,\bS, R} (\bx,y,\bs, r; \bgz) = f_{\bX}(\bx; \bgz) f_{\bS|\bX} (\bs|\bx; \bgz)\left\{\rho f_{Y|\bS,\bX}(y|\bs,\bx; \bgz)\right\}^r  \notag  \\
  & \qquad  \times (1-\rho)^{1-r}: \; \bgz\in\R^d \bigg\}. \label{def:parM}
\end{align}
The score vector of the parametric sub-model is
\begin{align}
  & \bgPs(\bX,Y,\bS, R) \notag \\
  = &  \left.\frac{\partial \log\{ f_{\bX,Y,\bS, R} (\bX,Y,\bS, R; \bgz)\}}{\partial \bgz}\right|_{\bgz = \bgz\subo}\notag \\
  = &\left.\frac{\partial \log\{ f_{\bX} (\bX; \bgz)\}}{\partial \bgz}\right|_{\bgz = \bgz\subo}
  + \left.\frac{\partial \log\{ f_{\bS|\bX} (\bS|\bX; \bgz)\}}{\partial \bgz}\right|_{\bgz = \bgz\subo}
  + R\left.\frac{\partial \log\{ f_{Y\mid \bS,\bX} (Y\mid \bS,\bX; \bgz)\}}{\partial \bgz}\right|_{\bgz = \bgz\subo} \notag\\
  = & \bgPs_{\bX}(\bX) +  \bgPs_{\bS}(\bS,\bX) + R \bgPs_{Y}(Y,\bS,\bX). \label{def:sp-score}
\end{align}
Next, we decompose the the Hilbert space of mean zero finite variance random variables measurable to $\sigma\{\bX,\bS,R,YR\}$,
denoted as $\Hcal$. The model tangent space spanned by the score \eqref{def:sp-score} is a linear sub-space of $\Hcal$,
\begin{gather}
\Lambda = \Lambda_{\bX} \oplus  \Lambda_{\bS} \oplus \Lambda_Y, \notag\\
\Lambda_{\bX} = \bigcup_{\Mpar} \mathrm{span}\{\bgPs_{\bX}(\bX)\} = \{h(\bX) \in \Hcal: \E[h(\bX)]=0\}, \notag\\
\Lambda_{\bS} = \bigcup_{\Mpar} \mathrm{span}\{ \bgPs_{\bS}(\bS,\bX)\} = \{h(\bS,\bX) \in \Hcal: \E[h(\bS,\bX)\mid \bX]=0\},\notag\\
\Lambda_{\bY} = \bigcup_{\Mpar} \mathrm{span}\{R \bgPs_{Y}(Y,\bS,\bX)\} = \{R h(Y,\bS,\bX) \in \Hcal: \E[h(Y,\bS,\bX)\mid \bS,\bX]=0\}.
\end{gather}
The orthogonal space of model tangent space $\Lambda$ is
\begin{equation}
    \Lambda^\perp = \{h(R, \bS,\bX) \in \Hcal: \E[h(R,\bS,\bX)\mid \bS,\bX]=0\}, \;
    \Hcal = \Lambda  \oplus \Lambda^\perp.
\end{equation}

Now, we verify that the supervised learning influence function
$$
\inflSL(\theta; \bbeta) =  \frac{R}{\rho} \bu\subo\trans\bX\{Y-g(\bbeta\trans\bX)\}
$$
is indeed an influence function for $\bx\substd\trans\bbeta$
by showing
$$
\E\{\inflSL(\theta_0; \bbeta_0) \bgPs(\bX,Y,\bS, R)\}
=  \left.\bx\substd \frac{d}{d \bgz} \bbeta(\bgz)\right|_{\bgz = \bgz\subo}.
$$
Since $\bbeta(\bgz)$ is an implicit function of $\bgz$ through the moment condition
$$
\E_{\bgz} [\bX \{g(\bbeta(\bgz)\trans\bX) - Y\} = 0,
$$
we solve for its derivative by differentiating
the moment condition
\begin{align*}
    \left. \frac{d}{d \bgz}\E_{\bgz} [\bX \{g(\bbeta(\bgz)\trans\bX) - Y\}] \right|_{\bgz=\bgz\subo} = & 0 \\
     \left. \frac{d}{d \bgz}\E_{\bgz} [\bX \{g(\bbeta\subo\trans\bX) - Y\}] \right|_{\bgz=\bgz\subo}+\E_{\bgz\subo} \{\bX\bX\trans g'(\bbeta\subo\trans\bX)\}  \left.\frac{d}{d \bgz} \bbeta(\bgz)\right|_{\bgz = \bgz\subo}
     = & 0\\
     -\Prec\E_{\bgz\subo} [\bX \{g(\bbeta\subo\trans\bX) - Y\}\{\bgPs_{\bX}(\bX) +  \bgPs_{\bS}(\bS,\bX) + \bgPs_{Y}(Y,\bS,\bX)\}]
     = &  \left.\frac{d}{d \bgz} \bbeta(\bgz)\right|_{\bgz = \bgz\subo}.
\end{align*}
Then, we verify that the supervised learning influence function is valid
\begin{align*}
  \left.\frac{d}{d \bgz} \bx\substd\trans\bbeta(\bgz)\right|_{\bgz = \bgz\subo}   = &   -\E_{\bgz\subo} \left[\frac{R}{\rho}\bu\subo\trans\bX \{g(\bbeta\subo\trans\bX) - Y\}\bgPs(\bX,Y,\bS, R)\right]  \\
 =&   \E\{\inflSL(\theta_0; \bbeta_0) \bgPs(\bX,Y,\bS, R)\}.
\end{align*}

Finally, we derive the efficient influence function by subtract from $\inflSL$ its projection onto $\Lambda^\perp = \Lambda_R$.
Let $\Pi[ h(\bD) \mid \Lambda]$ be the projection of $h(\bD) \in \Hcal$ to the space $\Lambda$.
We can easily calculate the projection of $\inflSL$ onto $\Lambda_R$,
\begin{align*}
   \Pi[ \inflSL(\theta_0; \bbeta_0) \mid \Lambda_R] = & \E\{\inflSL(\theta_0; \bbeta_0) \mid R,\bS, \bX\}
   -\E\{\inflSL(\theta_0; \bbeta_0) \mid \bS, \bX\} \\
   = & \frac{R}{\rho} \bu\subo\trans\bX\{\E(Y\mid \bS,\bX)-g(\bbeta\trans\bX)\}
   - \bu\subo\trans\bX\{\E(Y\mid \bS,\bX)-g(\bbeta\trans\bX)\}.
\end{align*}
The efficient influence function is thus obtained
\begin{align*}
\infleff(\theta_0; \bbeta_0) = & \inflSL(\theta_0; \bbeta_0) - \Pi[ \inflSL(\theta_0; \bbeta_0) \mid \Lambda_R ] \\
=& \frac{R}{\rho} \bu\subo\trans\bX\{Y-g(\bbeta\trans\bX)\} - \frac{R}{\rho} \bu\subo\trans\bX\{\E(Y\mid \bS,\bX)-g(\bbeta\trans\bX)\} \\
&   + \bu\subo\trans\bX\{\E(Y\mid \bS,\bX)-g(\bbeta\trans\bX)\} \\
= & \frac{R}{\rho} \bu\subo\trans\bX\{Y-\E(Y\mid \bS,\bX)\}
   + \bu\subo\trans\bX\{\E(Y\mid \bS,\bX)-g(\bbeta\trans\bX)\}.
\end{align*}
\end{proof}

\section{Auxiliary Results}\label{asection:aux}

\subsection{General}

\begin{lemma}\label{lemma:subg-g}
Under Assumptions \ref{assume:Y-tail}, \ref{assume:link}, \ref{assume:X-subG},
the residuals of the imputed loss are sub-Gaussian random variables ,
$$
\|g(\bbeta\trans \bX_i) - g({\bgamma}\trans \bW_i)\|_{\psi_2} \\
\le  4\max\{\nu_1,\nu_2, M\|\bbeta-\bbeta\subo\|_2 \sigma\submax/\sqrt{2},
  M\|{\bgamma}-\bgamma\subo\|_2 \sigma\submax/\sqrt{2}\}
$$
Similarly,
\begin{gather*}
\|Y_i - g({\bgamma}^{\top} \bW_i)\|_{\psi_2}
\le 2\max\{\nu_1,
  M\|{\bgamma}-\bgamma\subo\|_2\sigma\submax/\sqrt{2}\},\\
  \|g({\bgamma}^{\top} \bW_i) - g({\bgamma}\subo^{\top} \bW_i)\|_{\psi_2}
\le M\|{\bgamma}-\bgamma\subo\|_2\sigma\submax/\sqrt{2}, \\
    \|g({\bbeta}^{\top} \bX_i) - g({\bbeta}\subo^{\top} \bX_i)\|_{\psi_2}
\le
  M\|{\bbeta}-\bbeta\subo\|_2 \sigma\submax/\sqrt{2},\\
  \begin{aligned}
 \|\rho \cdot g(\bbeta\trans \bX_i) + (1-\rho)\cdot g({\bgamma}\trans \bW_i) &- Y\|_{\psi_2}
\le  4\max\{(1-\rho)\nu_1,\rho\nu_2, \\ &\rho M\|\bbeta-\bbeta\subo\|_2 \sigma\submax/\sqrt{2},
  (1-\rho) M\|{\bgamma}-\bgamma\subo\|_2 \sigma\submax/\sqrt{2}\},
  \end{aligned} \\
 \|g({\bbeta}^{\top} \bX_i) - g({\bbeta}\subo^{\top} \bX_i) -
 g({\bgamma}^{\top} \bW_i) + g({\bgamma}\subo^{\top} \bW_i)\|_{\psi_2}
\le \sqrt{2}M\sigma\submax \max\{\|{\bbeta}-\bbeta\subo\|_2, \|{\bgamma}-\bgamma\subo\|_2\}.
\end{gather*}
\end{lemma}

\begin{proof}[Proof of Lemma \ref{lemma:subg-g}]
  To establish the sub-exponential tail,
  we consider the following decomposition
  \begin{equation}\label{eq:dldag-term}
  \begin{aligned}
    g(\bbeta\trans \bX_i) - g(\bgamma\trans \bW_i)
    = & \{g(\bbeta\subo\trans \bX_i) - Y_i\} - \{g(\bgamma\subo\trans \bW_i)-Y_i\} \\
&    + \{g(\bbeta\trans \bX_i) - g(\bbeta\subo\trans \bX_i)\}
    - \{g(\bgamma\trans \bW_i) - g(\bgamma\subo\trans \bW_i)\}.
  \end{aligned}
  \end{equation}
  According to Assumption \ref{assume:Y-tail},
  the first two terms on the right-hand side of \eqref{eq:dldag-term}
  are sub-Gaussian,
  $$
  \|g(\bbeta\subo\trans \bX_i) - Y_i\|_{\psi_2} \le \nu_1, \;
  \|g(\bgamma\subo\trans \bW_i)-Y_i\|_{\psi_2} \le \nu_2.
  $$
  According to Assumption \ref{assume:link},
  the latter two terms on the right-hand side of \eqref{eq:dldag-term}
  are bounded by
  \begin{equation*}
    |g(\bbeta\trans \bX_i) - g(\bbeta\subo\trans \bX_i)| \le M |(\bbeta-\bbeta\subo)\trans \bX_i|, \;
   |g(\bgamma\trans \bW_i) - g(\bgamma\subo\trans \bW_i)| \le M |(\bgamma-\bgamma\subo)\trans \bW_i|.
  \end{equation*}
  Under Assumption \ref{assume:X-subG},
  $(\bbeta-\bbeta\subo)\trans \bX_i$ and $(\bgamma-\bgamma\subo)\trans \bW_i$
  are sub-Gaussian random variables, 
  \begin{gather*}
   \|(\bbeta-\bbeta\subo)\trans \bX_i\|_{\psi_2} \le \|\bbeta-\bbeta\subo\|_2 \sigma\submax/\sqrt{2} \\
    \|(\bgamma-\bgamma\subo)\trans \bW_i\|_{\psi_2} \le \|\bgamma-\bgamma\subo\|_2 \sigma\submax/\sqrt{2}.
  \end{gather*}
  By Lemma \ref{lemma:sub-mb},
    \begin{gather*}
  \|g(\bbeta\trans \bX_i) - g(\bbeta\subo\trans \bX_i)\|_{\psi_2} \le M\|(\bbeta-\bbeta\subo)\trans \bX_i\|_{\psi_2} \le M\|\bbeta-\bbeta\subo\|_2 \sigma\submax/\sqrt{2} \\
   \|g(\bgamma\trans \bW_i) - g(\bgamma\subo\trans \bW_i)\|_{\psi_2} \le \|(\bgamma-\bgamma\subo)\trans \bW_i\|_{\psi_2} \le M\|\bgamma-\bgamma\subo\|_2 \sigma\submax/\sqrt{2}.
  \end{gather*}
  Finally, we apply Lemma \ref{lemma:sub-add}
  \begin{align*}
  \|g(\bbeta\trans \bX_i) - g(\bgamma\trans \bW_i)\|_{\psi_2}
  \le & 4 \max\big\{\|g(\bbeta\subo\trans \bX_i) - Y_i\|_{\psi_2}, \|g(\bgamma\subo\trans \bW_i)-Y_i\|_{\psi_2}, \\
 & \|g(\bbeta\trans \bX_i) - g(\bbeta\subo\trans \bX_i)\|_{\psi_2}, \|g(\bgamma\trans \bW_i) - g(\bgamma\subo\trans \bW_i)\|_{\psi_2} \big\} \\
  \le & 4\max\left\{\nu_1,\nu_2, M\|\bbeta-\bbeta\subo\|_2 \sigma\submax/\sqrt{2},
  M\|\bgamma-\bgamma\subo\|_2 \sigma\submax/\sqrt{2}\right\}.
  \end{align*}
   Therefore, we have reached the conclusion.

   We may obtain the rest of bounds
   following the same derivation.
\end{proof}

\subsection{Inference}\label{asection:inf-remainders}

\subsubsection*{Analysis of Estimated Precision Matrix}

\begin{proof}[Proof of Lemma \ref{lemma:u-oracle}]

The definition of the cross-fitted loss functions $m\supkk$ and their derivatives
can be found at \eqref{def:m}.
By the definition of $\hat{\bu}\supk$, we have
\begin{equation*}
  \sum_{k' \neq k} \frac{N_{k'}}{N-N_k} m\supkk\left(\hat{\bu}\supk;\hat{\bbeta}\supkk\right) + \lamu\|\hat{\bu}\supk\|_1 \le \sum_{k'\neq k} \frac{N_{k'}}{N-N_k}  m\supkk\left(\bu\subo;\hat{\bbeta}\supkk\right) + \lamu\|\bu\subo\|_1.
\end{equation*}
Denote the standardized estimation error as $\bgd = (\hat{\bu}\supk-\bu\subo)/\|\hat{\bu}\supk-\bu\subo\|_2$.
Due to convexity of the loss function, we have for $t = \|\hat{\bu}\supk-\bu\subo\|_2 \wedge 1$
\begin{equation}\label{eq:u-loss-min-dir}
   \sum_{k'\neq k} \frac{N_{k'}}{N-N_k}  m\supkk\left(\bu\subo + t\bgd;\hat{\bbeta}\supkk\right) + \lamu\|\bu\subo + t\bgd\|_1 \le \sum_{k'\neq k} \frac{N_{k'}}{N-N_k}  m\supkk\left(\bu\subo;\hat{\bbeta}\supkk\right) + \lamu\|\bu\subo\|_1.
\end{equation}
By the triangle inequality $\|\bu\subo\|_1 -  \|\bu\subo + t\bgd\|_1 \le t\|\bgd\|_1$, we have from \eqref{eq:u-loss-min-dir}
\begin{equation}\label{eq:u-loss-dir-triangle}
  \sum_{k'\neq k} \frac{N_{k'}}{N-N_k} \left\{ m\supkk\left(\bu\subo + t\bgd;\hat{\bbeta}\supkk\right)-m\supkk\left(\bu\subo;\hat{\bbeta}\supkk\right)\right\} \le t \lamu \|\bgd\|_1
\end{equation}
Because the loss functions $m\supk$ are quadratic functions of $\bu$,
we can apply the restricted strong convexity event $\Omega\supk$ to obtain
\begin{align}
 & m\supkk\left(\bu\subo + t\bgd;\hat{\bbeta}\supkk\right)-m\supkk\left(\bu\subo;\hat{\bbeta}\supkk\right) - t\bgd^\top \dm\supkk\left(\bu\subo;\hat{\bbeta}\supkk\right) \notag \\
 = & t^2 \bgd^\top\ddm\supkk\left(\hat{\bbeta}\supkk\right)\bgd
  \notag \\
  \ge & t^2\crsc{1}^* -t\crsc{1}^*\crsc{2}^*\sqrt{\log(p)/N_{k'}}\|\bgd\|_1. \label{eq:u-rsc-imputed}
\end{align}
Applying \eqref{eq:u-rsc-imputed} to \eqref{eq:u-loss-dir-triangle}, we have with large probability
\begin{equation*}
\sum_{k'\neq k} \frac{N_{k'}}{N-N_k}\left\{t \bgd^\top \dm\supkk\left(\bu\subo;\hat{\bbeta}\supkk\right)
+ t^2\crsc{1}^* -t\crsc{1}^*\crsc{2}^*\sqrt{\log(p)/N_{k'}}\|\bgd\|_1\right\}
   \le t \lamu\|\bgd\|_1
\end{equation*}
where $\|\bgd\|_2 = 1$ from definition.
Thus, we have reach
\begin{equation}\label{eq:lemu-bound}
 t\crsc{1}^*
   \le  \lamu\|\bgd\|_1 - \sum_{k'\neq k} \frac{N_{k'}}{N-N_k}\left\{\bgd^\top \dm\supkk\left(\bu\subo;\hat{\bbeta}\supkk\right) -\crsc{1}^*\crsc{2}^*\sqrt{\log(p)/N_{k'}}\|\bgd\|_1\right\}.
\end{equation}

The target parameter $\bu\subo$ can be identify by $\E\left\{\dm\supkk\left(\bu\subo;\hat{\bbeta}\supkk\right)\mid \Dscr_{k'}^c\right\} = 0$.
We use the fact to do a careful analysis of $\bgd^\top \dm\supkk\left(\bu\subo;\hat{\bbeta}\supkk\right)$
by the decomposition
\begin{align}
  \left|\bgd^\top \dm\supkk\left(\bu\subo;\hat{\bbeta}\supkk\right)\right|
= &\bgd^\top\left[\dm\supkk\left(\bu\subo;\hat{\bbeta}\supkk\right) - \E\left\{\dm\supkk\left(\bu\subo;\hat{\bbeta}\supkk\right)\mid \Dscr_{k'}^c\right\}
 \right] \notag \\
 &+ \bgd^\top\left[\E\left\{\dm\supkk\left(\bu\subo;\hat{\bbeta}\supkk\right)\mid \Dscr_{k'}^c\right\}-\E\left\{\dm\supkk\left(\bu\subo;\bbeta\subo\right)\mid \Dscr_{k'}^c\right\}\right] \notag \\
  \le & \|\bgd\|_1 \left\|\dm\supkk\left(\bu\subo;\hat{\bbeta}\supkk\right) - \E\left\{\dm\supkk\left(\bu\subo;\hat{\bbeta}\supkk\right)\mid \Dscr_{k'}^c\right\}
 \right\|_\infty \notag \\
&  + \left\|\E_{i \in \Ical_{k'}\cup\Jcal_{k'}}[\bX_i\bX_i\trans\bu\subo\{g(\bbeta\subo^\top \bX_i) - g(\hat{\bbeta}\supkkt \bX_i)\}\mid \Dscr_{k'}^c]\right\|_2.
  \label{eq:lemu-decomp-score}
\end{align}
We establish the rate for $L_2$-norm of the population score at $\bu\subo$ through analyzing
$$
\sup_{\|\bv\|_2 = 1} \E_{i \in \Ical_{k'}\cup\Jcal_{k'}}[\bv\trans\bX_i\bX_i\trans\bu\subo\{g(\bbeta\subo^\top \bX_i) - g(\hat{\bbeta}\supkkt \bX_i)\}\mid \Dscr_{k'}^c],
$$
whose bound can be derived from Assumptions \ref{assume:link}, \ref{assume:X-subG}, \ref{assume:X-subG},
 the Cauchy-Schwartz inequality and Lemma \ref{lemma:sub-moment},
\begin{align*}
    & \E_{i \in \Ical_{k'}\cup\Jcal_{k'}}[\bv\trans\bX_i\bX_i\trans\bu\subo\{g(\bbeta\subo^\top \bX_i) - g(\hat{\bbeta}\supkkt \bX_i)\}\mid \Dscr_{k'}^c] \\
\le &     M \E_{i \in \Ical_{k'}\cup\Jcal_{k'}}[|\bv\trans\bX_i\bX_i\trans\bu\subo\{\left(\bbeta\subo-\hat{\bbeta}\supkk\right)\trans \bX_i\}|\mid \Dscr_{k'}^c] \\
\le & M \left[\E\{(\bv\trans\bX_i)^4\}\E\{(\bv\trans\bX_i)^4\}\right]^{1/4} \sqrt{\E_{i \in \Ical_{k'}\cup\Jcal_{k'}}\left[\left\{\left(\bbeta\subo-\hat{\bbeta}\supkk\right)\trans \bX_i\right\}^2\mid \Dscr_{k'}^c\right]} \\
\le &  M\sigma\submax^3 \|\bv\|_2 \|\bu\subo\|_2 \left\|\bbeta\subo-\hat{\bbeta}\supkk\right\|_2.
\end{align*}
Hence, we have shown
\begin{equation}\label{eq:dm-u0}
    \left\|\E_{i \in \Ical_{k'}\cup\Jcal_{k'}}\left[\bX_i\bX_i\trans\bu\subo\left\{g(\bbeta\subo^\top \bX_i) - g\left(\hat{\bbeta}\supkkt \bX_i\right)\right\}\mid \Dscr_{k'}^c\right]\right\|_2
    \le  M\sigma\submax^3  \|\bu\subo\|_2 \left\|\bbeta\subo-\hat{\bbeta}\supkk\right\|_2.
\end{equation}

By the bound for \eqref{eq:lemu-decomp-score} through
\eqref{eq:dm-u0} and the definition of $\lamb$,
we have the bound from \eqref{eq:lemu-bound}
\begin{equation}
  t\crsc{1}^*
   \le  2\lamu\|\bgd\|_1 +   M\sigma\submax^3 \|\bu\subo\|_2\sup_{k'\neq k}\left\|\bbeta\subo-\hat{\bbeta}\supkk\right\|_2.  \label{eq:lemu-case}
\end{equation}
Hence, we can reach
an immediate bound for estimation error
from \eqref{eq:lemu-case} without considering the sparsity of $\bu\subo$.
We shall proceed to derive a sharper bound
that involves the sparsity of $\bu\subo$.
We separately analyze two cases.

\textbf{Case 1: }

\begin{equation*}
   M\sigma\submax^3 \|\bu\subo\|_2 \sup_{k=1,\dots,\Kn}\left\|\bbeta\subo-\hat{\bbeta}\supkk\right\|_2
  \ge  \|\bgd\|_1\lamu/3
\end{equation*}
In this case, the estimation error is dominated by $\bbeta\subo-\hat{\bbeta}\supkk$.
We simply have from \eqref{eq:lemu-case}
\begin{gather*}
  t\crsc{1}^*
   \le  7 M\sigma\submax^3 \|\bu\subo\|_2 \sup_{k=1,\dots,\Kn}\left\|\bbeta\subo-\hat{\bbeta}\supkk\right\|_2, \\
   t\crsc{1}^* \|\bgd\|_1\lamu/3 \le 7 M \sigma\submax^6  \|\bu\subo\|_2^2 \sup_{k=1,\dots,\Kn}\left\|\bbeta\subo-\hat{\bbeta}\supkk\right\|_2^2.
   \end{gather*}
Thus, we have
\begin{gather}\label{eq:lemu-rate1}
  \|\hat{\bu}\supk-\bu\subo\|_2 \le 7 M\sigma\submax^3 \|\bu\subo\|_2\sup_{k' \neq k}\left\|\bbeta\subo-\hat{\bbeta}\supkk\right\|_2/\crsc{1}^*, \notag \\
   \|\hat{\bu}\supk-\bu\subo\|_1 \le 21 M \sigma\submax^6  \|\bu\subo\|_2^2 \sup_{k' \neq k}\left\|\bbeta\subo-\hat{\bbeta}\supkk\right\|_2^2/(\crsc{1}^*\lamu).
\end{gather}

\textbf{Case 2: }
\begin{equation}\label{eq:lemu-case2}
 M \sigma\submax^3  \|\bu\subo\|_2 \sup_{k' \neq k}\left\|\bbeta\subo-\hat{\bbeta}\supkk\right\|_2
  \le  \|\bgd\|_1\lamu/3
\end{equation}
In this case, the estimation error is comparable to the situation that we have the true $\bbeta\subo$
for the Hessian.
Thus, the sparsity of $\bu\subo$ may affect the estimation error.

Following the typical approach to establish the cone condition for $\bgd$,
we analyze the symmetrized Bregman's divergence,
\begin{align}
  &(\hat{\bu}\supk - \bu\subo)^\top  \sum_{k'=k}\frac{N_{k'}}{N-N_k}\left\{ \dm\supkk\left(\hat{\bu}\supk;\hat{\bbeta}\supkk\right)
  - \dm\supkk\left(\bu\subo;\hat{\bbeta}\supkk\right) \right\} \notag \\
  =& \|\hat{\bu}\supk - \bu\subo\|_2 \sum_{k'=k}\frac{N_{k'}}{N-N_k}\bgd^\top \left\{ \dm\supkk\left(\hat{\bu}\supk;\hat{\bbeta}\supkk\right)
  - \dm\supkk\left(\bu\subo;\hat{\bbeta}\supkk\right) \right\}.\label{eq:lemu-breg}
\end{align}
Due to the convexity of the quadratic loss $m\supkk\left(\cdot;\hat{\bbeta}\supkk\right)$,
the symmetrized Bregman's divergence \eqref{eq:lemu-breg} is nonnegative
through a mean-value theorem,
\begin{align*}
  & (\hat{\bu}\supk - \bu\subo)^\top  \sum_{k'=k}\frac{N_{k'}}{N-N_k}\left\{ \dm\supkk\left(\hat{\bu}\supk;\hat{\bbeta}\supkk\right)
  - \dm\supkk\left(\bu\subo;\hat{\bbeta}\supkk\right) \right\}\\
  = & \sum_{k'=k}\frac{N_{k'}}{N-N_k}\sum_{i \in \Ical_{k'}\cup \Jcal_{k'}} g'(\hat{\bbeta}\supkkt \bX_i) \{ (\hat{\bu}\supk - \bu\subo)^\top \bX_i\}^2 \notag \\
  \ge & 0.
\end{align*}
Denote the indices set of nonzero coefficient in $\bu\subo$ as
$\nzu = \{j: \bu_{0,j} \neq 0\}$.
We denote the $\bgd_{\nzu}$ and $\bgd_{\nzuc}$
as the sub-vectors for $\bgd$ at positions in $\nzu$ and at positions not in $\nzu$,
respectively.
The solution $\hat{\bu}\supk$ satisfies the KKT condition
\begin{gather*}
 \left\| \sum_{k'=k}\frac{N_{k'}}{N-N_k} \dm\supkk\left(\hat{\bu}\supk;\hat{\bbeta}\supkk\right)\right\|_\infty \le \lamu, \\
  \sum_{k'=k}\frac{N_{k'}}{N-N_k} \dot{m}\supkk\left(\hat{\bu}\supk;\hat{\bbeta}\supkk\right)_j = -\lamu\sgn(\hat{u}\supk_j), \, j:\hat{u}\supk_j \neq 0.
\end{gather*}
From the KKT condition and the definitions of $\bgd$ and $\nzu$, we have
\begin{gather}
  \delta_j \sum_{k'=k}\frac{N_{k'}}{N-N_k} \dot{m}\supkk\left(\hat{\bu}\supk;\hat{\bbeta}\supkk\right)_j \le |\delta_j| \lamu, \, j \in \nzb; \notag \\
  \delta_j \sum_{k'=k}\frac{N_{k'}}{N-N_k} \dot{m}\supkk\left(\hat{\bu}\supk;\hat{\bbeta}\supkk\right)_j = \frac{-\hat{u}\supk_j \lamu\sgn(\hat{u}\supk_j) }{\|\hat{\bu}\supk-\bu\subo\|_2} = -\lamu|\delta_j|, \, j \in \nzbc. \label{eq:u-KKT-breg}
\end{gather}
Applying the \eqref{eq:u-KKT-breg} to \eqref{eq:lemu-breg}, we have the upper bound,
\begin{align*}
  & \bgd^\top  \sum_{k'=k}\frac{N_{k'}}{N-N_k}\left\{ \dm\supkk\left(\hat{\bu}\supk;\hat{\bbeta}\supkk\right)
  - \dm\supkk\left(\bu\subo;\hat{\bbeta}\supkk\right) \right\} \\
  = & \sum_{j \in \nzu} \delta_j \sum_{k'=k}\frac{N_{k'}}{N-N_k} \dot{m}\supkk\left(\hat{\bu}\supk;\hat{\bbeta}\supkk\right)_j
  + \sum_{j \in \nzuc} \delta_j \sum_{k'=k}\frac{N_{k'}}{N-N_k} \dot{m}\supkk\left(\hat{\bu}\supk;\hat{\bbeta}\supkk\right)_j \notag \\
  &+ \sum_{k'=k}\frac{N_{k'}}{N-N_k} \bgd^\top\dm\supkk\left(\bu\subo;\hat{\bbeta}\supkk\right)  \notag \\
  \le &  \lamu\sum_{j \in \nzu} |\delta_j|
  -  \lamu\sum_{j \in \nzuc} |\delta_j|
  + \bgd^\top \dell^\dag(\bbeta\subo;\hat{\bgamma}) \notag \\
  \le & \lamu\|\bgd_{\nzu}\|_1 - \lamu\|\bgd_{\nzuc}\|_1
  + \left|\sum_{k'=k}\frac{N_{k'}}{N-N_k} \bgd^\top\dm\supkk\left(\bu\subo;\hat{\bbeta}\supkk\right)\right|.
\end{align*}
Then, we apply \eqref{eq:lemu-decomp-score}, the definition of $\lamu$ and \eqref{eq:lemu-case2},
\begin{align*}
  0 \le & \lamu\|\bgd_{\nzu}\|_1 - \lamu\|\bgd_{\nzuc}\|_1
  + \frac{2}{3} \lamu\|\bgd\|_1 , \quad \mbox{and}\quad
 \lamu\|\bgd_{\nzuc}\|_1  \le & 5 \lamu\|\bgd_{\nzu}\|_1.
\end{align*}
Therefore, we can bound the $L_1$ norm of $\bgd$ by the cone property,
\begin{equation}\label{eq:lemu-cone}
\|\bgd\|_1 \le 6 \lamu\|\bgd_{\nzu}\|_1
\le 6\sqrt{\su}\|\bgd\|_2 = 6\sqrt{\su}.
\end{equation}

Now, we apply the cone condition \eqref{eq:lemu-cone} and the case condition \eqref{eq:lemu-case2}
to the bound \eqref{eq:lemu-case},
\begin{equation*}
   t\crsc{1}^*
   \le  14\sqrt{\su}\lamu, \;
   t\crsc{1}^* \|\bgd\|_1 \le  84\su\lamu
\end{equation*}
Thus, we obtain the rate for estimation error
\begin{equation}\label{eq:lemu-rate2}
  \|\hat{\bu}\supk-\bu\subo\|_2 \le 14\sqrt{\su}\lamu/\crsc{1}^*, \;
  \|\hat{\bu}\supk-\bu\subo\|_1 \le 84\su\lamu/\crsc{1}^*.
\end{equation}

\textbf{Conclusion: }

Since Case 1 and Case 2 are the complement of each other,
one of them must occur.
Thus, the bound of estimation error is controlled by the larger bound in
the two cases,
\begin{gather*}
  \|\hat{\bu}\supk-\bu\subo\|_2 \le \max\left\{14\sqrt{\su}\lamu/\crsc{1}^*,
   7 M \sigma\submax^3  \|\bu\subo\|_2 \sup_{k'\neq k}\left\|\bbeta\subo-\hat{\bbeta}\supkk\right\|_2/\crsc{1}^*\right\}, \\
     \|\hat{\bu}\supk-\bu\subo\|_1 \le \max\left\{84\su\lamu/\crsc{1}^*,
   21 M \sigma\submax^6  \|\bu\subo\|_2^2 \sup_{k'\neq k}\left\|\bbeta\subo-\hat{\bbeta}\supkk\right\|_2^2/(\crsc{1}^*\lamu)\right\},
\end{gather*}
which is our oracle inequality.
\end{proof}

\subsubsection*{Analysis for Terms $T_1$-$T_5$ in Part 1}

To show
$$
T_1 = \sqrt{n}\left\{\bx\substd - \Hess(\tilde{\bbeta})\bu\subo\right\}(\hat{\bbeta}\supk - \bbeta\subo)
= O_p\left(\sqrt{n}\|\hat{\bbeta}\supk-\bbeta\subo\|_2^2\right),
$$
we rewrite the term as a conditional expectation
\begin{align*}
  T_1 = & \sqrt{n}\bu\subo\trans\left\{\Sigb- \Hess(\tilde{\bbeta})\right\}(\hat{\bbeta}\supk - \bbeta\subo) \\
  = & \sqrt{n}\E_{i\in \Jcal_k} \left[ \bu\subo\trans\bX_i (\hat{\bbeta}\supk - \bbeta\subo)\trans \bX_i
  \{g'(\tilde{\bbeta}\trans\bX_i) - g'(\hat{\bbeta}\supkt\bX_i)\}\mid \Dscr_k^c\right].
\end{align*}
Under Assumptions \ref{assume:link}, \ref{assume:X-subG}, we derive the bound for the expectation using the Cauchy-Schwartz inequality and Lemma \ref{lemma:sub-moment},
\begin{align}
 |T_1| \le & M \sqrt{n}\E_{i\in \Jcal_k} \left[ \bu\subo\trans\bX_i
 \{(\hat{\bbeta}\supk - \bbeta\subo)\trans \bX_i\}^2 \mid \Dscr_k^c\right] \notag  \\
 \le & M \sqrt{n \E_{i\in \Jcal_k} \left\{ (\bu\subo\trans\bX_i)^2 \mid \Dscr_k^c\right\}
 \E_{i\in \Jcal_k} \left[
 \{(\hat{\bbeta}\supk - \bbeta\subo)\trans \bX_i\}^4 \mid \Dscr_k^c\right]} \notag \\
 \le& M \sqrt{n 8 \|\bu\subo\trans\bX_i\|_{\psi_2}^2\|(\hat{\bbeta}\supk - \bbeta\subo)\trans \bX_i\|_{\psi_2}^4} \notag \\
 \le& \sqrt{n} M \|\bu\subo\|_2 \left\|\hat{\bbeta}\supk - \bbeta\subo\right\|_2^2 \sigma\submax^3. \label{eq:T1-exp}
\end{align}
Since $\|\bu\subo\|_2$ is bounded according to \eqref{eq:u0b-ukt},
we have established in
$$
|T_1|= O_p\left(\sqrt{n}\|\hat{\bbeta}\supk-\bbeta\subo\|_2^2\right)
$$
as declared.

To show
$$
T_2 = \sqrt{n}\left(\bu\subo-\hat{\bu}\supk\right)\trans\Hess(\tilde{\bbeta})(\hat{\bbeta}\supk - \bbeta\subo)
= O_p\left(\sqrt{n}\|\hat{\bbeta}\supk-\bbeta\subo\|_2\|\hat{\bu}\supk-\bu\subo\|_2\right),
$$
we rewrite the term as a conditional expectation
$$
  T_2 =  \sqrt{n}\E_{i\in \Jcal_k} \left[ \left(\bu\subo-\hat{\bu}\supk\right)\trans \bX_i (\hat{\bbeta}\supk - \bbeta\subo)\trans \bX_i
  g'(\tilde{\bbeta}\trans\bX_i)\mid \Dscr_k^c\right].
$$
Similar to \eqref{eq:T1-exp}, we derive the bound for the expectation under Assumptions \ref{assume:link}, \ref{assume:X-subG} through the Cauchy-Schwartz inequality and Lemma \ref{lemma:sub-moment}, \ref{lemma:sub-mg},
\begin{align*}
  |T_2| \le  & M \sqrt{n}\E_{i\in \Jcal_k} \left[ |\left(\bu\subo-\hat{\bu}\supk\right)\trans \bX_i (\hat{\bbeta}\supk - \bbeta\subo)\trans \bX_i| \mid \Dscr_k^c\right] \\
  \le & 2M\sqrt{n} \|\left(\bu\subo-\hat{\bu}\supk\right)\trans \bX_i (\hat{\bbeta}\supk - \bbeta\subo)\trans \bX_i\|_{\psi_1} \\
  \le & 2M\sqrt{n} \|\left(\bu\subo-\hat{\bu}\supk\right)\trans \bX_i\|_{\psi_2} \|(\hat{\bbeta}\supk - \bbeta\subo)\trans \bX_i\|_{\psi_2} \\
  \le & M\sqrt{n} \|\hat{\bu}\supk-\bu\subo\|_2 \|\hat{\bbeta}\supk - \bbeta\subo\|_2 \sigma\submax^2.
\end{align*}
This bound immediately implies
$$
T_2 = O_p\left(\sqrt{n}\|\hat{\bbeta}-\bbeta\subo\|_2\|\hat{\bu}\supk-\bu\subo\|_2\right).
$$

To show
\begin{align*}
  T_3 =  &  \sqrt{n} \hat{\bu}\supkt\left[ \E\{\dlkdag(\hat{\bbeta}\supk; \hat{\bgamma}\supk)\mid \Dscr_k^c\}
     - \left\{\dlkdag(\hat{\bbeta}\supk; \hat{\bgamma}\supk) - \dlkdag(\bbeta\subo; \bgamma\subo)\right\}\right] \\
     = & O_p\left(\rho\|\hat{\bbeta}\supk-\bbeta\subo\|_2+\sqrt{\rho(1-\rho)}\|\hat{\bgamma}\supk-\bgamma\subo\|_2\right),
\end{align*}
we rewrite the term as two empirical processes with diminishing summands
\begin{align*}
   T_3 =  & - \sqrt{n} \frac{1}{N_k}\sum_{i\in \Jcal_k} \bigg(\hat{\bu}\supkt\bX_i \{g(\hat{\bbeta}\supkt\bX_i) - g(\bbeta\subo\trans\bX_i)-g(\hat{\bgamma}\supkt\bW_i) + g(\bgamma\subo\trans\bW_i) \} \\
   & - \E_{i \in \Jcal_k}\left[\hat{\bu}\supkt\bX_i \{g(\hat{\bbeta}\supkt\bX_i) - g(\bbeta\subo\trans\bX_i)-g(\hat{\bgamma}\supkt\bW_i) + g(\bgamma\subo\trans\bW_i) \}\mid \Dscr_k^c\right]\bigg) \\
   & - \sqrt{n} \frac{1}{N_k}\sum_{i\in \Ical_k} \bigg(\hat{\bu}\supkt\bX_i \{g(\hat{\bbeta}\supkt\bX_i) - g(\bbeta\subo\trans\bX_i)\} \\
   & - \E_{i \in \Jcal_k}\left[\hat{\bu}\supkt\bX_i \{g(\hat{\bbeta}\supkt\bX_i) - g(\bbeta\subo\trans\bX_i)\}\mid \Dscr_k^c\right]\bigg).
\end{align*}
We have used the identity
$
\E\{\dlkdag(\bbeta\subo; \bgamma\subo)\mid \Dscr_k^c\} = 0
$
above.
Using Lemmas \ref{lemma:subg-g}, \ref{lemma:sub-bern} and Assumptions \eqref{assume:X-subG}
and \eqref{assume:X-subG}, we show that each summand is sub-exponential
\begin{gather*}
\begin{aligned}
   & \|\hat{\bu}\supkt\bX_i \{g(\hat{\bbeta}\supkt\bX_i) - g(\bbeta\subo\trans\bX_i)-g(\hat{\bgamma}\supkt\bW_i) + g(\bgamma\subo\trans\bW_i) \}\|_{\psi_1} \\
   \le & \|\hat{\bu}\supkt\bX_i\|_{\psi_2} \|g(\hat{\bbeta}\supkt\bX_i) - g(\bbeta\subo\trans\bX_i)-g(\hat{\bgamma}\supkt\bW_i) + g(\bgamma\subo\trans\bW_i)\|_{\psi_2} \\
   \le & M\sigma\submax^2 \|\hat{\bu}\supk\|_2\left(\|\hat{\bbeta}\supk-\bbeta\subo\|_2+\|\hat{\bgamma}\supk-\bgamma\subo\|_2\right), \end{aligned}\\
    \|\hat{\bu}\supkt\bX_i \{g(\hat{\bbeta}\supkt\bX_i) - g(\bbeta\subo\trans\bX_i)\}\|_{\psi_1}
   \le  M\sigma\submax^2 \|\hat{\bu}\supk\|_2\|\hat{\bbeta}\supk-\bbeta\subo\|_2/2.
\end{gather*}
Applying the Bernstein's inequality, we obtain
$$
T_3 =  O_p\left(\|\hat{\bu}\supk\|_2\left\{\sqrt{\rho}\|\hat{\bbeta}\supk-\bbeta\subo\|_2+\sqrt{\rho(1-\rho)}\|\hat{\bgamma}\supk-\bgamma\subo\|_2\right\}\right).
$$
We achieve the stated rate with the tightness of $\|\hat{\bu}\supk\|_2$ from \eqref{eq:u0b-ukt}.

To show
\begin{align*}
  T_4 =  & \sqrt{n} (1-\rho) \hat{\bu}\supkt\left[ \E\{\dlkimp(\hat{\bgamma}\supk)\mid \Dscr_k^c\}
     - \left\{\dlkimp(\hat{\bgamma}\supk) - \dlkimp(\bgamma\subo)\right\}\right] \\
     = & O_p\left((1-\rho)\|\hat{\bgamma}\supk-\bgamma\subo\|_2\right),
\end{align*}
we rewrite the term as the empirical process with diminishing summands
\begin{align*}
  T_4 = & - \sqrt{n} (1-\rho)\frac{1}{n_k}\sum_{i\in \Ical_k} \bigg(\hat{\bu}\supkt\bX_i \{g(\hat{\bgamma}\supkt\bW_i) - g(\bgamma\subo\trans\bW_i)\} \\
   & - \E_{i \in \Jcal_k}\left[\hat{\bu}\supkt\bX_i \{g(\hat{\bgamma}\supkt\bW_i) - g(\bgamma\subo\trans\bW_i)\}\mid \Dscr_k^c\right]\bigg).
\end{align*}
We have used the identity
$
\E\{\dlkimp(\bgamma\subo)\mid \Dscr_k^c\} = 0
$
above.
Similar to the analysis of $T_3$, we show that each summand is sub-exponential
$$
   \|\hat{\bu}\supkt\bX_i \{g(\hat{\bgamma}\supkt\bW_i) - g(\bgamma\subo\trans\bW_i)\}\|_{\psi_1}
   \le  M\sigma\submax^2 \|\hat{\bu}\supk\|_2\|\hat{\bgamma}\supk-\bgamma\subo\|_2/2.
$$
Applying the Bernstein's inequality, we obtain
$$
T_4 =  O_p\left((1-\rho)\|\hat{\bu}\supk\|_2\|\hat{\bgamma}\supk-\bgamma\subo\|_2\right).
$$
We achieve the stated rate with the tightness of $\|\hat{\bu}\supk\|_2$ from \eqref{eq:u0b-ukt}.

To show
$$
  T_5 =   \sqrt{n}\left(\bu\subo-\hat{\bu}\supk\right)\trans\left\{\dlkdag(\bbeta\subo; \bgamma\subo)+(1-\rho)\dlkimp(\bgamma\subo)\right\}
     =  O_p\left(\|\hat{\bu}\supk-\bu\subo\|_2\right),
$$
we rewrite the term as the empirical process with diminishing summands
\begin{align*}
   T_5 =  & - \sqrt{n} \frac{1}{N_k}\sum_{i\in \Jcal_k}(\hat{\bu}\supk-\bu\subo)\trans\bX_i \{g(\bbeta\subo\trans\bX_i)-g(\bgamma\subo\trans\bW_i) \} \\
   & - \sqrt{n} \frac{1}{n_k}\sum_{i\in \Ical_k} (\hat{\bu}\supk-\bu\subo)\trans \bX_i \{\rho \cdot g(\bbeta\subo\trans\bX_i)+ (1-\rho)\cdot g(\bgamma\subo\trans\bW_i) -Y_i \}.
\end{align*}
The summands have zero mean because
\begin{align*}
  & \E_{i \in \Jcal_k}\left[(\hat{\bu}\supk-\bu\subo)\trans\bX_i \{g(\bbeta\subo\trans\bX_i)-g(\bgamma\subo\trans\bW_i)\}\mid \Dscr_k^c \right] \\
= & (\hat{\bu}\supk-\bu\subo)\trans\E\left[\bX_i \{g(\bbeta\subo\trans\bX_i)-g(\bgamma\subo\trans\bW_i)\}\right] \\
= & 0, \\
  & \E_{i \in \Ical_k}\left[(\hat{\bu}\supk-\bu\subo)\trans \bX_i \{\rho \cdot g(\bbeta\subo\trans\bX_i)+ (1-\rho)\cdot g(\bgamma\subo\trans\bW_i) -Y_i \} \mid \Dscr_k^c \right] \\
= & (\hat{\bu}\supk-\bu\subo)\trans\left(\rho\E\left[\bX_i \{g(\bbeta\subo\trans\bX_i)-Y_i \}\right]
+ (1-\rho)\E\left[\bX_i \{g(\bgamma\subo\trans\bW_i) -Y_i \}\right] \right)\\
= & 0.
\end{align*}
Similar to the analysis of $T_3$, we show that each summand is sub-exponential
\begin{gather*}
  \left\|(\hat{\bu}\supk-\bu\subo)\trans\bX_i \{g(\bbeta\subo\trans\bX_i)-g(\bgamma\subo\trans\bW_i)\}\right\|_{\psi_1}
  \le \sqrt{2}\sigma\submax(\nu_1\vee\nu_2)\|\hat{\bu}\supk-\bu\subo\|_2 \\
    \left\|(\hat{\bu}\supk-\bu\subo)\trans \bX_i \{\rho \cdot g(\bbeta\subo\trans\bX_i)+ (1-\rho)\cdot g(\bgamma\subo\trans\bW_i) -Y_i \}\right\|_{\psi_1}
  \le \sqrt{2}\sigma\submax(\nu_1\vee\nu_2)\|\hat{\bu}\supk-\bu\subo\|_2
\end{gather*}
Applying the Bernstein's inequality, we obtain
$$
T_5 =  O_p\left(\|\hat{\bu}\supk-\bu\subo\|_2\left\{\sqrt{\rho(1-\rho)} + 1 \right\}\right)
= O_p\left(\|\hat{\bu}\supk-\bu\subo\|_2\right).
$$

\subsubsection*{Analysis for Terms $T'_1$-$T'_4$ in Part 2}

Conditionally on the out-of-fold data, the term $T'_1$ is the empirical average of i.i.d. mean zero random variables,
\begin{align*}
T_1 = &  \sum_{k=1}^{\Kn}\frac{n_k}{n}\Bigg( \frac{1}{n_k}\sum_{i\in\Ical_k} (\hat{\bu}\supkt  \bX_i)^2 \{(1-\rho)\cdot g(\hat{\bgamma}\supkt  \bW_i)+\rho \cdot g(\hat{\bbeta}\supkt \bX_i)-Y_i\}^2 \\
& \qquad -\E_{i\in\Ical_k}\left[(\hat{\bu}\supkt  \bX_i)^2 \{(1-\rho)\cdot g(\hat{\bgamma}\supkt  \bW_i)+\rho \cdot g(\hat{\bbeta}\supkt \bX_i)-Y_i\}^2\mid \Dscr_k^c\right]\Bigg).
\end{align*}
We bound the variance of each summand by the Cauchy-Schwartz inequality and Lemmas \ref{lemma:sub-moment}, \ref{lemma:sub-add}, \ref{lemma:sub-mb},
\begin{align}
   & \Var_{i\in\Ical_k}\left[(\hat{\bu}\supkt  \bX_i)^2 \{(1-\rho)\cdot g(\hat{\bgamma}\supkt  \bW_i)+\rho \cdot g(\hat{\bbeta}\supkt \bX_i)-Y_i\}^2\mid \Dscr_k^c\right] \notag \\
   \le & \E_{i\in\Ical_k}\left[(\hat{\bu}\supkt  \bX_i)^4 \{(1-\rho)\cdot g(\hat{\bgamma}\supkt  \bW_i)+\rho \cdot g(\hat{\bbeta}\supkt \bX_i)-Y_i\}^4\mid \Dscr_k^c\right] \notag \\
   \le &\sqrt{\E_{i\in\Ical_k}\left\{(\hat{\bu}\supkt  \bX_i)^8\mid \Lscr_k^c\right\}
   \E_{i\in\Ical_k}\left[\{(1-\rho)\cdot g(\hat{\bgamma}\supkt  \bW_i)+\rho \cdot g(\hat{\bbeta}\supkt \bX_i)-Y_i\}^8\mid \Dscr_k^c\right] } \notag \\
   \le & \sqrt{ 48 \|\hat{\bu}\supkt  \bX_i\|_{\psi_2}^8 48 \left(\rho\|g(\hat{\bbeta}\supkt \bX_i)-Y_i\|_{\psi_2}\vee (1-\rho)\|g(\hat{\bgamma}\supkt  \bX_i)-Y_i\|_{\psi_2}  \right)^8}
   \label{eq:Tpp1-var}
\end{align}
Under Assumption \ref{assume:Y-tail}, \ref{assume:link}, \ref{assume:X-subG},
we have
$$
\|\hat{\bu}\supkt  \bX_i\|_{\psi_2} \le (\|\bu\subo\|_2+\|\hat{\bu}\supk-\bu\subo\|_2)\sigma\submax/\sqrt{2}
= O_p\left(1+\|\hat{\bu}\supk-\bu\subo\|_2\right).
$$
We apply Lemma \ref{lemma:subg-g} to obtain
\begin{gather*}
\|g(\hat{\bbeta}\supkt \bX_i)-Y_i\|_{\psi_2} = O_p\left(1+\|\hat{\bbeta}\supk-\bbeta\subo\|_2\right), \\
\|g(\hat{\bgamma}\supkt \bW_i)-Y_i\|_{\psi_2} = O_p\left(1+ \|\hat{\bgamma}\supk-\bgamma\subo\|_2\right).
\end{gather*}
We have shown that the variance in \eqref{eq:Tpp1-var} is of order
$$
O_p\left(1+\|\hat{\bu}\supk-\bu\subo\|^4_2+\rho^4\|\hat{\bbeta}\supk-\bbeta\subo\|^4_2+(1-\rho)^4\|\hat{\bgamma}\supk-\bgamma\subo\|^4_2\right).
$$
Thus by the Tchebychev's inequality, we obtain
$$
T'_1 = O_p\left(\left\{1+\|\hat{\bu}\supk-\bu\subo\|_2+\rho\|\hat{\bbeta}\supk-\bbeta\subo\|_2+(1-\rho)\|\hat{\bgamma}\supk-\bgamma\subo\|_2\right\}/\sqrt{n}\right)
$$
Applying the consistency of $\hat{\bgamma}\supk$, $\hat{\bbeta}\supk$ and $\hat{\bu}\supk$ from
\eqref{assume:est-rate}
$$
T'_1 = O_p\left(n^{-1/2}\right) = o_p(1).
$$

To analyze $T'_2$, we consider the decomposition in which the estimators are replaced by the estimands one by one,
\begin{align*}
   &T'_2 \\
   = &\sum_{k=1}^{\Kn}\frac{n_k}{n}\bigg( \E_{i\in\Ical_k}\left[(\hat{\bu}\supkt  \bX_i)^2 \{(1-\rho)\cdot g(\hat{\bgamma}\supkt  \bW_i)+\rho \cdot g(\hat{\bbeta}\supkt \bX_i)-Y_i\}^2\mid \Dscr_k^c\right] \\
& \qquad -\E\left[(\bu\subo\trans \bX_i)^2 \{(1-\rho)\cdot g(\bgamma\subo\trans  \bW_i)+\rho \cdot g(\bbeta\subo\trans \bX_i)-Y_i\}^2\right]\bigg)\\
= &\sum_{k=1}^{\Kn}\frac{n_k}{n} \E_{i\in\Ical_k}\left[\{(\hat{\bu}\supk-\bu\subo)\trans  \bX_i\}\hat{\bu}\supkt  \bX_i \{(1-\rho)\cdot g(\hat{\bgamma}\supkt  \bW_i)+\rho \cdot g(\hat{\bbeta}\supkt \bX_i)-Y_i\}^2\mid \Dscr_k^c\right]\\
&+\sum_{k=1}^{\Kn}\frac{n_k}{n}\E_{i\in\Ical_k}\left[\{(\hat{\bu}\supk-\bu_0)\trans  \bX_i\}\bu\subo\trans  \bX_i \{(1-\rho)\cdot g(\hat{\bgamma}\supkt  \bW_i)+\rho \cdot g(\hat{\bbeta}\supkt \bX_i)-Y_i\}^2\mid \Dscr_k^c\right] \\
&+\sum_{k=1}^{\Kn}\frac{n_k}{n} \E_{i\in\Ical_k}\Big[(\bu\subo\trans  \bX_i)^2 \{(1-\rho)\cdot g(\hat{\bgamma}\supkt  \bW_i)+\rho \cdot g(\hat{\bbeta}\supkt \bX_i)-Y_i\} \\
& \qquad \times
[(1-\rho) \{g(\hat{\bgamma}\supkt  \bW_i)-g(\bgamma\subo\trans  \bW_i)\}+\rho \{ g(\hat{\bbeta}\supkt \bX_i)- g(\bbeta\subo\trans \bX_i)\}]\mid \Dscr_k^c\Big]\\
&+\sum_{k=1}^{\Kn}\frac{n_k}{n}\E_{i\in\Ical_k}\Big[(\bu\subo\trans  \bX_i)^2 \{(1-\rho)\cdot g(\bgamma\subo\trans   \bW_i)+\rho \cdot g(\bbeta\subo\trans \bX_i)-Y_i\} \\
& \qquad \times
[(1-\rho) \{g(\hat{\bgamma}\supkt  \bW_i)-g(\bgamma\subo\trans  \bW_i)\}+\rho \{ g(\hat{\bbeta}\supkt \bX_i)- g(\bbeta\subo\trans \bX_i)\}]\mid \Dscr_k^c\Big].
\end{align*}
Following the same calculation as in \eqref{eq:Tpp1-var},
we can bound the expectations
$$
T'_2 = O_p\left(\|\hat{\bu}\supk-\bu\subo\|_2+\rho\|\hat{\bbeta}\supk-\bbeta\subo\|_2+(1-\rho)\|\hat{\bgamma}\supk-\bgamma\subo\|_2\right)
$$
Applying the consistency of $\hat{\bgamma}$, $\hat{\bbeta}$ and $\hat{\bu}$ from Lemma \ref{lemma:est-gamma}, Theorem \ref{thm:consistency} and Part 1 in the proof of Theorem \ref{thm:inference},
we have established
$$
T'_2 = o_p(1).
$$

Repeating the analyses for $T'_1$ and $T'_2$,
we can show
\begin{gather*}
T'_3 = O_p\left(\rho\sqrt{(1-\rho)/N}\left\{1+\|\hat{\bu}\supk-\bu\subo\|_2+\|\hat{\bbeta}\supk-\bbeta\subo\|_2+\|\hat{\bgamma}\supk-\bgamma\subo\|_2\right\}\right)
= o_p(1), \\
T'_4 = O_p\left(\rho(1-\rho)\left\{\|\hat{\bu}\supk-\bu\subo\|_2+\|\hat{\bbeta}\supk-\bbeta\subo\|_2+\|\hat{\bgamma}\supk-\bgamma\subo\|_2\right\}\right)
= o_p(1)
\end{gather*}

\section{Additional Technical Details}\label{asection:detail}

\subsection{Definitions}

We adopt the following definition of sub-Gaussian and sub-exponential random variables.
\begin{definition}[Sub-Gaussian and Sub-Exponential Random Variables]\label{def:subge}
The sub-Gaussian parameter
for a random variable $V$ is defined as
\begin{equation*}
  \|V\|_{\psi_2} = \inf\left\{\sigma >0: \E(e^{V^2/\sigma^2}) \le 2 \right\}.
\end{equation*}
The random variable $V$ is sub-Gaussian if $\|V\|_{\psi_2}$ is finite.
The sub-Gaussian parameter
for a random vector $\bU$ is defined as
\begin{equation*}
  \|\bU\|_{\psi_2} = \sup_{\|\bv\|_2 = 1} \|\bv\trans\bU\|_{\psi_2}.
\end{equation*}
The sub-Gaussian parameter
for a random variable $V$ is defined as
\begin{equation*}
  \|V\|_{\psi_1} = \inf\left\{\nu >0: \E(e^{|V|/\nu}) \le 2 \right\}.
\end{equation*}
The random variable $V$ is sub-exponential if $\|V\|_{\psi_1}$ is finite.
The more general Orlicz norm for $\alpha \in (0,1)$ is defined as
\begin{equation*}
  \|V\|_{\psi_\alpha} = \inf\left\{\nu >0: \E\left[e^{(|V|/\nu)^\alpha}\right] \le 2 \right\}.
\end{equation*}
\end{definition}

Mimicking the (minimal) Restricted Eigenvalue condition on the minimal
eigenvalue of matrix over a cone \citep{BickelRT09}, we define the maximal Restricted Eigenvalue
in Definition \ref{def:maxRE}.
\begin{definition}[Maximal Restricted Eigenvalue]\label{def:maxRE}
For a cone-set of the indices set $\Ocal \subset \{1,\dots,p\}$
\begin{equation}
\Ccal_{\smbgamma} (\xi,\Ocal):=
\left\{v \in \R^{p+q+1}: \|v_{\Ocal^c}\|_1 \le \xi \|v_{\Ocal}\|_1  \right\},
\end{equation}
we define the maximal Restricted Eigenvalue of a matrix $\Sigma$ as
\begin{equation}
\mathrm{RE}\submax(\xi,\Ocal; \Sigma ) = \sup_{v \in \Ccal_{\smbgamma}(\xi,\Ocal)\setminus \{0\}}
  \frac{\sqrt{v\trans \Sigma v}}
  {\|v\|_2}.
\end{equation}
\end{definition}

\subsection{Statements of Existing Results}

The properties in Lemmas \ref{lemma:subge} and \ref{lemma:subg-var} are covered in \cite{vershynin2018HDP} Chapter 2 and 4.
\begin{lemma}[Properties of sub-Gaussian and sub-exponential random variables]\ \label{lemma:subge}
\begin{enumerate}[label = \alph*) , ref =\ref{lemma:subge}-\alph* ]
  \item \label{lemma:sub-tail} Tail-probability:
  \begin{gather*}
  \P(|V| \ge x)  \le 2 e^{-x/\|V\|_{\psi_1}}, \\
  \P(|V| \ge x)  \le 2 e^{-x^2/\|V\|_{\psi_2}^2};
  \end{gather*}

  \item\label{lemma:sub-moment} Moments: $\E(|V|^r) \le \min\{\csubge{1}\|V\|_{\psi_1}^r,\csubge{2}\|V\|_{\psi_2}^r\}$
  with $\csubge{1} = r!2$ and $\csubge{2} = \Gamma(r/2)r$,
  and $\E(|V|) \le \sqrt{\pi}\|V\|_{\psi_2}$;
  \item\label{lemma:subg-sube} Hierarchy: $\|V\|_{\psi_1} \le \|V\|_{\psi_2}$;
  \item\label{lemma:sub-add} Arbitrary addition: $\left\|\sum_{i=1}^{m}V_i\right\|_{\psi_2} \le m\max_{i=1,\dots,m}\|V_i\|_{\psi_2}$
  and $\left\|\sum_{i=1}^{m}V_i\right\|_{\psi_1} \le m\max_{i=1,\dots,m}\|V_i\|_{\psi_1}$;
  \item\label{lemma:sub-mb} Multiplication with bounded random variable: $\|V_1V_2\|_{\psi_2} \le \|V_1\|_{\psi_2} K$, $\|V_1V_2\|_{\psi_1} \le \|V_1\|_{\psi_1} K$ for $|V_2| \le K$ almost surely;
  \item\label{lemma:sub-mg} Multiplication between sub-Gaussian random variables:
  $\|V_1V_2\|_{\psi_1} \le \|V_1\|_{\psi_2} \|V_2\|_{\psi_2}$,
  in particular, $\|V_1\|_{\psi_1} \le \|V_1\|_{\psi_2}/\sqrt{\log(2)}$;
  \item\label{lemma:sub-hoef} Hoeffding's inequality: $V_1, \dots, V_m$ are independent mean zero sub-Gaussian random
  variables. For $t>0$,
  $$
  \P\left(\left|\sum_{i=1}^{m}V_i\right| \ge t \right) \le 4 \exp\left(-\frac{ t^2}{\csubge{3}\sum_{i=1}^{m}\|V_i\|_{\psi_2}^2}\right), \;
  \csubge{3} =8.
  $$
  \item\label{lemma:sub-bern} Bernstein's inequality: $V_1, \dots, V_m$ are independent mean zero sub-exponential random
  variables. For $t>0$,
  $\csubge{4} = 16$ and $\csubge{5} = 4$
  $$
  \P\left(\left|\sum_{i=1}^{m}V_i\right| \ge t \right) \le 2 \exp\left[- \min\left\{t^2\left(\csubge{4}\sum_{i=1}^{m}\|V_i\|_{\psi_1}^2\right)^{-1}, t \left(\csubge{5}\max_{i=1,\dots,m}\|V_i\|_{\psi_1}\right)^{-1}\right\}\right].
  $$
\end{enumerate}
\end{lemma}
\begin{lemma}\label{lemma:subg-var}
Let $V_1,\dots,V_m$ be i.i.d sub-Gaussian vectors in $\R^p$ such that
$$
\|v\trans V\|_{\psi_2}^2 \le K^2 \E\{(v\trans V)^2\}
$$
for some $1 \le K < \infty$.
Then,
\begin{align*}
\left\|\frac{1}{m}\sum_{i=1}^{m}V_iV_i\trans - \E(VV\trans) \right\|_2
= O_p\left(p/m + \sqrt{p/m}\right).
\end{align*}
\end{lemma}

From \cite{NRWY2010TR} and \cite{HuangZhang12} among other literatures,
we have the following results concerning the LASSO under the generalized
linear models.
\begin{lemma}\label{lemma:rsc}
Under Assumptions \ref{assume:link}, \ref{assume:X-subG} and \ref{assume:sigmin},
\begin{align*}
\P\bigg(& \implik(\bgamma\subo+\bgD) - \implik(\bgamma\subo)
  - \bgD^\top\impscore(\bgamma\subo)  \\
  & \ge \crsc{1} \|\bgD\|_2\{\|\bgD\|_2-\crsc{2}\sqrt{\log(p+q)/n}\|\bgD\|_1\},  \forall \|\bgD\|_2 \le 1 \bigg) \ge 1 - \crsc{3}e^{-\crsc{4}n} ; \\
  \P\bigg(& \ell\subPL(\bbeta\subo+\bgD ) - \ell\subPL(\bbeta\subo)
  - \bgD^\top\dell\subPL(\bbeta\subo)  \\
  & \ge \crsc{1} \|\bgD\|_2\{\|\bgD\|_2-\crsc{2}\sqrt{\log(p)/N}\|\bgD\|_1\},  \forall \|\bgD\|_2 \le 1 \bigg) \ge 1 - \crsc{3}e^{-\crsc{4}N}.
\end{align*}
The negative log-likelihoods are defined in \eqref{def:gamma-dev} and \eqref{def:beta-dev},
and their gradients defined in \eqref{def:dot-grad}.
 See Definition \ref{def:cond-exp} for the definition of conditional expectation notation.
The constants are all absolute.
\end{lemma}
  The two inequalities in Lemma \ref{lemma:rsc} are direct application of \cite{NRWY2010TR} Proposition 2 page 22.
We can construct an auxiliary loss function to prove the following lemma.
\begin{lemma}\label{lemma:u-rsc}
Under Assumptions \ref{assume:link}, \ref{assume:X-subG} and \ref{assume:sigmin},
\begin{align*}
&\P\bigg(\frac{1}{N_{k'}}\sum_{i\in \Ical_{k'}\cup\Jcal_{k'}} g'\left(\hat{\bbeta}\supkkt\bX_i\right) (\bgD\trans\bX_i)^2 \\
&\qquad \ge 2\crsc{1}^* \|\bgD\|_2^2-\crsc{1}^*\crsc{2}^*\sqrt{\log(p)/N}\|\bgD\|_2\|\bgD\|_1,  \forall \|\bgD\|_2 \le 1\bigg) \\
\ge & \P\left(\left\|\hat{\bbeta}\supk - \bbeta\subo\right\|_2 \le \frac{\sigma\submin^2}{2\sigma\submax^3}\right) - \crsc{3}^*e^{-\crsc{4}^*N}.
\end{align*}
The constants are all absolute.
\end{lemma}

\begin{proof}[Proof of Lemma \ref{lemma:u-rsc}]
First, we show $\sqrt{g'\left(\hat{\bbeta}\supkkt\bX_i\right)} \bX_i$ is a sub-Gaussian random vector whose second moment has all eigenvalues
bounded away from infinity and zero.
Under Assumptions \ref{assume:link} and \ref{assume:X-subG}, we may apply Lemma \ref{lemma:sub-mb},
$$
\left\|\bv\trans\sqrt{g'\left(\hat{\bbeta}\supkkt\bX_i\right)} \bX_i\right\|_{\psi_2} \le \sqrt{M} \|\bv\trans\bX_i\|_{\psi_2}
\le \sqrt{M}\sigma\submax \|\bv\|_2/\sqrt{2}.
$$
Thus, $\sqrt{g'\left(\hat{\bbeta}\supkkt\bX_i\right)} \bX_i$ is a sub-Gaussian random vector.
Under Assumptions \ref{assume:link} and \ref{assume:X-subG}, we can bound the maximal eigenvalue of its second moment,
$$
\bv\trans \E_{i \in \Ical_{k'}\cup\Jcal_{k'}}\left\{g'\left(\hat{\bbeta}\supkkt\bX_i\right)\bX_i\bX_i\trans\mid \Dscr_{k'}^c\right\}\bv
\le M \E\{(\bv\trans\bX_i)^2\} \le M \|\bv\|_2^2 \sigma\submax^2.
$$
We derive the lower bound for the  minimal eigenvalue of its second moment from Assumptions \ref{assume:link}, \ref{assume:X-subG}, \ref{assume:X-subG},
\ref{assume:sigmin-X},
 the Cauchy-Schwartz inequality and Lemma \ref{lemma:sub-moment},
\begin{align*}
& \bv\trans \E_{i \in \Ical_{k'}\cup\Jcal_{k'}}\left\{g'\left(\hat{\bbeta}\supkkt\bX_i\right)\bX_i\bX_i\trans\mid \Dscr_{k'}^c\right\}\bv \\
\ge & \bv\trans \E_{i \in \Ical_{k'}\cup\Jcal_{k'}}\{g'(\bbeta\subo\trans\bX_i)\bX_i\bX_i\trans\mid \Dscr_{k'}^c\}\bv \\
&   -\E_{i \in \Ical_{k'}\cup\Jcal_{k'}}\left[(\bv\trans\bX_i)^2\left\{g(\bbeta\subo^\top \bX_i) - g\left(\hat{\bbeta}\supkkt \bX_i\right)\right\}\mid \Dscr_{k'}^c\right] \\
\ge &  \|\bv\|_2^2\sigma\submin^2-   M \E_{i \in \Ical_{k'}\cup\Jcal_{k'}}\left[\left|(\bv\trans\bX_i)^2\left\{\left(\bbeta\subo-\hat{\bbeta}\supkk\right)\trans \bX_i\right\}\right|\mid \Dscr_{k'}^c\right] \\
\ge & \|\bv\|_2^2\sigma\submin^2- M \sqrt{\E\{(\bv\trans\bX_i)^4\}\E_{i \in \Ical_{k'}\cup\Jcal_{k'}}\left[\left\{\left(\bbeta\subo-\hat{\bbeta}\supkk\right)\trans \bX_i\right\}^2\mid \Dscr_{k'}^c\right]} \\
\ge & \|\bv\|_2^2\left(\sigma\submin^2-   M\sigma\submax^3 \left\|\hat{\bbeta}\supkk - \bbeta\subo\right\|_2\right).
\end{align*}
Whenever $\left\|\hat{\bbeta}\supk - \bbeta\subo\right\|_2 \le \frac{\sigma\submin^2}{2\sigma\submax^3}$,
we have
$$
\bv\trans \E_{i \in \Ical_{k'}\cup\Jcal_{k'}}\left\{g'\left(\hat{\bbeta}\supkkt\bX_i\right)\bX_i\bX_i\trans\mid \Dscr_{k'}^c\right\}\bv
\ge \|\bv\|_2^2 \sigma\submin^2/2.
$$

Second, we construct an auxiliary least square loss to apply \cite{NRWY2010TR}.
Let $\varepsilon_i$ be independent standard normal random variables.
Construct the loss function
$$
\Lcal\supkk(\bv) = \frac{1}{N_{k'}}\sum_{i\in \Ical_{k'}\cup\Jcal_{k'}} \left\{ \varepsilon_i + (\bv\subo - \bv)\trans \sqrt{g'\left(\hat{\bbeta}\supkkt\bX_i\right)} \bX_i   \right\}^2.
$$
By the design, we have
$$
\Lcal\supkk(\bv\subo+\bgD) - \Lcal\supkk (\bv\subo)
  - \bgD^\top\frac{\partial}{\partial \bv}\Lcal\supkk(\bv\subo+\bgD)
  = \frac{1}{N_{k'}}\sum_{i\in \Ical_{k'}\cup\Jcal_{k'}} g'\left(\hat{\bbeta}\supkkt\bX_i\right) (\bgD\trans\bX_i)^2.
$$
We apply Proposition 2 in \cite{NRWY2010TR} for $\Lcal\supkk(\bv)$ conditionally on out-of-fold data $\Dscr_{k'}^c$
and the event $\left\{\left\|\hat{\bbeta}\supkk - \bbeta\subo\right\|_2 \le \frac{\sigma\submin^2}{2\sigma\submax^3}\right\}$
to finish the proof.

\end{proof}

\begin{lemma}\label{lemma:est-gamma}
For a constant $\quantdlginf(n,p,q,\epsp) \asymp \sqrt{s_{\smbgamma} \log(p+q)/n}$,
the event
$$
\eventdlginf = \left\{\|\impscore(\bgamma\subo)\|_\infty = \left\|\frac{1}{n}\sum_{i=1}^{n}\bW_i\{g(\bgamma\subo\trans \bW_i) - Y_i\}\right\|_\infty \le \quantdlginf(n,p,q,\epsp)
\right\}
$$
occur with probability greater than $1-\epsp$
under Assumptions \ref{assume:Y-tail} and \ref{assume:X-subG}.
Setting $\lamg = 2$,
we have on event $\eventdlginf$ that
$$
\hat{\bgamma}-\bgamma\subo \in \Ccal_{\smbgamma}(3,\mathrm{supp}(\bgamma\subo))
=
\left\{v \in \R^{p+q+1}: \|v_{\Ocal_{\scriptscriptstyle \bgamma}^{\scriptscriptstyle \sf c}}\|_1 \le 3 \|v_{\Ocal_{\scriptscriptstyle \bgamma}}\|_1  \right\},
$$
where $\Ocal_{\scriptscriptstyle \bgamma} = \{j: \gamma_j \neq 0\}$
is the indices set for nonzero coefficient in $\bgamma\subo$.
Moreover, we have
$$
\|\hat{\bgamma}-\bgamma\subo\|_2 = O_p\left(\sqrt{s_{\smbgamma} \log(p+q)/n}\right).
$$
\end{lemma}
The concentration on the event $\eventdlginf$ is established by
the union bound of element wise concentration,
which is in turn obtained by the Bernstein inequality for sub-exponential random variables
(Lemma \ref{lemma:sub-bern}).
The rest of Lemma \ref{lemma:est-gamma} follows \cite{HuangZhang12} Lemma 1 page 5 (page 1843 of the issue)
and \cite{NRWY2010TR} Corollary 5 page 23.

\end{appendix}
\end{document}